\documentclass{article}
\usepackage{geometry}
\usepackage{amsmath, amsfonts, amssymb,bm}
\usepackage{amsthm,graphicx,float}
\graphicspath{{figures/}}
\usepackage[dvipsnames]{xcolor}
\usepackage{mathtools}

\usepackage[unicode=true,
 bookmarks=false,
 breaklinks=false,pdfborder={0 0 1},colorlinks=false] 
 {hyperref}
\hypersetup{
 colorlinks,citecolor=blue,filecolor=blue,linkcolor=blue,urlcolor=blue}

\allowdisplaybreaks

\usepackage{natbib}
\usepackage{float}
\usepackage{multirow}
\usepackage{footnote}
\usepackage{dsfont}
\usepackage{booktabs}

\usepackage[linesnumbered,ruled,vlined]{algorithm2e}
\usepackage{algorithmic}
 \usepackage[normalem]{ulem}
\usepackage{makecell}

\definecolor{yxc}{RGB}{255,0,0}
\definecolor{yc}{RGB}{190,0,255} 
\definecolor{ytw}{RGB}{255,0,127}
\definecolor{dacong}{RGB}{10,103,68}
\definecolor{cc}{RGB}{1,11,111}

\newcommand{\soft}[1]{{#1}_{\tau}}

\newcommand{\ex}[2]{\mathbb{E}_{#1}\left[#2\right]}

\newcommand{\V}{V} % V value function
\newcommand{\Q}{Q} % Q value function
\newcommand{\ssp}{\mathcal{S}} % state space
\newcommand{\asp}{\mathcal{A}} % action space
\newcommand{\real}{\mathbb{R}} % real numbers
\newcommand{\disct}{\gamma} % discount facter
\newcommand{\prn}[1]{\left({#1}\right)} % parentheses
\newcommand{\prnbig}[1]{\big({#1}\big)} % parentheses
\newcommand{\brk}[1]{\left[{#1}\right]} % bracket
\newcommand{\norm}[1]{\left\|{#1}\right\|} % norm
\newcommand{\abs}[1]{\left|{#1}\right|} % norm
\newcommand{\innprod}[1]{\big\langle{#1} \big\rangle} % inner product
\newcommand{\cS}{\mathcal{S}}
\newcommand{\cA}{\mathcal{A}}
\newcommand{\cB}{\mathcal{B}}

\newcommand{\KL}[2]{\mathsf{KL}\prnbig{{#1}\,\|\,{#2}}}
\newcommand{\KLbig}[2]{\mathsf{KL}\prnbig{{#1}\,\|\,{#2}}}

\newcommand{\simplexA}{\Delta(\mathcal{A})}
\newcommand{\simplexB}{\Delta(\mathcal{B})}

\newcommand{\best}[1]{#1_\tau^\star}

\newtheorem{theorem}{Theorem}
\newtheorem{lemma}{Lemma}
\newtheorem{definition}{Definition}
\newtheorem{prop}{Proposition}
\newtheorem{remark}{Remark}

\newtheorem{corollary}{Corollary}

\title{Fast Policy Extragradient Methods for Competitive Games \\ with Entropy Regularization}
\author{Shicong Cen\thanks{Department of Electrical and Computer Engineering, Carnegie Mellon University; email: \texttt{shicongc@andrew.cmu.edu}. } \\
	Carnegie Mellon University    \\
	\and
	Yuting Wei\thanks{Department of Statistics and Data Science, The Wharton School, University of Pennsylvania; email: \texttt{ytwei@wharton.upenn.edu}.}\\
	University of Pennsylvania \\
	\and
	Yuejie Chi\thanks{Department of Electrical and Computer Engineering, Carnegie Mellon University; email: \texttt{yuejiechi@cmu.edu}. }\\
	Carnegie Mellon University\\
}

\date{May 29, 2021; Revised \today}

\begin{document}
	\maketitle

\begin{abstract}

This paper investigates the problem of computing the equilibrium of competitive games, 
which is often modeled as a constrained saddle-point optimization problem with probability simplex constraints. Despite recent efforts in understanding the last-iterate convergence of extragradient methods in the unconstrained setting, the theoretical underpinnings of these methods in the constrained settings, especially those using multiplicative updates, remain highly inadequate, even when the objective function is bilinear. 
Motivated by the algorithmic role of entropy regularization in single-agent reinforcement learning and game theory, we develop provably efficient extragradient methods to find the quantal response equilibrium (QRE)---which are solutions to zero-sum two-player matrix games with entropy regularization---at a linear rate. 
The proposed algorithms can be implemented in a decentralized manner, where each player executes symmetric and multiplicative updates iteratively using its own payoff without observing the opponent's actions directly. In addition, by controlling the knob of entropy regularization, the proposed algorithms can locate an approximate Nash equilibrium of the unregularized matrix game at a sublinear rate without assuming the Nash equilibrium to be unique. 
Our methods also lead to efficient policy extragradient algorithms for solving (entropy-regularized) zero-sum Markov games at similar rates. All of our convergence rates are nearly dimension-free, which are independent of the size of the state and action spaces up to logarithm factors, 
highlighting the positive role of entropy regularization for accelerating convergence.

\end{abstract}

\noindent \textbf{Keywords:} zero-sum Markov game, matrix game, entropy regularization, global convergence, multiplicative updates, no-regret learning, extragradient methods 

\setcounter{tocdepth}{2}
	\tableofcontents

\section{Introduction} \label{sec:intro}

Finding the equilibrium of competitive games, which can be viewed as constrained saddle-point optimization problems with probability simplex constraints, lies at the heart of modern machine learning and decision making paradigms such as Generative Adversarial Networks (GANs) \citep{goodfellow2014generative}, competitive reinforcement learning (RL) \citep{littman1994markov}, game theory \citep{shapley1953stochastic}, adversarial training \citep{mertikopoulos2018cycles}, to name a few.   

In this paper, we study one of the most basic forms of competitive games, namely two-player zero-sum games, in both the matrix setting and the Markov setting. Our goal is to find the equilibrium policies of both players in an {\em independent} and {\em decentralized} manner \citep{daskalakis2020independent,wei2021last} with guaranteed {\em last-iterate convergence}. Namely, each player will execute symmetric and independent updates iteratively using its own payoff without observing the opponent's actions directly, and the final policies of the iterative process should be a close approximation to the equilibrium up to any prescribed precision. This kind of algorithms is more advantageous and versatile especially in federated environments, as it requires neither prior coordination between the players like two-timescale algorithms, nor a central controller to collect and disseminate the policies of all the players, which  are often unavailable due to privacy constraints.

\subsection{Last-iterate convergence in competitive games}
In recent years, there have been significant progresses in understanding the last-iterate convergence of simple iterative algorithms for {\em unconstrained} saddle-point optimization, where one is interested in bounding the sub-optimality of the last iterate of the algorithm, rather than say, the ergodic iterate --- which is the average of all the iterations --- that are commonly studied in the earlier literature. This shift of focus is motivated, for example, by the infeasibility of averaging large machine learning models in training GANs \citep{goodfellow2014generative}. While vanilla Gradient Descent\,/\,Ascent (GDA) may diverge or cycle even for bilinear matrix games \citep{daskalakis2018training}, quite remarkably, small modifications lead to guaranteed last-iterate convergence to the equilibrium in a non-asymptotic fashion. A flurry of algorithms is proposed, including Optimistic Gradient Descent Ascent (OGDA) \citep{rakhlin2013optimization,daskalakis2018limit,wei2020linear}, predictive updates \citep{yadav2017stabilizing}, implicit updates \citep{liang2019interaction}, and more. Several unified analyses of these algorithms have been carried out (see, e.g. \cite{mokhtari2020unified,liang2019interaction} and references therein), where these methods in principle all make clever extrapolation of the local curvature in a predictive manner to accelerate convergence. With slight abuse of terminology, in this paper, we refer to this ensemble of algorithms as extragradient methods \citep{korpelevich1976extragradient,tseng1995linear,mertikopoulos2018optimistic,harker1990finite}.

However, saddle-point optimization in the {\em constrained setting}, which includes competitive games as a special case, remains largely under-explored even for bilinear matrix games. While it is possible to reformulate constrained bilinear games to unconstrained ones using softmax parameterization of the probability simplex, this approach falls short of preserving the bilinear structure and convex-concave properties in the original problem, which are crucial to the convergence of gradient methods. Therefore, there is a strong necessity of understanding and  developing improved extragradient methods in the constrained setting, where existing analyses in the unconstrained setting do not generalize straightforwardly. 
\citet{daskalakis2018last} proposed the optimistic variant of the multiplicative weight updates (MWU) method \citep{arora2012multiplicative}---which is extremely natural and popular for optimizing over probability simplexes---called Optimistic Multiplicative Weight Updates (OMWU), and established the asymptotic last-iterate convergence of OMWU for matrix games. Very recently, \citet{wei2020linear} established non-asymptotic last-iterate convergences of OMWU. However, these last-iterate convergence results require the Nash equilibrium to be unique, and cannot be applied to problems with multiple Nash equilibria.

%such as softmax and neural network parameterization

\subsection{Our contributions}

Motivated by the algorithmic role of entropy regularization in single-agent RL \citep{neu2017unified,geist2019theory,cen2020fast} as well as its wide use in game theory to account for imperfect and noisy information \citep{mckelvey1995quantal,savas2019entropy}, we initiate the design and analysis of extragradient algorithms using {\em multiplicative updates} for finding the so-called quantal response equilibrium (QRE), which are solutions to competitive games with entropy regularization \citep{mckelvey1995quantal}. While finding QRE is of interest in its own right, by controlling the knob of entropy regularization, the QRE provides a close approximation to the Nash equilibrium (NE), and in turn acts as a smoothing scheme for finding the NE. Our contributions are summarized below, with the detailed problem formulations provided in Section~\ref{sec:formulation_matrix} for matrix games and Section~\ref{sec:formulation_markov} for Markov games, respectively.

\begin{itemize}
\item {\em Near dimension-free last-iterate convergence to QRE of entropy-regularized matrix games.} We propose two policy extragradient algorithms to solve entropy-regularized matrix games, namely the Predictive Update (PU) and OMWU methods, where both players execute symmetric and multiplicative updates without knowing the entire payoff matrix nor the opponent's actions. Encouragingly, we show that the last iterate of the proposed algorithms converges to the unique QRE at a linear rate that is almost independent of the size of the action spaces. Roughly speaking, to find an $\epsilon$-optimal QRE in terms of Kullback-Leibler (KL) divergence, it takes no more than
$$\widetilde{\mathcal{O}} \left( \frac{1}{\eta \tau} \log\left(\frac{1}{\epsilon}\right) \right)$$
iterations, where $\widetilde{\mathcal{O}}(\cdot)$ hides logarithmic dependencies. Here, $\tau$ is the regularization parameter, and $\eta $ is the learning rate of both players no larger than $\mathcal{O}(1/(\tau + \norm{A}_\infty))$, where $\|A\|_{\infty} =\max_{i,j} |A_{i,j}|$ is the $\ell_\infty$ norm of the payoff matrix $A$. Optimizing the learning rate, the iteration complexity is bounded by $\widetilde{\mathcal{O}} \left(\|A\|_{\infty}\tau^{-1}\log(1/\epsilon)\right)$.

%$\cA$ and $\cB$ are the action space of the players, 
 
\item {\em Last-iterate convergence to $\epsilon$-NE of unregularized matrix games without uniqueness assumption.} The QRE provides an accurate approximation to the NE by setting the entropy regularization $\tau$ sufficiently small, therefore our result directly translates to finding a NE with last-iterate convergence guarantee. Roughly speaking, to find an $\epsilon$-NE (measured in terms of the duality gap), it takes no more than
$$ \widetilde{\mathcal{O}} \left(\frac{\|A\|_{\infty} }{\epsilon}    \right)$$
iterations with optimized learning rates, again independent of the size of the action spaces up to logarithmic factors. Unlike prior literature \citep{daskalakis2018last,wei2020linear}, our last-iterate convergence guarantee does not require the NE to be unique.   

\item {\em No-regret learning of entropy-regularized OMWU.} We further establish that under a decaying learning rate, the proposed OMWU method achieves a {\em logarithmic} regret for the entropy-regularized matrix game --- on the order of $\mathcal{O}(\log T)$ --- even when only one player follows the algorithm against arbitrary plays of the opponent. By setting $\tau$ appropriately, this translates to a regret of $\mathcal{O}((T\log T)^{1/2})$ for the unregularized matrix game, therefore matching the regret in \cite{rakhlin2013optimization} without the need of mixing in an auxiliary uniform distribution for exploration.

\item {\em Extensions to two-player zero-sum Markov games.} By connecting value iteration with  matrix games, we propose a policy extragradient method for solving infinite-horizon discounted entropy-regularized zero-sum Markov games, which finds an $\epsilon$-optimal minimax soft Q-function --- in terms of $\ell_{\infty}$ error --- in at most $\widetilde{\mathcal{O}}\prn{\frac{1}{\tau(1-\gamma)^2}\log^2 \left( \frac{1}{\epsilon}\right) }$ iterations, where $\gamma \in (0,1)$ is the discount factor. By setting $\tau$ sufficiently small, the proposed method finds an $\epsilon$-approximate NE (measured in terms of the duality gap) of the unregularized Markov game within $\widetilde{\mathcal{O}}\left(\frac{1}{(1-\gamma)^3 \epsilon} \right)$ iterations, which is independent of the dimension of the state-action space up to logarithmic factors. 

\end{itemize}

 \begin{table}[!t]

\begin{center}
\resizebox{\textwidth}{!}{
\begin{tabular}{c|c|c|c|c}
\toprule  
\multirow{2}{*}{\makecell{Equilibrium\\ type}}  & \multirow{2}{*}{Method} & \multirow{2}{*}{\makecell{Convergence rate } \vphantom{$\frac{1^{7^{7}}}{1^{7^{7}}}$}} & \multirow{2}{*}{\makecell{Dimension-free } \vphantom{$\frac{1^{7^{7}}}{1^{7^{7}}}$}}  & \multirow{2}{*}{\makecell{Require \\unique NE}} \tabularnewline
 &  &  & &  \tabularnewline	
\toprule  
	\multirow{2}{*}{$\epsilon$-QRE} 
	& \multirowcell{2}{PU \& OMWU \\ \textbf{(this work)}} 	& \multirow{2}{*}{linear} & \multirow{2}{*}{{\color{blue} yes}}  & \multirow{2}{*}{n/a} \tabularnewline
	&  &  &   & \tabularnewline
\midrule 
    \multirow{6}{*}{$\epsilon$-NE} & \multirowcell{2}{OMWU\\ \citep{daskalakis2018last}} & \multirow{2}{*}{asymptotic } & \multirow{2}{*}{no} & \multirow{2}{*}{yes}   \tabularnewline
	 &  &  & &    \tabularnewline[1pt]
\cline{2-5}
     & \multirowcell{2}[-2pt]{OMWU\\ \citep{wei2020linear}} 	& \multirowcell{2}[-2pt]{sublinear + linear  } & \multirowcell{2}[-2pt]{no}  & \multirowcell{2}[-2pt]{yes} \tabularnewline
     &  &  &  & \tabularnewline[2pt]
\cline{2-5}
     & \multirowcell{2}[-2pt]{PU \& OMWU  \\ \textbf{(this work)}} 	& \multirowcell{2}[-2pt]{sublinear} & \multirowcell{2}[-2pt]{{\color{blue} yes}} & \multirowcell{2}[-2pt]{{\color{blue} no}} \tabularnewline
     &  &  & &  \tabularnewline
\bottomrule 
\end{tabular}
}
\end{center}

\caption{
Comparisons of last-iterate convergence of the proposed entropy-regularized PU and OMWU methods with prior results for finding $\epsilon$-QRE or $\epsilon$-NE of competitive matrix games. We note that the convergence rates of unregularized OMWU established in \citet{wei2020linear} are problem-dependent, and scale at least polynomially on the size of the action spaces. Desirable features in the last two columns are highlighted in {\color{blue}blue}. 
}\label{table:comparison}
\end{table}

  \begin{table}[!t]

\begin{center}
\resizebox{\textwidth}{!}{
\begin{tabular}{c|c|c|c|c|c}
\toprule  
Equilibrium   & \multirow{2}{*}{Method} &  Convergence \vphantom{$\frac{1^{7^{7}}}{1^{7^{7}}}$} & Dimension-free \vphantom{$\frac{1^{7^{7}}}{1^{7^{7}}}$}  &  Symmetric  &  Last-iterate  \tabularnewline
type &  & rate  & rate & updates  & convergence  \tabularnewline	
\toprule  
	\multirow{2}{*}{$\epsilon$-QRE} 
	& \multirowcell{2}{ \textbf{(this work)}} 	& \multirow{2}{*}{linear} & \multirow{2}{*}{{\color{blue} yes}}  & \multirow{2}{*}{{\color{blue}yes}}  &  \multirow{2}{*}{{\color{blue}yes}}  \tabularnewline
	&  &  &   &  & \tabularnewline
\midrule 
\multirow{10}{*}{$\epsilon$-NE}& \multirowcell{2}[-2pt]{ \citep{perolat2015approximate}} 	& \multirowcell{2}[-2pt]{sublinear  } & \multirowcell{2}[-2pt]{no}  & \multirowcell{2}[-2pt]{no} &  \multirow{2}{*}{{\color{blue}yes}} \tabularnewline
     &  &  &  & & \tabularnewline[2pt]
\cline{2-6}
     & \multirowcell{2}{ \citep{zhao2021provably}} & \multirow{2}{*}{sublinear } & \multirow{2}{*}{no} & \multirow{2}{*}{no}  &  \multirow{2}{*}{{\color{blue}yes}} \tabularnewline
	 &  &  & &   &  \tabularnewline[1pt]
\cline{2-6}
     & \multirowcell{2}[-2pt]{ \citep{daskalakis2020independent}} 	& \multirowcell{2}[-2pt]{sublinear  } & \multirowcell{2}[-2pt]{no}  & \multirowcell{2}[-2pt]{no}  & \multirowcell{2}[-2pt]{no}  \tabularnewline
     &  &  &  & & \tabularnewline[2pt]
\cline{2-6}
     & \multirowcell{2}[-2pt]{ \citep{wei2021last}} 	& \multirowcell{2}[-2pt]{sublinear  } & \multirowcell{2}[-2pt]{no}  & \multirowcell{2}[-2pt]{{\color{blue} yes}}& \multirowcell{2}[-2pt]{{\color{blue} yes}}  \tabularnewline
     &  &  &  & & \tabularnewline[2pt]
\cline{2-6}
     & \multirowcell{2}[-2pt]{ \textbf{(this work)}} 	& \multirowcell{2}[-2pt]{sublinear} & \multirowcell{2}[-2pt]{{\color{blue} yes}} & \multirowcell{2}[-2pt]{{\color{blue} yes}} & \multirowcell{2}[-2pt]{{\color{blue} yes}}  \tabularnewline
     &  &  & & & \tabularnewline
\bottomrule 
\end{tabular}
}
\end{center}

\caption{
Comparisons of the proposed policy extragradient method with competitive algorithms for finding $\epsilon$-NE of two-player zero-sum Markov games.  We note that the convergence rates in the prior arts all depend on various notions of concentrability coefficient and therefore not dimension-free.  Desirable features in the last three columns are highlighted in {\color{blue}blue}.  
}\label{table:comparison_markov}
\end{table}

To the best of our knowledge, our paper is the first one that develops policy extragradient algorithms for solving entropy-regularized competitive games with multiplicative updates and dimension-free linear last-iterate convergence, and demonstrates entropy regularization as a smoothing technique to find $\epsilon$-NE without the uniqueness assumption.  Table~\ref{table:comparison} and Table~\ref{table:comparison_markov} provide detailed comparisons of the proposed methods with prior arts for solving competitive games. Our results highlight the positive role of entropy regularization for accelerating convergence and safeguarding against imperfect payoff information in competitive games. 

%\paragraph{Reproducible research.} To reproduce the experiments in this paper, the codes can be found at: \yc{add a link to the codes}
%\begin{center}
%\url{https://github.com/..}.
%\end{center}

 %that introduced the quantal response equilibria. 
 
%In this paper, we consider entropy regularization which is introduced to model the information uncertainty, and accelerate convergence .

% Our results on matrix games in turn lead to efficient algorithms for solving infinite-horizon discounted Markov games. The convergence results are both dimensional-free, symmetric (namely, both players act in a symmetrical fashion), rational (namely, converging to the opponent?s best strategy under the stationary policy), underscoring the important role of entropy regularization for solving competitive games.

%  
%\begin{itemize}
% 
% 
%
%\item We establish an explicit convergence rate of game value \eqref{eq:EG_conv_f_last} \eqref{eq:EG_conv_f}, which enables us to characterize the iteration complexity of solving markov games by approximate value iteration with EG and OMWU as black box solvers. \dacong{Is there any similar result in the literature?}
% 
%\end{itemize}
%

	\subsection{Related works}

Our work lies at the intersection of saddle-point optimization, game theory, and reinforcement learning. In what follows, we discuss a few topics that are closely related to ours.

	% \begin{itemize}
	% 	\item Normal-form game = matrix game = $1$ state.
	% \end{itemize}
	\paragraph{Unregularized matrix game.} 
	\cite{freund1999adaptive} showed that many standard methods such as GDA and MWU have a converging average duality gap at the rate of $O(1/\sqrt{T})$, which is improved to $O(1/T)$ by considering optimistic variants of these methods, such as OGDA and OMWU \citep{rakhlin2013optimization,daskalakis2011near,syrgkanis2015fast}. However, the last-iterate convergence of these methods are less understood until recently \citep{daskalakis2018last,wei2020linear}. In particular, under the assumption that the NE is unique for the unregularized matrix game, \cite{daskalakis2018last} showed the asymptotic convergence of the last iterate of OMWU to the unique equilibrium, and \cite{wei2020linear} showed the last iterate of OMWU achieves a linear rate of convergence after an initial phase of sublinear convergence, however the rates therein can be highly pessimistic in terms of the problem dimension, while our rate for entropy-regularized OMWU is dimension-free up to logarithmic factors.
	 In terms of no-regret analysis, \cite{rakhlin2013optimization} established a no-regret learning rate of $O(\log T/{T}^{1/2})$ with an auxiliary mixing of a uniform distribution at each update, which is later improved to $O(1/{T}^{1/2})$ in \cite{kangarshahi2018let} with a slightly different algorithm.
	% \begin{itemize}
	% 	\item \cite{rakhlin2013optimization}: Regret analysis of some variant of OMWU.
	% 	\item \cite{daskalakis2018last}: Asymptotic convergence of OMWU.
	% 	\item \cite{wei2020linear}: linear convergence of OGDA \& OMWU. OMWU's result needs extra assumption (unique Nash equilibrium). OGDA's result can be applied to strongly-convex-strongly-concave problems, with a convergence rate of $(1 - O(\eta^2))^T$. (ours is $(1 - \eta \tau)^T$)
	% \end{itemize}

\paragraph{Saddle-point optimization.}  
Considerable progress has been made towards understanding OGDA and extragradient (EG) methods in the {unconstrained} convex-concave saddle-point optimization with general objective functions \citep{mokhtari2020unified,mokhtari2020convergence,nemirovski2004prox,liang2019interaction}. However, most works have focused on either average-iterate convergence (also known as ergodic convergence) \citep{nemirovski2004prox}, or the characterization of \textit{Euclidean update} rules \citep{mokhtari2020unified,mokhtari2020convergence,liang2019interaction}, where parameters are updated in an additive manner. These analyses do not generalize in a straightforward manner to \textit{non-Euclidean} updates.  As a result, the last-iterate convergence of {non-Euclidean} updates for saddle-point optimization still lacks theoretical understanding in general, and most works fall short of characterizing a finite-time convergence result. In particular, \cite{mertikopoulos2018optimistic} demonstrated the asymptotic last-iterate convergence of EG, and \citet{hsieh2019convergence} investigated similar questions for single-call EG algorithms. \cite{lei2021last} showed that OMWU converges to the equilibrium locally without an explicit rate. \cite{wei2020linear} showed that the last-iterate of OGDA converges linearly for strongly-convex strongly-concave constrained saddle-point optimization with an explicit rate.
%at a rate of $(1-O(\eta^2))^t$ 
% \cite{mokhtari2020unified} \cite{mokhtari2020convergence} \cite{nemirovski2004prox} constrained \cite{wei2020linear} \cite{lei2021last}
 
%  \cite{liang2019interaction}
 
%  \cite{mertikopoulos2018optimistic}

 \paragraph{Entropy regularization in RL and games.} In single-agent RL, the role of entropy regularization as an algorithmic mechanism to encourage exploration and accelerate convergence has been investigated extensively  \citep{neu2017unified,geist2019theory,mei2020global,cen2020fast,lan2021policy,zhan2021policy}. Turning to the game setting, entropy regularization is used to account for imperfect information in the seminal work of \cite{mckelvey1995quantal} that introduced the QRE, and a few representative works on entropy and more general regularizations in games include but are not limited to \cite{savas2019entropy,hofbauer2002global,mertikopoulos2016learning,cen2022independent}.
   
\paragraph{Zero-sum Markov games.} There have been a significant recent interest in developing provably efficient self-play algorithms for Markov games, including model-based algorithms \citep{perolat2015approximate,sidford2020solving,zhang2020model,li2022minimax,cui2021minimax}, value-based algorithms \citep{bai2020provable,xie2020learning,mao2022improving}, and policy-based algorithms \citep{daskalakis2020independent,wei2021last,zhao2021provably}. Our approach can be regarded as a policy-based algorithm to approximate value iteration, which can be implemented in a decentralized manner with symmetric and multiplicative updates from both players, and the iteration complexity is almost independent of the size of the state-action space. The iteration complexities in prior works \citep{perolat2015approximate,daskalakis2020independent,wei2021last,zhao2021provably} depend on various notions of concentrability coefficient and therefore can scale quite pessimistically with the problem dimension. In addition, while the last-iterate convergence guarantees in \citep{perolat2015approximate,daskalakis2020independent,zhao2021provably} are applicable to the duality gap, \cite{wei2021last} proves the last-iterate convergence in terms of the Euclidean distance to NE, together with an average convergence in terms of the duality gap. 

It is worth mentioning that the study of model-based and value-based algorithms typically focuses on the statistical issues in terms of sample complexity under the generative model or the online model of data collection; on the other end, the study of policy-based algorithms highlights the optimization issues by sidestepping the statistical issues using exact gradient evaluations, and later translating to sample complexity guarantees by leveraging model-based or value-based policy evaluation algorithms. Indeed, after the appearance of the initial version of this paper, \citet{chen2021sample} has built on our algorithm to develop a sample-efficient version of policy extragradient methods in the online setting using bandit feedback.

\subsection{Notation}
We denote by $\Delta(\cA)$ the probability simplex over the set $\cA$. We overload the functions such as $\log(\cdot)$ and $\exp(\cdot)$ to take vector inputs with the understanding that the function is applied in an entrywise manner. For instance, given any vector $z=[z_i]_{1\leq i\leq n}\in \mathbb{R}^n$, the notation $\exp(z)$ denotes $\exp(z) \coloneqq [\exp(z_i)]_{1\leq i\leq n}$; other functions are defined analogously. Given two probability distributions $\mu$ and $\mu'$ over $\cA$,  the KL divergence from $\mu'$ to $\mu$ is defined by $\KL{\mu}{\mu'} \coloneqq \sum_{a\in\cA}\mu(a) \log\frac{\mu(a)}{\mu'(a)}$.  Given a matrix $A$, $\norm{A}_\infty$ is used to denote entrywise maximum norm, namely, $\norm{A}_\infty=\max_{i,j}|A_{i,j}|$. 
The all-one vector is denoted as $\mathbf{1}$.

%For any vectors $z=[z_i]_{1\leq i\leq n}$ and $w=[w_i]_{1\leq i\leq n}$, the notation $z \geq w$ (resp.~$ z \leq  w$) means $z_i \geq w_i$ (resp.~$z_i\leq w_i$) for all $1\leq i\leq n$. 

\section{Zero-sum matrix games with entropy regularization} 
\label{sec:background}

In this section, we consider a two-player zero-sum game with bilinear objective and probability simplex constraints, 
and demonstrate the positive role of entropy regularization in solving this problem. 
Throughout this paper, let $\cA=\{1,\ldots,m\}$ and $\cB =\{1,\ldots, n\}$ be the action spaces of each player. The proofs for this section are collected in Appendix~\ref{sec:matrix-games-proof-main}.
%with $\Delta(\cA) \coloneqq \big\{ z\in \mathbb{R}^m \mid z\geq 0, 1^{\top} z=1 \big\}$ the $m$-dimensional probability simplex, and $\Delta(\cB)$ defined similarly. 

\subsection{Background and problem formulation}
\label{sec:formulation_matrix}

\paragraph{Zero-sum two-player matrix game.}

The focal point of this subsection is a constrained two-player zero-sum matrix game, 
which can be formulated as the following min-max problem (or saddle point optimization problem):
\begin{align}
    \max_{\mu \in \Delta(\cA)}\min_{\nu \in \Delta(\cB) } f(\mu, \nu) \coloneqq \mu^\top A \nu  ,
	\label{eq:max-min-problem-games}
\end{align}
where $A\in\mathbb{R}^{m\times n}$ denotes the payoff matrix, $\mu\in \Delta(\cA)$ and $\nu \in \Delta(\cB)$ stand for the mixed/randomized policies of each player, defined respectively as distributions over the probability simplex $\Delta(\cA)$ and $\Delta(\cB)$. It is well known since \citet{neumann1928theorie} that the max and min operators in \eqref{eq:max-min-problem-games} can be exchanged without affecting the solution. A pair of policies $(\mu^{\star},\nu^{\star})$ is said to be a {\em Nash equilibrium (NE)} of \eqref{eq:max-min-problem-games} if
\begin{align}
	f(\mu^{\star}, \nu) \geq f(\mu^{\star}, \nu^{\star }) \geq f(\mu, \nu^{\star }) \qquad \text{for all }  (\mu, \nu) \in \Delta(\cA) \times \Delta(\cB).  
	\label{eq:defn-Nash-equiv-matrix}
\end{align}
In words, the NE corresponds to when both players play their best-response strategies against their respective opponents.

\paragraph{Entropy-regularized zero-sum two-player matrix game.} 
There is no shortage of scenarios where the payoff matrix $A$ might not be known perfectly. 
In an attempt to accommodate imperfect knowledge of $A$, \citet{mckelvey1995quantal} proposed a seminal extension to the Nash equilibrium called the {\em quantal response equilibrium (QRE)} when the payoffs are perturbed by  Gumbel-distributed noise. Formally, this amounts to solving the following matrix game with entropy regularization \citep{mertikopoulos2016learning}:
\begin{equation}\label{eq:QRE}
    \max_{\mu \in \Delta(\cA)}\min_{\nu \in \Delta(\cB) } f_{\tau}(\mu, \nu) \coloneqq \mu^\top A \nu + \tau \mathcal{H}(\mu) - \tau \mathcal{H}(\nu) ,
\end{equation}
where  $\mathcal{H}(\pi) \coloneqq - \sum _i \pi_i \log(\pi_i)$ denotes the Shannon entropy of a distribution $\pi$, and $\tau\geq 0$ is the regularization parameter. 
As is well known, the optimal solution $(\best{\mu}, \best{\nu})$ to \eqref{eq:QRE}, dubbed as the QRE, is unique whenever $\tau>0$ (due to the presence of strong concavity/convexity), 
which satisfies the following fixed point equations:
\begin{align}
    \begin{cases}
	    \best{\mu}(a) = \frac{ \exp( [A \best{\nu}]_a/\tau  ) }{\sum_{a=1}^m \exp( [A \best{\nu}]_a/\tau  )}  \propto \exp( [A \best{\nu}]_a/\tau  ) , \qquad  & \text{for all } a \in \mathcal{A} ,  \\[0.2cm]
	    \best{\nu}(b) = \frac{ \exp( -[A^\top \best{\mu}]_b/\tau) } {\sum_{b=1}^n \exp( -[A^\top \best{\mu}]_b/\tau) } \propto \exp( -[A^\top \best{\mu}]_b/\tau) ,\qquad & \text{for all } b\in \mathcal{B} .\\      
    \end{cases}
	\label{eq:QRE-matrix}
\end{align}

\paragraph{Goal.} We aim to efficiently compute the QRE of the entropy-regularized matrix game in a decentralized manner, and investigate how an efficient solver of QRE can be leveraged to find a NE of the unregularized matrix game \eqref{eq:max-min-problem-games}.
Namely, we only assume access to ``first-order information'' as opposed to full knowledge of the payoff matrix $A$ or the actions of the opponent. The information received by each player is formally 
described in the following sampling oracle. 
\begin{definition}[Sampling oracle for matrix games] 
	\label{defn:sampling-oracle-matrix}
	For any policy pair $(\mu,\nu)$ and payoff matrix $A$, the sampling oracle returns the exact values of $\mu^{\top}A$ and $A\nu$. 
\end{definition}

\paragraph{Additional notation.} For notational convenience, we let $\zeta$ represent the concatenation of $\mu \in \real^{|\mathcal{A}|} $ and $\nu\in \real^{|\mathcal{B}|}$, namely,  
$\zeta = (\mu, \nu)$. 
The solution to \eqref{eq:QRE}, which is specified in \eqref{eq:QRE-matrix}, is denoted by $\best{\zeta} = (\best{\mu}, \best{\nu})$.
% $\zeta = (\mu, \nu)$. 
For any $\zeta = (\mu, \nu)$ and $\zeta' = (\mu', \nu')$, we shall often abuse the notation and let $$\KL{\zeta}{\zeta'}=\KL{\mu}{\mu'} + \KL{\nu}{\nu'}.$$
The duality gap of the entropy-regularized matrix game \eqref{eq:QRE} at $\zeta = (\mu,\nu)$ is defined as 
    \begin{align}
    \label{eqn:def-duality-gap}
        \mathsf{DualGap}_\tau(\zeta) & = \max_{\mu'\in \Delta(\mathcal{A})} f_\tau(\mu', \nu) - \min_{\nu'\in\Delta(\mathcal{B})}f_\tau(\mu, \nu')
   %      = \max_{\mu'\in \Delta(\mathcal{A}),\nu'\in\Delta(\mathcal{B})} \left\{ f_\tau(\mu', \nu) - f_\tau(\mu, \nu')\right\},
    \end{align}
    which is clearly nonnegative and $\mathsf{DualGap}_{\tau}(\zeta_{\tau}^{\star})=0$. Similarly, let the optimality gap of the entropy-regularized matrix game \eqref{eq:QRE} at $\zeta = (\mu,\nu)$ be 
    $\mathsf{OptGap}(\zeta) = 	\big| f_\tau( {\mu}, {\nu}) - f_\tau(\best{\mu}, \best{\nu}) \big| $.

    %Last, the optimality gap of \eqref{eq:QRE} at $\zeta = (\mu,\nu)$ is defined as

\subsection{Proposed extragradient methods: PU and OMWU}

% In this section, we shall establish the proof of Theorem~\ref{thm:EG-gurantees-last-iterate} and \ref{thm:EG-gurantees-Ave} together. 
% Towards this, it is useful to first consider the following updating rule that underlies both EG and OMWU

%By replacing $\mu$, $\nu$ with fixed policy $z_1 \in \simplexA$, $z_2 \in \simplexB$ separately,
 
To begin, assume we are given a pair of policies $z_1 \in \simplexA$, $z_2 \in \simplexB$ employed by each player respectively. If we proceed with fictitious play, i.e. player 1 (resp. player 2) aims to optimize its own policy by assuming the opponent's policy is fixed as $z_2$ (resp. $z_1$), the saddle-point optimization problem \eqref{eq:QRE}  is then decoupled into two independent min/max optimization problems:
\[
        \max_{\mu \in \Delta(\cA)} \; \mu^\top A z_2 + \tau \mathcal{H}(\mu) - \tau \mathcal{H}(z_2)  \qquad \mbox{and}\qquad
        \min_{\nu \in \Delta(\cB)} \; z_1^\top A \nu + \tau \mathcal{H}(z_1) - \tau \mathcal{H}(\nu),
\]
which are naturally solved via mirror descent\,/\,ascent with KL divergence. Specifically, one step of mirror descent\,/\,ascent takes the form
\begin{equation*}   
    \begin{cases}
        \mu^{(t+1)} = \arg\max_{\mu \in \Delta(\cA)} \; (A z_2 - \tau \log \mu^{(t)})^\top \mu - \frac{1}{\eta}\KL{\mu}{\mu^{(t)}}\\[0.2cm]
        \nu^{(t+1)} = \arg\min_{\nu \in \Delta(\cB)} \; (A^\top z_1 + \tau \log \nu^{(t)})^\top \nu + \frac{1}{\eta}\KL{\nu}{\nu^{(t)}}
    \end{cases},
    % \label{eq:mirror_update}
\end{equation*}
where $\eta$ is the learning rate, or equivalently
\begin{equation}  \label{eq:mirror_update}
    \begin{cases}
        \mu^{(t+1)}(a) \propto {\mu^{(t)}(a)}^{1-\eta\tau}\exp(\eta [A z_2]_a) , \qquad  & \text{for all } a \in \mathcal{A},\\[0.2cm]
        \nu^{(t+1)}(b) \propto {\nu^{(t)}(b)}^{1-\eta\tau}\exp(-\eta [A^\top z_1]_b) , \qquad  & \text{for all } b \in \mathcal{B}.\\
    \end{cases}
\end{equation}
The above update rule forms the basis of our algorithm design.

\paragraph{Motivation: a form of implicit updates with linear convergence.} To begin with, we select the policy pair $(z_1,z_2) = \zeta^{(t+1)}: = ( \mu^{(t+1)},  \nu^{(t+1)} )$ as the solution to the following equations, and call the conceptual update rule as the Implicit Update (IU) method:
\begin{equation}\label{eq:implicit_updates}
\mbox{Implicit Update:} \qquad
    \begin{cases}
        \mu^{(t+1)}(a) \propto {\mu^{(t)}(a)}^{1-\eta\tau}\exp(\eta [A \nu^{(t+1)}]_a), \qquad  & \text{for all } a \in \mathcal{A},\\[0.2cm]
        \nu^{(t+1)}(b) \propto {\nu^{(t)}(b)}^{1-\eta\tau}\exp(-\eta [A^\top \mu^{(t+1)}]_b), \qquad  & \text{for all } b \in \mathcal{B}. \\
    \end{cases} 
\end{equation} 
Though unrealistic --- since it uses the future updates and denies closed-form solutions --- it leads to a one-step convergence to the QRE when $\eta= 1/\tau$ (see the optimality condition in \eqref{eq:QRE-matrix}). Encouragingly, we have the following linear convergence guarantee of IU when adopting a general learning rate.
\begin{prop}[Linear convergence of IU] \label{prop:iu}
 Assume $0< \eta \leq 1/\tau$, then for all $t\geq 0$, the iterates $\zeta^{(t)}: = ( \mu^{(t)},  \nu^{(t)} )$ of the IU method in \eqref{eq:implicit_updates} satisfy
\begin{align*}
    \KL{\best{\zeta}}{\zeta^{(t)}} \le (1-\eta\tau)^t \KL{\best{\zeta}}{\zeta^{(0)}}.
\end{align*}
\end{prop}
In words, the IU method achieves an appealing linear rate of convergence that is independent of the problem dimension. Motivated by this observation, we seek to design algorithms where the policies $(z_1,z_2)$ employed in \eqref{eq:mirror_update} serve as good predictions of $(\mu^{(t+1)},\nu^{(t+1)})$, such that the resulting algorithms are both practical and retain the appealing convergence rate of IU.
% Guided by this intuition, the proof of Theorem~\ref{thm:EG-gurantees-last-iterate} and \ref{thm:EG-gurantees-Ave} are perturbed version of the above analysis. 

\paragraph{Proposed algorithms.}
We propose two extragradient algorithms for solving the entropy-regularized matrix game, namely the {\em Predictive Update (PU)} method and the {\em Optimistic Multiplicative Weights Update (OMWU)} method, where the latter is adapted from \cite{rakhlin2013optimization}.
Detailed procedures can be found in Algorithm~\ref{alg:ex_grad} and Algorithm~\ref{alg:omwu}, respectively. On a high level, both algorithms maintain two intertwined sequences $\{ (\mu^{(t)}, \nu^{(t)} ) \}_{t\geq 0}$ and $\{ (\bar{\mu}^{(t)}, \bar{\nu}^{(t)}) \}_{t\geq 0}$, and 
in each iteration $t=0,1,\ldots$, proceed in two steps:
\begin{itemize}
\item The midpoint $(\bar{\mu}^{(t+1)}, \bar{\nu}^{(t+1)})$ serves as a prediction of $(\mu^{(t+1)},\nu^{(t+1)})$ by running one step of mirror descent\,/\,ascent (cf.~\eqref{eq:mirror_update}) from either $(z_1,z_2) = (\mu^{(t)},\nu^{(t)})$ (for PU) or $(z_1,z_2) = (\bar{\mu}^{(t)}, \bar{\nu}^{(t)})$ (for OMWU). 

\item The update of $(\mu^{(t+1)}, \nu^{(t+1)})$ then mimics the implicit update \eqref{eq:implicit_updates} using the prediction $(\bar{\mu}^{(t+1)}, \bar{\nu}^{(t+1)})$ obtained above.
\end{itemize} 
%rely on first-order information pertaining to the midpoints (i.e., $A\bar{\nu}^{(t)}$ and $A^{\top}\bar{\mu}^{(t)}$).  

When the proposed algorithms converge, both $(\mu^{(t)}, \nu^{(t)})$ and $(\bar{\mu}^{(t)}, \bar{\nu}^{(t)})$ converge to the same point. 
The two players are completely symmetric and adopt the same learning rate, and require \emph{only} first-order information provided by the sampling oracle. While the two algorithms resemble each other in many aspects, a key difference lies in the query and use of the sampling oracle:  
in each iteration, OMWU makes a single call to the sampling oracle for gradient evaluation, 
while PU calls the sampling oracle twice. It is worth noting that, when $\tau=0$ (i.e., no entropy regularization is enforced), the OMWU method in Algorithm~\ref{alg:omwu} reduces to the method analyzed in \cite{rakhlin2013optimization,daskalakis2018last,wei2020linear} without entropy regularization.

\begin{remark}
It is worth highlighting that the proposed algorithms are {\em different} from the mirror prox algorithm \citep{nemirovski2004prox} or the optimistic mirror descent method \citep{mertikopoulos2018optimistic}, as the extragradient is only applied to the bilinear term but not the entropy regularization term. This seemingly small, but important, difference leads to a more concise closed-form update rule and a cleaner analysis, as shall be seen momentarily. 
\end{remark}

\begin{figure}[t]
\begin{minipage}{0.47\textwidth}

\begin{algorithm}[H]
    \DontPrintSemicolon
       \textbf{initialization:}  $\mu^{(0)}$, $\nu^{(0)}$.  \\
       \textbf{parameters:} learning rate $\eta_t$.
        % number of iterations $T$, 
      %\textbf{initialization:} $Q_0=0$. \\
   
       \For{$t=0,1,2,\cdots$}
        {
             Update $\bar\mu$ and $\bar\nu$ according to
            \[
                \begin{cases}
                    \bar{\mu}^{(t+1)}(a) \propto {\mu^{(t)}(a)}^{1-\eta_t\tau}\exp(\eta_t [A \nu^{(t)}]_a),\\
                    \bar{\nu}^{(t+1)}(b) \propto {\nu^{(t)}(b)}^{1-\eta_t\tau}\exp(-\eta_t [A^\top \mu^{(t)}]_b).
                \end{cases}
            \]
            
            Update $\mu$ and $\nu$ according to
            \[
                \begin{cases}
                    {\mu}^{(t+1)}(a) \propto {\mu^{(t)}(a)}^{1-\eta_t\tau}\exp(\eta_t [A \bar{\nu}^{(t+1)}]_a),\\
                    {\nu}^{(t+1)}(b) \propto {\nu^{(t)}(b)}^{1-\eta_t\tau}\exp(-\eta_t [A^\top \bar{\mu}^{(t+1)}]_b).
                \end{cases}
            \]
   
        }
        \caption{The PU method}
     \label{alg:ex_grad}

    \end{algorithm}

\end{minipage}
\begin{minipage}{0.02\textwidth}
\qquad
\end{minipage}
\begin{minipage}{0.47\textwidth}

\begin{algorithm}[H]
    \DontPrintSemicolon
       \textbf{initialization:} $\mu^{(0)}=\bar{\mu}^{(0)}$, $\nu^{(0)}=\bar{\nu}^{(0)} $. \\
        \textbf{parameters:} learning rate $\eta_t$.
        % number of iterations $T$, 
      %\textbf{initialization:} $Q_0=0$. \\
   
       \For{$t=0,1,2,\cdots$}
        {
Update $\bar\mu$ and $\bar\nu$ according to
            \[
                \begin{cases}
                    \bar{\mu}^{(t+1)}(a) \propto {\mu^{(t)}(a)}^{1-\eta_t\tau}\exp(\eta_t [A \bar{\nu}^{(t)}]_a),\\
                    \bar{\nu}^{(t+1)}(b) \propto {\nu^{(t)}(b)}^{1-\eta_t\tau}\exp(-\eta_t [A^\top \bar{\mu}^{(t)}]_b).
                \end{cases}
            \]
            
            Update $\mu$ and $\nu$ according to
            \[
                \begin{cases}
                    {\mu}^{(t+1)}(a) \propto {\mu^{(t)}(a)}^{1-\eta_t\tau}\exp(\eta_t [A \bar{\nu}^{(t+1)}]_a),\\
                    {\nu}^{(t+1)}(b) \propto {\nu^{(t)}(b)}^{1-\eta_t\tau}\exp(-\eta_t [A^\top \bar{\mu}^{(t+1)}]_b).
                \end{cases}
            \]
   
        }
        \caption{The OMWU method}
     \label{alg:omwu}
    \end{algorithm}

\end{minipage}

\end{figure}

\subsection{Last-iterate linear convergence guarantees}
\label{sec:theory-matrix-game}

We are now positioned to present our main theorem concerning the last-iterate convergence of PU and OMWU for solving \eqref{eq:QRE}. 
Its proof can be found in Section~\ref{Sec:pf-matrix-game}. 
\begin{theorem}[Last-iterate convergence of PU and OMWU]
\label{thm:EG-gurantees-last-iterate}
	Suppose that the learning rates $\eta_t = \eta = \eta_{\mathsf{PU}} $ of PU in Algorithm~\ref{alg:ex_grad} and $\eta_t = \eta = \eta_{\mathsf{OMWU}} $ of OMWU in Algorithm~\ref{alg:omwu} satisfy 
	\begin{align}\label{eq:learning_rates_req}
	0 < \eta_{\mathsf{PU}}    \le \frac{1}{\tau + 2\norm{A}_\infty}, \quad \mbox{and}\quad
	0< \eta_{\mathsf{OMWU}}   \le \min \left\{\frac{1}{2\tau + 2\norm{A}_\infty  },\, \frac{1}{4\norm{A}_\infty} \right\}.
	\end{align}
	Then for any $t\geq 0$, the iterates $\zeta^{(t)}=(\mu^{(t)},\nu^{(t)})$ and $\bar{\zeta}^{(t)}=(\bar{\mu}^{(t)},\bar{\nu}^{(t)})$  of both PU and OMWU achieve
\begin{subequations}
\label{eq:EG-conv-all}
\begin{itemize}
\item {\bf Linear convergence of policies in KL divergence and entrywise log-ratios:}
\begin{align}
	\max\left\{ \KLbig{\best{\zeta}}{\zeta^{(t)}},\; \tfrac{1}{2} \KL{\best{\zeta}}{\bar{\zeta}^{(t+1)}}  \right\} & \le (1-\eta\tau)^t \KLbig{\best{\zeta}}{\zeta^{(0)}},
    \label{eq:EG_conv_KL} \\
%\end{align}
%and
%\begin{align}
    \norm{\log \frac{\zeta^{(t)}}{\best{\zeta}}  }_\infty \le 2(1-\eta\tau)^{t}\norm{\log \frac{\zeta^{(0)}}{ \best{\zeta}}}_\infty & + \frac{8\norm{A}_\infty}{\tau} (1-\eta\tau)^{t/2} {\KLbig{\best{\zeta}}{\zeta^{(0)}}}^{1/2}.
    \label{eq:EG_conv_log}
\end{align}
\item {\bf Linear convergence of values in optimality and duality gaps:}
\begin{align}
\mathsf{OptGap}_\tau(\bar{\zeta}^{(t)})   &  \le 3\eta^{-1}(1-\eta\tau)^{t}  \KLbig{\best{\zeta}}{\zeta^{(0)}},    \label{eq:EG_conv_f_last} \\
   %\end{align}
%and
%\begin{align}
 \mathsf{DualGap}_\tau(\bar{\zeta}^{(t)})  &\le \prn{\eta^{-1}+2\tau^{-1}{\norm{A}_\infty^2}}(1-\eta\tau)^{t-1}\KLbig{\best{\zeta}}{\zeta^{(0)}}.
    \label{eq:EG_conv_gap_last} 
\end{align}
\end{itemize}
%
%
% In addition, if $0 < \eta =\eta_{\mathsf{PU}} \le \tfrac{1}{\tau + 2\norm{A}_\infty}$, then $\zeta^{(t)}$ of PU in Algorithm~\ref{alg:ex_grad} \dacong{and OMWU} further satisfies

\end{subequations}
\end{theorem}

\begin{remark} To further understand the term $\KLbig{\best{\zeta}}{\zeta^{(0)}}$ in \eqref{eq:EG-conv-all},
    setting $\mu^{(0)}$ and $\nu^{(0)}$ to be uniform policies leads to a universal bound
    $$ \KLbig{\best{\zeta}}{\zeta^{(0)}} = \log|\mathcal{A}| + \log |\mathcal{B}| - \mathcal{H}(\best{\mu}) - \mathcal{H}(\best{\nu}) 
    \le \log|\mathcal{A}| + \log |\mathcal{B}|$$ 
    regardless of $\best{\zeta}= (\best{\mu},\best{\nu})$. 
\end{remark}

\begin{remark}
    Similar results continue to hold even when the two players use different regularization parameters $\tau_\mu, \tau_\nu>0$ in \eqref{eq:QRE}, as long as the regularization parameter $\tau$ is replaced by  $\max\{\tau_\mu, \tau_\nu\}$ in the upper bounds of the learning rate,  and the contraction parameter is replaced by $1 - \min\{\tau_\mu, \tau_\nu\}\eta$.
\end{remark}
%
%\begin{remark}
%Inequalities \eqref{eq:EG_conv_KL} and \eqref{eq:EG_conv_f_last} for PU hold for a slightly wider range of stepsize namely, $0 < \eta_{\mathsf{PU}} \le \frac{1}{\tau + \norm{A}_\infty}$ as can be seen in the proof. 
%% \ytw{it is true for \eqref{eq:EG_conv_gap_last} right?}
%\end{remark}
%n particular, we derive performance guarantees both in terms of the KL divergence and the value function with universal constant stepsizes. 

Theorem~\ref{thm:EG-gurantees-last-iterate} characterizes the convergence of the {\em last-iterates} $\zeta^{(t)}$ and $\bar{\zeta}^{(t)}$ of PU and OMWU as long as the learning rate lies within the specified ranges. While PU doubles the number of calls to the sampling oracle, it also allows roughly as large as twice the learning rate compared with OMWU (cf. \eqref{eq:learning_rates_req}). Compared with the vast literature analyzing the average-iterate performance of variants of extragradient methods, our results contribute towards characterizing the last-iterate convergence of multiplicative update methods in the presence of entropy regularization and simplex constraints, which to the best of our knowledge, are the first of its kind. 
Several remarks are in order. 
\begin{itemize}
 %   \item {\bf Last-iterate convergence.} 

\item {\bf Linear convergence to QRE.}
To achieve an $\epsilon$-accurate estimate of the QRE in terms of the KL divergence, the bound \eqref{eq:EG_conv_KL} tells that it is sufficient to take 
$$\frac{1}{\eta \tau} 
\log\left(\frac{\log|\mathcal{A}| + \log|\mathcal{B}|}{\epsilon}\right)$$ 
 iterations using either PU or OMWU. 
Notably, this iteration complexity does not depend on any hidden constants and only depends double logarithmically on the cardinality of action spaces, which is almost dimension-free. Maximizing the learning rate, the iteration complexity is bounded by $\left(1 + \|A\|_{\infty}/\tau \right)\log(1/\epsilon)$ (modulo log factors), which only depends on the ratio $\|A\|_{\infty}/\tau$.% which further decreases as we increase the amount of regularization. 
% A similar discussion holds also for the value gaps of the players in \eqref{eq:EG_conv_gap_last}.

\item {\bf Entrywise error of the policy log-ratios.} 
Both PU and OMWU enjoy strong entrywise guarantees in the sense we can guarantee the convergence of the $\ell_\infty$ norm of the log-ratios between the learned policy pair and the QRE at the same dimension-free linear rate (cf.~\eqref{eq:EG_conv_log}), which suggests the policy pair converges in a somewhat uniform manner across the entire action space.% much stronger than prior literature has to offer. 

\item {\bf Linear convergence of optimality and duality gaps.} 
Our theorem also establishes the last-iterate convergence of the game values in terms of the optimality gap (cf.~\eqref{eq:EG_conv_f_last}) and the duality gap (cf.~\eqref{eq:EG_conv_gap_last})
for both PU and OMWU. In particular, as will be seen, bounding the optimality gap of matrix games turns out to be the key enabler for generalizing our algorithms to Markov games, and bounding the duality gap allows to directly translate our results to finding a NE of unregularized matrix games. %finding 

%We point out that relation~\eqref{eq:EG_conv_gap_last} again demonstrates a near dimension-free linear convergence for both PU and OMWU which allows for more detailed comparisons with prior work on convergence of duality gaps. 

%As with the duality gaps, by definition, $\max_{\mu'\in \Delta(\mathcal{A})} f(\mu', {\nu}^{(t)}) - \min_{\nu'\in\Delta(\mathcal{B})}f({\mu}^{(t)}, \nu')$
%always stays non-negative.

%\item {\bf Near dimension-free linear convergence to game values.} Our theorem also establishes an explicit linear rate of convergence in terms of the sub-optimality gap of the game value \eqref{eq:EG_conv_f_last}. It turns out to be the key enabler of the analysis to approximate value iteration for solving Markov games using policy extragradient methods, as shall be seen in Section~\ref{sec:markov_game}.

%While PU requires an additional call to the sampling oracle every iteration compared to OMWU, it allows for a wider range of choice for the stepsize, making it more flexible in practice.  

%    Having entropy regularization ensures uniqueness of the optimal policy, making it feasible to consider the convergence of the policy directly.
 %  Under slightly stringent condition of the stepsize while still being an universal constant, we are able to establish the convergence of log policy in term of the $\|\cdot\|_{\infty}$ norm. 

\end{itemize}

\begin{figure}[t]
    \centering
    \begin{tabular}{cc}
    \includegraphics[width=0.46\linewidth]{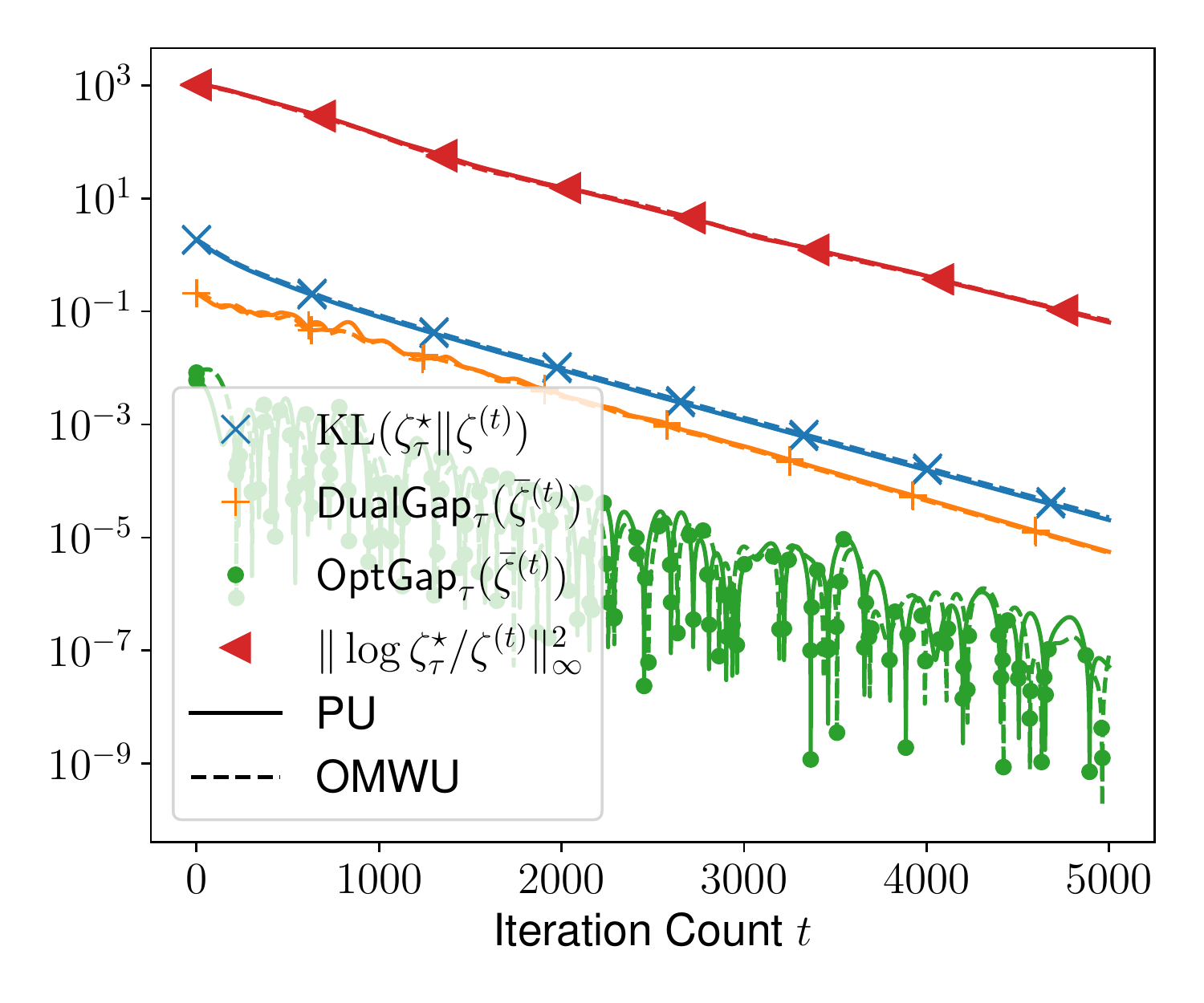} &
    \includegraphics[width=0.46\linewidth]{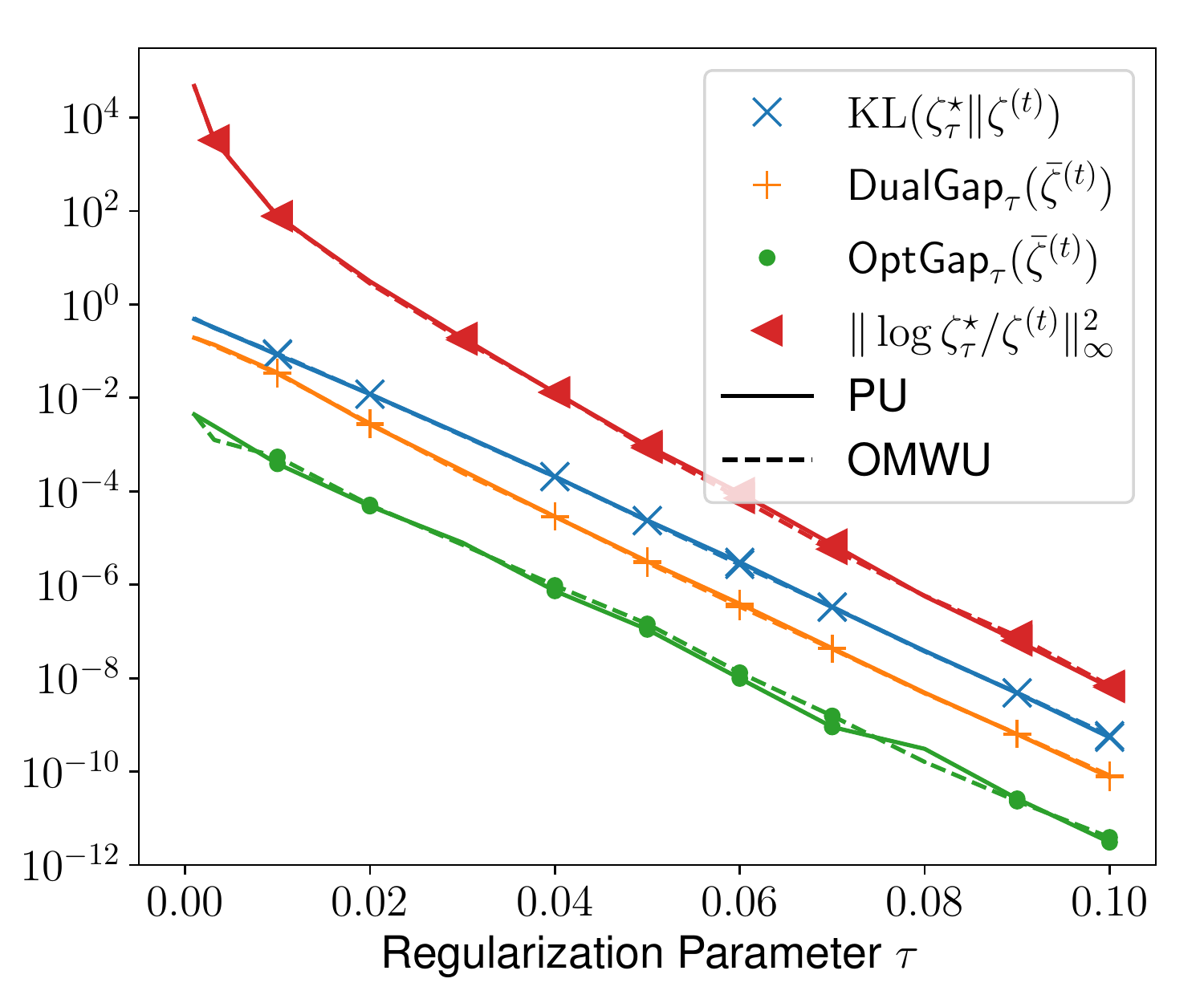}
    \end{tabular}
    \caption{Performance illustration of the PU and OMWU methods for solving entropy-regularized matrix games with $|\mathcal{A}| = |\mathcal{B}| = 100$, where the entries of the payoff matrix $A$ is generated independently from the uniform distribution on $[-1, 1]$. The learning rates are fixed as $\eta = 0.1$. The left panel plots various error metrics of convergence w.r.t. the iteration count with the entropy regularization parameter $\tau = 0.01$, while the right panel plots these error metrics at $1000$-th iteration with different choices of $\tau$. Due to their similar nature, PU and OMWU yield almost identical convergence behaviors and overlapping plots.}
    \label{fig:extragradient}
\end{figure}

Figure~\ref{fig:extragradient} illustrates the performance of the proposed PU and OMWU methods
for solving randomly generated entropy-regularized matrix games. It is evident that both algorithms converge linearly, and achieve faster convergence rates when the regularization parameter increases.

 \paragraph{Last-iterate convergence to approximate NE.} 
 The entropy-regularized matrix game can be thought as a smooth surrogate of the unregularized matrix game \eqref{eq:max-min-problem-games}; in particular, it is possible to find an $\epsilon$-NE by setting $\tau$ sufficiently small in \eqref{eq:QRE}.  According to \cite[Definition 2.1]{zhang2020model}, a policy pair $\zeta =(\mu,\nu)$ is an $\epsilon$-NE if it satisfies 
 $$  \mathsf{DualGap}(\zeta)  := \max_{\mu'\in \Delta(\mathcal{A})} f(\mu', \nu) - \min_{\nu'\in\Delta(\mathcal{B})}f(\mu, \nu' ) \leq \epsilon .$$
 
 Observe that setting $\tau = \frac{\epsilon/4}{\log|\mathcal{A}| + \log|\mathcal{B}|}$ guarantees  
$$ \left|f_{\tau}(\mu, \nu)  -f(\mu, \nu) \right| <\epsilon/4 \qquad \mbox{for all} \quad (\mu, \nu) \in \simplexA \times \simplexB$$ 
in view of the boundedness of the Shannon entropy $\mathcal{H}(\cdot)$. 
Theorem~\ref{eq:EG-conv-all} (cf.~\eqref{eq:EG_conv_gap_last}) also ensures that our proposed algorithms find an approximate QRE $\bar{\zeta}^{(T)}$ such that $\mathsf{DualGap}_\tau(\bar{\zeta}^{(T)}) \leq \epsilon/2$ after taking $T = \widetilde{O}\left(\frac{1}{\eta \epsilon} \right)$ iterations, which is no more than 
$$\widetilde{O} \left(\frac{\|A\|_{\infty} }{\epsilon}    \right)$$ 
iterations with optimized learning rates. 
It follows immediately that
 \begin{align}
    \mathsf{DualGap}(\bar{\zeta}^{(T)}) &  \leq  \mathsf{DualGap}_\tau(\bar{\zeta}^{(T)}) + \max_{\mu',\nu'}   \left|f_\tau(\mu', \bar{\nu}^{(T)}) - f_\tau(\bar{\mu}^{(T)}, \nu') - (f(\mu', \bar{\nu}^{(T)}) - f(\bar{\mu}^{(T)}, \nu')) \right|   \leq \epsilon, % \nonumber \\
 %  &\leq \epsilon/2 +\max_{\mu'}   |f_\tau(\mu', \bar{\nu}^{(T)}) - f(\mu', \bar{\nu}^{(T)})+ + \max_{\nu'}  | f_\tau(\bar{\mu}^{(T)}, \nu') - f(\bar{\mu}^{(T)}, \nu')|  .
    \label{eq:matrix_dualgap_unreg}
\end{align}
and therefore $\bar{\zeta}^{(T)}$ is an $\epsilon$-NE. Intriguingly, unlike prior work \citep{daskalakis2018last,wei2020linear} that analyzed the last-iterate convergence of OMWU in the unregularized setting ($\tau=0$), our last-iterate convergence does not require the NE of \eqref{eq:max-min-problem-games} to be unique. See Table~\ref{table:comparison} for further comparisons.

\begin{remark}
For simplicity, we have set the regularization parameter $\tau$ on the order of the final accuracy $\epsilon$. In practice, it might be desirable to use an annealing schedule of $\tau$ similar to the doubling trick, see e.g. \cite{yang2020finding,li2020sample}. We omit such straightforward generalizations for conciseness.  
\end{remark}

\paragraph{Rationality.} Another attractive feature of the algorithms developed above is being \emph{rational} (as introduced in \cite{bowling2001rational}) in the sense that the algorithm returns the best-response policy of one player when the opponent takes any \emph{fixed} stationary policy.  
More specially, in terms of matrix games, when player 2 sticks to a stationary policy $\nu$, the update of player 1 reduces to 
\begin{align}
\label{eqn:single-agent}
     \mu^{(t+1)}(a) \propto {\mu^{(t)}(a)}^{1-\eta\tau}\exp(\eta [A \nu]_a).
\end{align}
In this case, Theorem~\ref{thm:EG-gurantees-last-iterate} 
% and \ref{thm:EG-gurantees-Ave} 
can be established in exactly the same fashion by restricting attention only to the updates of $\mu^{(t)}$.

\subsection{No-regret learning of entropy-regularized OMWU}
%

%\begin{algorithm}[ht]
%    \DontPrintSemicolon
%       \textbf{inputs:} initialize $\bar{\nu}^{(0)} = \nu^{(0)}$. \\
%        % number of iterations $T$, 
%      %\textbf{initialization:} $Q_0=0$. \\
%    
%    
%       \For{$t=0,1,2,\cdots$}
%        {
%            Play $\bar{\nu}^{(t)}$ and observe $A^\top \bar{\mu}^{(t)}$.  
%
%            Update $\nu$ and $\bar{\nu}$ according to  
%            \begin{align*}
%                    {\nu}^{(t)}(b) & \propto {\nu^{(t-1)}(b)}^{1-\eta_{t-1}\tau}\exp(-\eta_{t-1} [A^\top \bar{\mu}^{(t)}]_b),   \\
%                    \bar{\nu}^{(t+1)}(b) & \propto {\nu^{(t)}(b)}^{1-\eta_{t}\tau}\exp(-\eta_{t} [A^\top \bar{\mu}^{(t)}]_b).
%            \end{align*} 
%        }
%        \caption{Single Agent OMWU Method for Entropy-regularized Matrix Games}
%     \label{alg:OMWU_single}
%\end{algorithm}
Besides convergence to equilibria, in game-theoretical settings, it is often desirable to design and implement no-regret algorithms, which are capable of providing black-box guarantees over arbitrary sequences played by the opponent \citep{cesa2006prediction,rakhlin2013optimization}. Therefore, no-regret algorithms provide a sort of robustness especially when operating in contested environments, when the opponent is potentially adversarial. Fortunately, it turns out that entropy regularization not only accelerates the convergence, but also enables no-regret learning somewhat ``for free'': it encourages exploration by putting a positive mass on every action, therefore guards against  adversaries. Moreover, since algorithms that call the sampling oracle more than once per iteration often incurs a linear regret \citep{golowich2020tight}, we focus on OMWU (Algorithm~\ref{alg:omwu}) and establish it as a no-regret algorithm.

\paragraph{No-regret learning of OMWU for the entropy-regularized matrix game.} We begin by formally defining the notion of regret. Suppose that player 2 plays according to Algorithm~\ref{alg:omwu} to update $\nu^{(t)}$ and $\bar\nu^{(t)}$ based on the received payoff sequence $\{A^{\top}\bar\mu^{(t)}\}$, $t=0,1,\ldots$, whose construction potentially is deviated from the update rule of Algorithm~\ref{alg:omwu}, and even adversarial. Let 
\[
    f_\tau^{(t)}(\nu) = \bar{\mu}^{(t)\top} A \nu  + \tau \mathcal{H}(\bar{\mu}^{(t)}) - \tau \mathcal{H}(\nu),
\]
which is the regularized game value upon fixing the policy of player 1 as $\mu = \bar{\mu}^{(t)}$. The regret experienced by player 2 is then defined as
\begin{equation}\label{eq:regret_def}
    \mathsf{Regret}_\lambda(T) = \sum_{t=0}^T f_\tau^{(t)}(\bar{\nu}^{(t)}) - \min_{\nu \in \Delta(\mathcal{B})}\sum_{t=0}^T f_\tau^{(t)}(\nu) ,
\end{equation}
which measures the gap between the actual performance and the performance in hindsight.
The following theorem shows that with appropriate choices of the learning rate, OMWU achieves a logarithmic regret bound $O(\log T)$.
\begin{theorem}\label{thm:no_regret}
Suppose only one player (say, player 2) follows the entropy-regularized OMWU method in Algorithm \ref{alg:omwu}. Setting the learning rate as $\eta_{t} = \frac{1}{(t+1)\tau}$ and the initialization policy as the uniform policy, i.e. $\nu^{(0)}(b) = 1/|\cB|, \forall b \in \cB$, the regret against any sequence $\{A^{\top}\bar\mu^{(t)}\}_{t=0}^T$ of play is bounded by
\[
    \mathsf{Regret}_\lambda(T)  \le \tau^{-1}(\log T + 1)(\tau\log |\mathcal{B}| + 5\norm{A}_\infty)^2 + 4\norm{A}_\infty.
\]
\end{theorem}
 
 Theorem~\ref{thm:no_regret} implies that the average regret satisfies
 $$ \frac{1}{T} \mathsf{Regret}_\lambda(T) \lesssim \frac{\log T}{T} ,$$
 which goes to zero as $T$ increases, therefore implies the entropy-regularized OMWU method is no-regret.

\paragraph{No-regret learning for the unregularized matrix game.} Similar to earlier discussions, one can still hope to control the regret of the unregularized matrix game, by appropriately setting $\tau$ sufficiently small.       It is easily seen that the regret of the unregularized matrix game is given by
    \begin{align*}
        \mathsf{Regret}_0(T) &= \sum_{t=0}^T f_0^{(t)}(\bar{\nu}^{(t)}) - \min_{\nu \in \Delta(\mathcal{B})}\sum_{t=0}^T f_0^{(t)}(\nu) = \max_{\nu \in \Delta(\mathcal{B})} \sum_{t=0}^T  \left [ f_0^{(t)}(\bar{\nu}^{(t)}) -   f_0^{(t)}(\nu) \right] \\
        &\le \mathsf{Regret}_\lambda(T)  + \max_{\nu \in \Delta(\mathcal{B})}\Big\{-\tau\sum_{t=0}^T \mathcal{H}(\bar{\nu}^{(t)}) + \tau T \mathcal{H}(\nu)\Big\}\\
        &\le \mathsf{Regret}_\lambda(T) + 2\tau T \log |\mathcal{B}|\\
        &\le \tau^{-1}(\log T + 1)(\tau\log |\mathcal{B}| + 5\norm{A}_\infty)^2 + 4\norm{A}_\infty + 2\tau T \log |\mathcal{B}|.
    \end{align*}
Therefore, by setting $\tau = O((\log T / T)^{1/2})$, we can ensure that the regret with regard to the unregularized problem is bounded by $\mathsf{Regret}_0(T) = O((T\log T)^{1/2})$, which is on par with the regret established in \cite{rakhlin2013optimization}. It is worthwhile to highlight that the OMWU method in \citep{rakhlin2013optimization} requires blending in a uniform distribution every iteration to guarantee no-regret learning, while a similar blending is enabled in ours without extra algorithmic steps.

\section{Zero-sum Markov games with entropy regularization}
\label{sec:markov_game}
Leveraging the success of PU and OMWU in solving the entropy-regularized matrix games, this section extends our current analysis to solve the zero-sum two-player Markov game, which is again formulated as finding the equilibrium of a saddle-point optimization problem.  
We start by introducing its basic setup, along with the entropy-regularized Markov game, which will be followed by the proposed policy extragradient method with its theoretical guarantees. The proofs for this section are collected in Appendix~\ref{sec:markov-games-proof-main}.

\subsection{Background and problem formulation}
\label{sec:formulation_markov}

\paragraph{Zero-sum two-player Markov game.} 
We consider a discounted Markov Game (MG) which is defined as $\mathcal{M} = \{\ssp, \mathcal{A}, \mathcal{B}, P, r, \gamma\}$, with discrete state space $\ssp$, action spaces of two players $\mathcal{A}$ and $\mathcal{B}$, transition probability $P$, reward function $r: \ssp\times\mathcal{A}\times\mathcal{B} \to [0, 1]$ and discount factor $\gamma \in [0, 1)$. 
A policy $\mu: \ssp \to \Delta (\mathcal{A})$ (resp. $\nu: \ssp \to \Delta (\mathcal{B})$) defines how player 1 (resp. player 2) reacts to a given state $s$, where the probability of taking action $a \in \mathcal{A}$ (resp. $b\in\mathcal{B}$) is $\mu(a|s)$ (resp. $\nu(b|s)$). The transition probability kernel $P: \ssp\times \mathcal{A} \times \mathcal{B} \to \Delta(\ssp)$ defines the dynamics of the Markov game, where $P(s'|s,a,b)$ specifies the probability of transiting to state $s'$ from state $s$ when the players take actions $a$ and $b$ respectively. The state value of a given policy pair $(\mu, \nu)$ is evaluated by the expected discounted cumulative reward: 
\[
    V^{\mu, \nu}(s) = \mathbb{E} \left[ \sum_{t=0}^{\infty} \gamma^t r(s_t, a_t, b_t ) \,\big|\, s_0 =s \right],
\]
where the trajectory $(s_0, a_0, b_0, s_1, \cdots)$ is generated by the MG $\mathcal{M}$ under the policy pair $(\mu, \nu)$, starting from the state $s_0$. Similarly, the Q-function captures the expected discounted cumulative reward with an initial state $s$ and initial action pair $(a, b)$ for a given policy pair $(\mu, \nu)$:
\[
    Q^{\mu, \nu}(s, a, b) = \mathbb{E} \left[ \sum_{t=0}^{\infty} \gamma^t r(s_t, a_t, b_t ) \,\big|\, s_0 =s, a_0 = a, b_0 = b \right].
\]

In a zero-sum game, one player seeks to maximize the value function while the other player wants to minimize it. The minimax game value on state $s$ is defined by
\[
    V^\star(s) = \max_{\mu}\min_{\nu}V^{\mu, \nu}(s) = \min_{\nu}\max_{\mu}V^{\mu, \nu}(s).
\]
Similarly, the minimax Q-function $Q^\star(s, a, b)$ is defined by
\begin{equation} 
\label{eq:Q_star_game}
Q^\star(s, a, b) = r(s, a, b) + \gamma\mathbb{E}_{s'\sim P(\cdot|s, a, b)}V^\star(s').
\end{equation}
It is proved by \cite{shapley1953stochastic} that a pair of stationary policy $(\mu^{\star}, \nu^{\star})$ attaining the minimax value on state $s$ attains the minimax value on all states as well \citep{filar2012competitive}, and is called the NE of the MG.
%Therefore, we can pick any initial state distribution $\rho$ over the state space $\ssp$, and frame the problem of finding the minimax value as a sadde-point optimization problem. 
%With a slight abuse of notation, denoting $V^{\mu,\nu}(\rho):=\mathbb{E}_{s\sim \rho}[V^{\mu,\nu}(s)]$, a policy pair $(\mu^{\star},\nu^{\star})$ is said to be a NE of the MG if it solves
% the following saddle-point optimization problem:
%\begin{align}
% \min_{\mu}\max_{\nu}V^{\mu, \nu}(\rho).
%\end{align}

\paragraph{Entropy-regularized zero-sum two-player Markov game.} Motivated by entropy regularization in Markov decision processes (MDP) \citep{geist2019theory,cen2020fast}, we consider an entropy-regularized variant of MG, where the value function is modified as
	\begin{align}
		\label{defn:V-tau}
		\soft{\V}^{\mu,\nu}(s) & := \mathbb{E} \left[ \sum_{t=0}^{\infty} \gamma^t \left( r(s_t, a_t, b_t )  - \tau \log \mu(a_t |s_t) + \tau \log \nu( b_t|s_t) \right) \,\big|\, s_0 =s \right],
	\end{align}
	where the quantity $\tau \ge 0$ denotes the regularization parameter, and the expectation is evaluated over the randomness of the transition kernel as well as the policies. The regularized Q-function $\soft{\Q}^{\mu,\nu}$ of a policy pair $(\mu,\nu)$ is related to $\soft{\V}^{\mu,\nu}$ as   
	\begin{align}
\soft{\Q}^{\mu,\nu} (s,a,b) & = r(s,a,b) + \gamma \mathbb{E}_{s'\sim P(\cdot|s,a,b)} \big[ V_{\tau}^{\mu,\nu}(s') \big]. 		\label{eq:defn-regularized-Q} 
%	\forall s \in \cS:\quad \qquad V^{\pi}_{\tau}(s) &=  \mathbb{E}_{a \sim \pi(\cdot|s)}  \big[ -\tau\log \pi(a|s) + \soft{\Q}^{\pi} (s,a) \big].\label{eq:regularized-V-to-Q} 
	\end{align}  
We will call $\soft{\V}^{\mu,\nu}$ and $\soft{\Q}^{\mu,\nu}$ the \emph{soft value function} and \emph{soft Q-function}, respectively. A policy pair $(\mu_{\tau}^{\star},\nu_{\tau}^{\star})$ is said to be the quantal response equilibrium (QRE) of the entropy-regularized MG, if its value attains the minimax value of the entropy-regularized MG over all states $s\in\cS$, i.e. 
\[
    V_\tau^\star(s) = \max_{\mu}\min_{\nu}V_\tau^{\mu, \nu}(s) = \min_{\nu}\max_{\mu}V_\tau^{\mu, \nu}(s): = V_\tau^{\mu_{\tau}^{\star},\nu_{\tau}^{\star}}(s) ,
\]
where  $V_{\tau}^{\star} $ is called 
the optimal minimax soft value function, and similarly $Q_{\tau}^{\star}: =Q_{\tau}^{\mu_{\tau}^{\star},\nu_{\tau}^{\star}}  $ is called the optimal minimax soft Q-function. %Similarly, we define $\soft{\V}^{\mu,\nu}(\rho)$ analogously when the initial state distribution is $\rho$.

%t is the solution 
%to   the following saddle-point optimization problem:
%\begin{align*}
% \min_{\mu}\max_{\nu} V_{\tau}^{\mu, \nu}(\rho) .
%\end{align*} 
\paragraph{Goal.} Our goal is to find the QRE of the entropy-regularized MG in a decentralized manner where the players only observe its own reward without accessing the opponent's actions, and leverage the QRE to find an approximate NE of the unregularized MG.

\subsection{From value iteration to policy extragradient methods}

\paragraph{Entropy-regularized value iteration.} It is known that classical dynamic programming approaches such as value iteration can be extended to solve MG \citep{perolat2015approximate}, where each iteration amounts to solving a series of matrix games for each state. Similar to the single-agent case \citep{cen2020fast}, we can extend these approaches to solve the entropy-regularized MG. Setting the stage, let us introduce the per-state Q-value matrix $Q(s): = Q(s, \cdot, \cdot) \in \mathbb{R}^{|\mathcal{A}|\times|\mathcal{B}|}$ for every $s\in \mathcal{S}$, where the element indexed by the action pair $(a, b)$ is $Q(s, a, b)$. 
Similarly, we define the per-state policies $\mu(s):=\mu(\cdot|s)\in \simplexA$ and $\nu(s):=\nu(\cdot|s)\in\simplexB$ for both players. 

In parallel to the original Bellman operator, we denote the \emph{soft Bellman operator} $\mathcal{T}_{\tau}$ as
  \begin{align} \label{eq:soft_bellman_op}
    \mathcal{T}_{\tau}(Q)(s,a,b) 
    & := r(s,a,b) + \gamma \mathbb{E}_{s'\sim P(\cdot|s,a,b)}\left[ \max_{\mu(s') \in \Delta(\asp)}\min_{\nu(s') \in \Delta(\cB)} f_{\tau} \big(Q(s'); \mu(s'),\nu(s') \big) \right],
  \end{align}
  where for each per-state Q-value matrix $Q(s)$, we introduce an entropy-regularized matrix game in the form 
of 
\begin{align*}
    \max_{\mu\in \simplexA}\min_{\nu\in \simplexB} \; f_{\tau}\big(Q(s) ; \mu(s),\nu(s) \big) :=  \mu(s)^{\top}Q(s) \nu(s) - \tau \mathcal{H}(\mu(s)) + \tau \mathcal{H}(\nu(s)).
\end{align*}
The entropy-regularized value iteration then proceeds as
 \begin{align}
 \label{eq:entropy_vi}
 Q^{(t+1)} = \mathcal{T}_{\tau}(Q^{(t)}),
\end{align}
where $Q^{(0)}$ is an initialization. 
By definition, the optimal minimax soft Q-function obeys $\mathcal{T}_{\tau}(Q_{\tau}^{\star})= Q_{\tau}^{\star}$ and therefore corresponds to the fix point of the soft Bellman operator.
Given the above entropy-regularized value iteration, the following lemma states its iterates contract linearly to the optimal minimax soft Q-function at a rate of the discount factor $\gamma$. 

\begin{prop}
\label{prop:MG-cont}
The entropy-regularized value iteration \eqref{eq:entropy_vi} converges at a linear rate, i.e.
\begin{align}
\label{eqn:MG-cont}
\|  Q^{(t)}  - Q_{\tau}^{\star} \|_{\infty} \leq \gamma^{t} \|  Q^{(0)}  - Q_{\tau}^{\star} \|_{\infty}.
\end{align}
\end{prop}

% Starting from an arbitrary $V^{(0)}\in \mathbb{R}^{|\ssp|}$, we can find $Q^\star$ by applying the following update:
%\begin{algorithm}[ht]
%    \DontPrintSemicolon
%       \textbf{inputs:} initialization $V^{(0)}$. \\
%        % number of iterations $T$, 
%      %\textbf{initialization:} $Q_0=0$. \\
%       \For{$t=0,1,2,\cdots$}
%        {
%            Compute the associated $Q^{(t)}$:
%            \[
%                Q^{(t)}(s, a, b) = r(s, a, b) + \gamma\mathbb{E}_{s'\sim P(\cdot|s, a, b)}V^{(t)}(s').
%            \]\\
%            Solve the following entropy-regularized matrix game for every state $s$:
%            \[
%                V^{(t+1)}(s) = \min_{\mu(s) \in \Delta(\asp)}\max_{\nu(s) \in \Delta(\asp)} f_{Q^{(t)}(s)}(\mu(s), \nu(s)),
%            \]
%            where
%            \[f_{Q^{(t)}(s)}(\mu(s), \nu(s)) := \mu(s)^\top Q^{(t)}(s) \nu(s) - \tau \mathcal{H}(\mu(s)) + \tau \mathcal{H}(\nu(s)).\]
%        }
%        \caption{Entropy-regularized Value Iteration}
%     \label{alg:value_iter}
%    \end{algorithm}
    % The update can be written compactly as
    % \[
    %     V^{(t+1)} = \min_{\mu}\max_{\nu} \mathcal{T}_{\mu, \nu} V^{(t)},
    % \]
    % where
    % \[
    %     \mathcal{T}_{\mu, \nu} V = r_{\mu, \nu}
    % \] 

\begin{algorithm}[!t]
    \DontPrintSemicolon
    %    \textbf{inputs:} . \\
        % number of iterations $T$, 
      \textbf{initialization:} $Q^{(0)}=0$. \\

       \For{$t=0,1,2,\cdots, T_{\text{\rm main}}$}
        {
            Let $Q^{(t)}$ denote
            \begin{equation}
                Q^{(t)}(s, a, b) = r(s, a, b) + \gamma\mathbb{E}_{s'\sim P(\cdot|s, a, b)}V^{(t)}(s'). 
                \label{eq:VI_update}
            \end{equation}

            Invoke PU (Algorithm \ref{alg:ex_grad}) or OMWU (Algorithm~\ref{alg:omwu}) for $T_{\rm sub}$ iterations to solve the following entropy-regularized matrix game for every state $s$, where the initialization is set as uniform distributions:  
            \[
                \max_{\mu(s) \in \Delta(\asp)}\min_{\nu(s) \in \Delta(\cB)} f_{\tau} \big( Q^{(t)}(s); \mu(s), \nu(s) \big).
            \]
             Return the last iterate $\bar\mu^{(t, T_{\rm sub})}(s),\bar\nu^{(t, T_{\rm sub})}(s)$. \\
            
            Set $V^{(t+1)}(s) = f_{\tau} \big( Q^{(t)}(s) ; \bar\mu^{(t, T_{\rm sub})}(s), \bar\nu^{(t, T_{\rm sub})}(s) \big)$.
            % Set $V^{(t+1)}(s) = f_{Q^{(t)}(s)}(\bar\mu^{(t, k)}(s), \bar\nu^{(t, k)}(s))$, where $k$ is the random variable defined in Theorem 1.
        }
        \caption{Policy Extragradient Method applied to Value Iteration for Entropy-regularized Markov Game}
     \label{alg:approx_value_iter}
    \end{algorithm}

\paragraph{Approximate value iteration via policy extragradient methods.} 
Proposition~\ref{prop:MG-cont} suggests that the optimal minimax soft Q-function of the entropy-regularized MG can be found by solving a series of entropy-regularized matrix games induced by $\{Q^{(t)}\}_{t\geq 0}$ in \eqref{eq:entropy_vi}, a task that can be accomplished by adopting the fast extragradient methods developed earlier. To proceed, we first define the following first-order oracle, which makes it rigorous that the proposed algorithm does not require access to the Q-function of the entire MG, but only its own single-agent Q-function when playing against the opponent's policy. 

\begin{definition}
[first-order oracle for Markov games] Given any policy pair $\mu(s), \nu(s)$ and Q-value matrix $Q(s)$ for any $s\in \mathcal{S}$, the first-order oracle returns
        \begin{align*}
            [Q(s)\nu(s)]_a &= \mathbb{E}_{b\sim \nu(s)}\left[Q(s,a,b)\right], \qquad \mbox{and}\qquad
            [Q(s)^\top\mu(s)]_b = \mathbb{E}_{a\sim \mu(s)}\left[Q(s,a,b)\right]
        \end{align*}
        for any $a\in\mathcal{A}$ and $b\in\mathcal{B}$.
\end{definition}

%For presentation simplicity, we start by assuming the proposed policy extragradient method (detailed in Algorithm \ref{alg:approx_value_iter}) samples from the first-order oracle in an exact manner, for example by using perfect knowledge of the MG model. 
Algorithm \ref{alg:approx_value_iter} describes the proposed policy extragradient method.
Encouragingly, by judiciously setting the number of iterations in both the outer loop (for updating the Q-value matrices) and the inner loop (for updating the QRE of the corresponding Q-value matrix), we are guaranteed to find the QRE of the entropy-regularized MG in a small number of iterations, as dictated by the following theorem. 

%We defer the proof of this result to Section~\ref{Sec:pf-mg-result}.

\begin{theorem}
\label{Thm:MG-result}
   Assume $|\mathcal{A}| \ge |\mathcal{B}|$ and $\tau \leq 1$. 
    Setting $\eta_t = \eta = \frac{1-\gamma}{2(1 + \tau (\log |\cA|+1-\gamma))}$, the total iterations (namely, the product $T_{\text{\rm main}}\cdot T_{\text{\rm sub}}$) required for Algorithm~\ref{alg:approx_value_iter} to achieve $\norm{Q^{(T_{\text{\rm main}})} - Q_\tau^\star}_\infty \le \epsilon$ is at most 
    \[
         O\prn{\frac{ 
         (  \log |\mathcal{A}| +1/\tau)}{ (1-\gamma)^2}
        \prn{\log \frac{\log |\mathcal{A}|}{(1-\gamma)\epsilon}}^2}.
    \]
\end{theorem}
Theorem~\ref{Thm:MG-result} ensures that 
within $\widetilde{O}\prn{\frac{1}{\tau(1-\gamma)^2}\log^2 \left( \frac{1}{\epsilon}\right) }$ iterations, Algorithm \ref{alg:approx_value_iter} finds a close approximation of the optimal minimax soft Q-function $Q_{\tau}^{\star}$ in an entrywise manner to a prescribed accuracy $\epsilon$. Remarkably, the iteration complexity is independent of the dimensions of the state space and the action space (up to log factors). In addition, the iteration complexity becomes smaller when the amount of regularization increases.

\begin{figure}[t]
    \centering
    \begin{tabular}{cc}
    \includegraphics[width=0.46\linewidth]{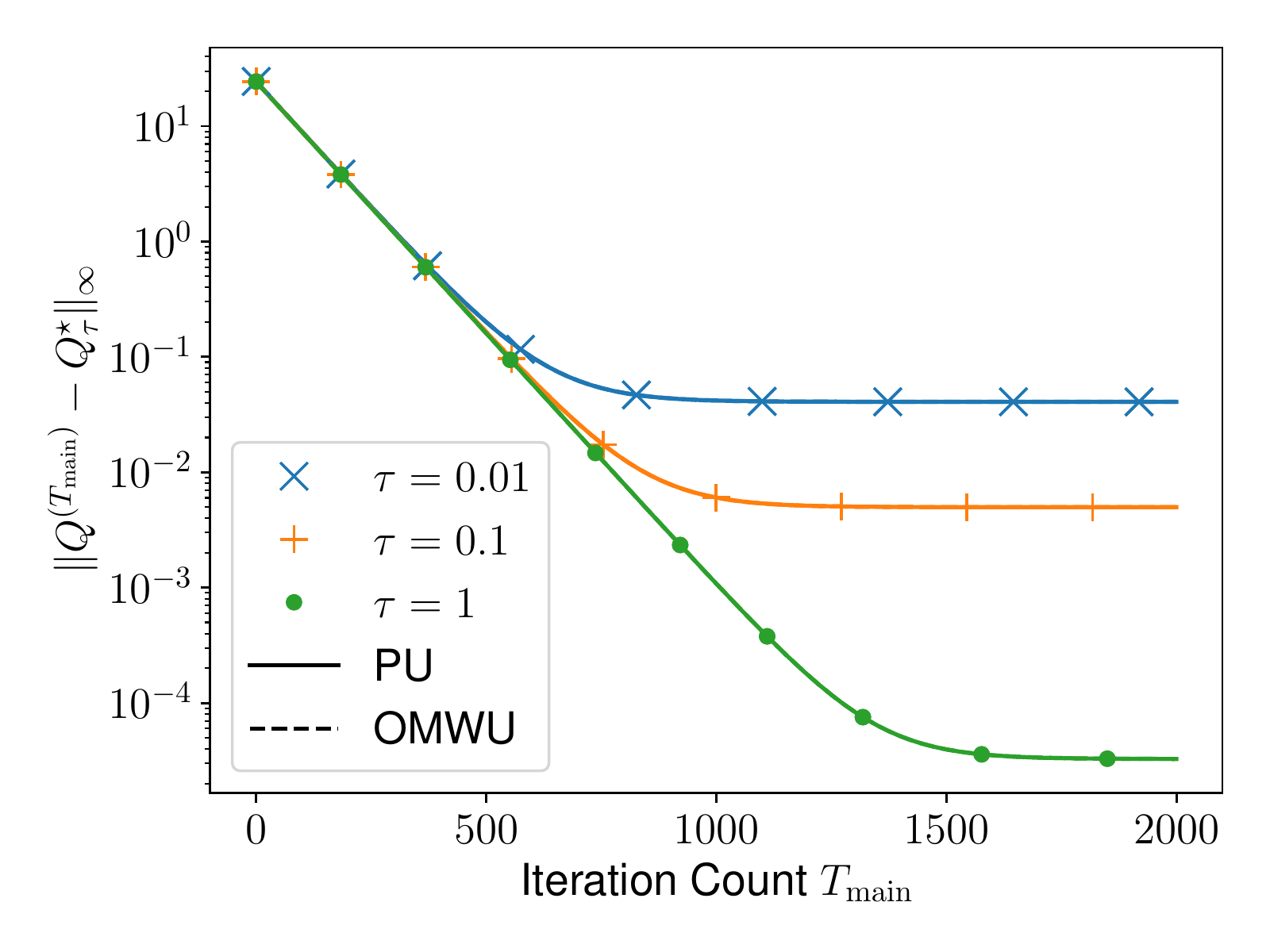} &
    \includegraphics[width=0.46\linewidth]{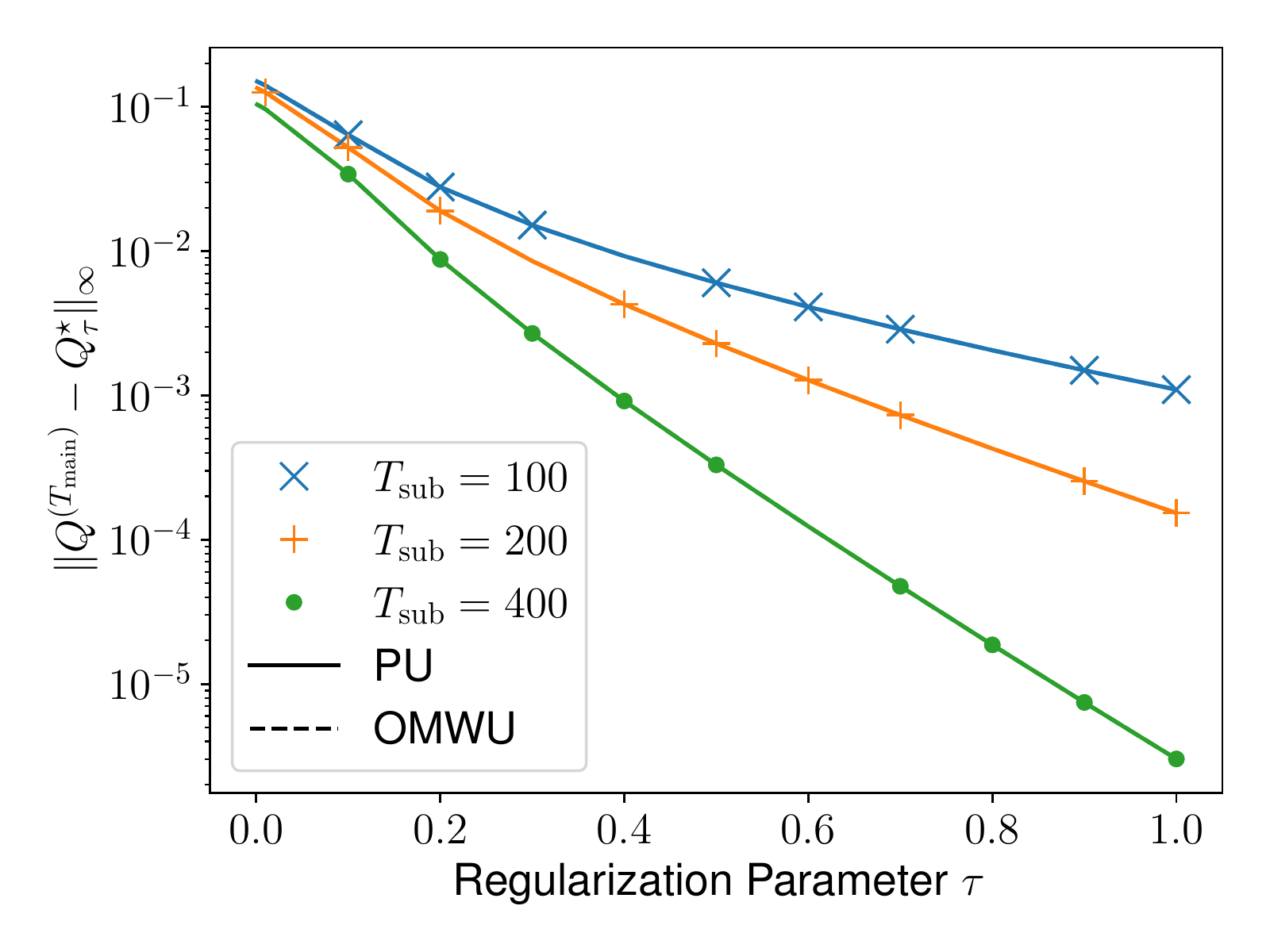}
    \end{tabular}
    \caption{Performance illustration of Algorithm~\ref{alg:approx_value_iter} for solving a random generated entropy-regularized Markov game with $|\mathcal{A}| = |\mathcal{B}| = 20$, $|\cS| = 100$ and $\gamma = 0.99$. The learning rates of both players are fixed as $\eta = 0.005$. The left panel plots $\norm{Q^{(T_{\text{\rm main}})} - Q_\tau^\star}_\infty$ w.r.t. the iteration count $T_{\rm main}$ when $T_{\rm sub}=400$ under various entropy regularization parameters, while the right panel plots  $\norm{Q^{(T_{\text{\rm main}})} - Q_\tau^\star}_\infty$ w.r.t. the regularization parameter $\tau$ when $T_{\rm main} = 2000$ with different choices of $T_{\rm sub}$. Due to their similar nature, PU and OMWU yield almost identical convergence behaviors and overlapping plots.}
    \label{fig:markov}
\end{figure}

Figure~\ref{fig:markov} illustrates the performance of Algorithm~\ref{alg:approx_value_iter} for solving a randomly generated entropy-regularized Markov game with $|\mathcal{A}| = |\mathcal{B}| = 20$, $|\cS| = 100$ and $\gamma = 0.99$ with varying choices of $T_{\rm main}$, $T_{\rm sub}$ and $\tau$. Here, the transition probability kernel and the reward function are generated as follows. For each state-action pair $(s,a,b)$, we randomly select $10$ states to form a set $\mathcal{S}_{s,a,b}$, and set $P(s'|s, a, b) \propto U_{s,a,b,s'}$ if $s' \in \mathcal{S}_{s,a}$, and $0$ otherwise, where $\{U_{s,a,b,s'}\sim U[0,1]\}$ are drawn independently from the uniform distribution over $[0, 1]$. The reward function is generated by $r(s, a,b) \sim  U_{s, a, b}\cdot U_{s}$, where $U_{s,a,b}$ and $ U_{s}$ are  drawn independently from the uniform distribution over $[0, 1]$. It is seen that the convergence speed of the $\ell_{\infty}$ error on $\norm{Q^{(T_{\text{\rm main}})} - Q_\tau^\star}_\infty $ improves as we increase the regularization parameter $\tau$, which corroborates our theory. 
 
\subsection{Last-iterate convergence to approximate NE} 

Similar to the case of matrix games, solving the entropy-regularized MG provides a viable strategy to find an $\epsilon$-approximate NE of the unregularized MG, where the optimality of a policy pair is typically gauged by the duality gap \citep{zhang2020model}.
To begin, define the duality gap of the entropy-regularized MG at a policy pair $\zeta = (\mu,\nu)$ as 
    \begin{align}
    \label{eqn:def-duality-gap}
        % {}^{\rm mar}_{\rm kov}
        \mathsf{DualGap}_\tau^{\mathsf{markov}}(\zeta) & = \max_{\mu'\in \Delta(\mathcal{A})^{|\cS|}} V^{\mu', \nu}_\tau(\rho) - \min_{\nu'\in\Delta(\mathcal{B})^{|\cS|}}V_\tau^{\mu, \nu'}(\rho),
   %      = \max_{\mu'\in \Delta(\mathcal{A}),\nu'\in\Delta(\mathcal{B})} \left\{ f_\tau(\mu', \nu) - f_\tau(\mu, \nu')\right\},
    \end{align}
where $\rho$ is an arbitrary distribution over the state space $\cS$, and $V_\tau^{\mu, \nu}(\rho) := \ex{s\sim \rho}{V_\tau^{\mu, \nu}(s)}$.\footnote{For notation simplicity, we omitted the dependency with $\rho$ in $ \mathsf{DualGap}_\tau^{\mathsf{markov}}(\zeta)$.} Similarly, the duality gap of a policy pair $\zeta = (\mu, \nu)$ for the unregularized MG  is defined as
\[
    \mathsf{DualGap}^{\mathsf{markov}}(\zeta) = \max_{\mu'\in \Delta(\mathcal{A})^{|\cS|}} V^{\mu', \nu}(\rho) - \min_{\nu'\in\Delta(\mathcal{B})^{|\cS|}}V^{\mu, \nu'}(\rho),
\]
where  $V^{\mu, \nu}(\rho) := \ex{s\sim \rho}{V^{\mu, \nu}(s)}$. Following \cite{zhang2020model}, a pair of policy $\zeta$ is said to be $\epsilon$-approximate NE if $\mathsf{DualGap}^{\mathsf{markov}}(\zeta)  \leq \epsilon$. Notice that for any policy pair $(\mu,\nu)$, it is straightforward that 
\[
\abs{V_\tau^{\mu,\nu}(\rho)  -  V^{\mu,\nu}(\rho)} \leq \frac{\tau}{1-\gamma} (\log |\cA| +\log|\cB|).
\] 

Encouragingly, the following corollary ensures that Algorithm \ref{alg:approx_value_iter} yields a policy pair with $\epsilon$-optimal duality gap for the entropy-regularized MG.
\begin{corollary}\label{cor:unreg_MG}
    Assume $|\mathcal{A}| \ge |\mathcal{B}|$ and $\tau \leq 1$. Setting $\eta_t = \eta = \frac{1-\gamma}{2(1 + \tau (\log |\cA|+1-\gamma))}$, Algorithm \ref{alg:approx_value_iter} takes no more than 
    $\widetilde{O} \left(\frac{1}{\tau (1-\gamma)^2}\log^2\left(\frac{1}{\epsilon}\right) \right)$ 
    iterations to achieve $\mathsf{DualGap}_\tau^{\mathsf{markov}}(\zeta) \le \epsilon$.
\end{corollary} 
With Corollary~\ref{cor:unreg_MG} in place, setting the regularization parameter sufficiently small, i.e. $\tau = O\big(\frac{(1-\gamma)\epsilon}{\log |\cA|} \big)$, and invoking similar discussions as \eqref{eq:matrix_dualgap_unreg} allows us to find an $\epsilon$-approximate NE of the unregularized MG within 
$$\widetilde{O}\left(\frac{1}{(1-\gamma)^3 \epsilon} \right)$$ 
iterations. See Table~\ref{table:comparison_markov} for further comparisons with 
\cite{perolat2015approximate,wei2021last,daskalakis2020independent,zhao2021provably}. To the best of our knowledge, the proposed method is the only one that simultaneously possesses symmetric updates, problem-independent rates, and last-iterate convergence.

\section{Conclusions}

This paper develops provably efficient policy extragradient methods (PU and OMWU) for entropy-regularized matrix games and Markov games, whose last iterates are guaranteed to converge linearly to the quantal response equilibrium at a linear rate. Encouragingly, the rate of convergence is independent of the dimension of the problem, i.e. the sizes of the space space and the action space. In addition, the last iterates of the proposed algorithms can also be used to locate Nash equilibria for the unregularized competitive games without assuming the uniqueness of the Nash equilibria by judiciously tuning the amount of regularization. 

This work opens up interesting opportunities for further investigations of policy extragradient methods for solving competitive games. For example, can we develop a two-time-scale policy extragradient algorithms for Markov games where the Q-function is updated simultaneously with the policy but potentially at a different time scale, using samples, such as in an actor-critic algorithm \citep{konda2000actor}? A recent work by \cite{cen2022faster} partially answered this question under exact gradient evaluation. Can we generalize the proposed algorithms to handle more general regularization terms, similar to what has been accomplished in the single-agent setting \citep{lan2021policy,zhan2021policy}? Can we generalize the proposed algorithm to other type of games  \citep{ao2022asynchronous}? We leave the answers to future work.

%\begin{itemize}
%\item \yc{can we say something about the last iterate convergence here for the unregularized problem?} \dacong{Yet another unfulfilled dream..}
%
%\item Extragradient applied to NPG updates / two time-scale updates instead of approximate value iteration
%\end{itemize}

	\section*{Acknowledgments}
 
 S.~Cen and Y.~Chi are supported in part by the grants ONR N00014-18-1-2142 and N00014-19-1-2404, ARO W911NF-18-1-0303, NSF CCF-1901199, CCF-2007911 and CCF-2106778. 
S.~Cen is also gratefully supported by Wei Shen and Xuehong Zhang Presidential Fellowship, and Nicholas Minnici Dean's Graduate Fellowship in Electrical and Computer Engineering at Carnegie Mellon University.
Y.~Wei is supported in part by the NSF grants CCF-2007911, 
DMS-2147546/2015447, CCF-2106778 and CAREER award DMS-2143215. 
% Y.~Chen is supported in part by the grants AFOSR YIP award FA9550-19-1-0030,
% ONR N00014-19-1-2120, ARO YIP award W911NF-20-1-0097, ARO W911NF-18-1-0303, NSF CCF-1907661, IIS-1900140 and DMS-2014279,  
% and the Princeton SEAS Innovation Award.  

	\bibliographystyle{apalike}
	\bibliography{bibfileRL,bibfileGame}

	\appendix

	 \section{Analysis for entropy-regularized matrix games}
\label{sec:matrix-games-proof-main}

Before embarking on the main proof, it is useful to first consider the update rule \eqref{eq:mirror_update} that underlies both PU and OMWU, which is reproduced below for convenience:
\begin{equation}  \label{eq:mirror_update2}
    \begin{cases}
        \mu^{(t+1)}(a) \propto {\mu^{(t)}(a)}^{1-\eta\tau}\exp(\eta [A z_2]_a) , \qquad  & \text{for all } a \in \mathcal{A},\\[0.2cm]
        \nu^{(t+1)}(b) \propto {\nu^{(t)}(b)}^{1-\eta\tau}\exp(-\eta [A^\top z_1]_b) , \qquad  & \text{for all } b \in \mathcal{B},\\
    \end{cases}
\end{equation}
where $z_1 \in \Delta(\mathcal{A})$ and $z_2 \in \Delta(\mathcal{B})$.  These updates satisfy the following property, whose proof is provided in Appendix~\ref{Sec:Pf-key-lemma}.

\begin{lemma}
\label{lemma:mirror_update}
Denote $\zeta^{(t)} = ( \mu^{(t)},  \nu^{(t)} )$ and $\zeta(z) = (z_1 ,z_2 )$.  
    The update rule \eqref{eq:mirror_update2} satisfies:
    \begin{subequations}
    \begin{align}
        \innprod{\log \mu^{(t+1)} - (1-\eta\tau)\log \mu^{(t)} - \eta\tau \log \best{\mu}, z_1 - \best{\mu}} & = \eta (\best{\mu} - z_1)^{\top} A ( \best{\nu} - z_2 ),         \label{eq:ohhh_mu} \\
        \innprod{\log \nu^{(t+1)} - (1-\eta\tau)\log \nu^{(t)} - \eta\tau \log \best{\nu}, z_2 - \best{\nu}} &= - \eta ( \best{\nu} - z_1)^{\top} A ( \best{\nu} - z_2 ),\label{eq:ohhh_nu}
    \end{align}
    \end{subequations}
    and
    \begin{equation}
        \innprod{\log \zeta^{(t+1)} - (1-\eta\tau)\log \zeta^{(t)} - \eta\tau \log \best{\zeta}, \zeta(z) - \best{\zeta}} = 0.
        \label{eq:ohhh}
    \end{equation}
\end{lemma}

As we shall see, the above lemma plays a crucial role in establishing the claimed convergence results.
The next lemma gives some basic decompositions related to the game values that are helpful.
\begin{lemma} 
\label{lemma:gap_KL_diff_decomp}
For every $(\mu, \nu) \in \Delta(\cA) \times \Delta(\cB)$, the following relations hold 
    \begin{subequations}
    \begin{align}
        f_{\tau}(\best{\mu}, \nu) - f_{\tau}(\mu, \best{\nu}) & = \tau \KL{\zeta}{\best{\zeta}}, \label{eq:gap_KL} \\
         f_{\tau}(\mu, \nu) - f_{\tau}(\best{\mu}, \best{\nu}) &= ( \best{\mu} - \mu)^{\top} A ( \best{\nu}-\nu) + \tau \KL{\nu}{\best{\nu}} - \tau \KL{\mu}{\best{\mu}}. \label{eq:f_diff_decomp}
    \end{align}
    \end{subequations}
\end{lemma}

 In addition, we also make record of the following elementary lemma that is used frequently.
 \begin{lemma}
    \label{lem:log_pi_gap}
    For any $\mu_1, \mu_2 \in \Delta(\mathcal{A})$ satisfying 
    $$\mu_1(a) \propto \exp(x_1(a))\quad \mbox{and} \quad
        \mu_2(a) \propto \exp(x_2(a))    $$
    for some $x_1, x_2 \in \mathbb{R}^{|\mathcal{A}|}$,
    we have
    \[
        \norm{\log \mu_1 - \log \mu_2}_\infty \le 2 \norm{x_1 - x_2}_\infty,
          \quad \mbox{and}\quad    
         \KL{\mu_1}{\mu_2} \le \frac{1}{2}\norm{x_1-x_2- c\cdot \mathbf{1}}_\infty^2, 
    \]
      where the latter inequality holds for all $ c \in \real$.
\end{lemma}

%\begin{lemma}[Lemma 25, \cite{mei2020global}]
%    \label{lemma:KL_logit}
%    For any $\nu_1$, $\nu_2 \in \Delta(\mathcal{B})$ satisfying
%    \[
%        \nu_1(b) \propto \exp(x_1(b))\quad\text{and}\quad\nu_2(b) \propto \exp(x_2(b))
%    \]
%    for some $x_1, x_2 \in \real^{|\mathcal{B}|}$, we have
%    \[
%        \KL{\nu_1}{\nu_2} \le \frac{1}{2}\norm{x_1 - x_2 - c\cdot \mathbf{1}}_\infty^2, \forall c \in \real.
%    \]
%\end{lemma}

\subsection{Proof of Proposition~\ref{prop:iu}}
 
Setting $\zeta(z) = \zeta^{(t+1)}$ in Lemma~\ref{lemma:mirror_update}, we have
\begin{equation}
    \innprod{\log \zeta^{(t+1)} - (1-\eta\tau)\log \zeta^{(t)} - \eta\tau \log \best{\zeta},  \zeta^{(t+1)}- \best{\zeta}} = 0. 
    \label{eq:implicitohh}
\end{equation}
By the definition of the KL divergence, one has
\begin{align}
    &-\innprod{\log \zeta^{(t+1)} - (1-\eta\tau)\log \zeta^{(t)} - \eta\tau \log \best{\zeta},  \best{\zeta}}\nonumber\\
    &= -(1-\eta\tau)\innprod{\log \best{\zeta} - \log \zeta^{(t)},  \best{\zeta}} + \innprod{\log \best{\zeta} - \log \zeta^{(t+1)},  \best{\zeta}}\nonumber\\
    &= -(1-\eta\tau) \KL{\best{\zeta}}{\zeta^{(t)}} + \KL{\best{\zeta}}{\zeta^{(t+1)}},
    \label{eq:little_alg_1}
\end{align}
and similarly, 
\begin{align*}
    &\innprod{\log \zeta^{(t+1)} - (1-\eta\tau)\log \zeta^{(t)} - \eta\tau \log \best{\zeta},  {\zeta}^{(t+1)}}\nonumber\\
    &= (1-\eta\tau)\innprod{\log \zeta^{(t+1)} - \log \zeta^{(t)},  \zeta^{(t+1)}} + \eta\tau\innprod{\log \zeta^{(t+1)} - \log \best{\zeta},  \zeta^{(t+1)}}\nonumber\\
    &= (1-\eta\tau) \KL{\zeta^{(t+1)}}{\zeta^{(t)}} + \eta\tau\KL{\zeta^{(t+1)}}{\best{\zeta}}.
    % \label{eq:little_alg_2}
\end{align*}
Combining the above two equalities with \eqref{eq:implicitohh}, we arrive at
\begin{equation}
   \KL{\best{\zeta}}{\zeta^{(t+1)}} + \eta\tau\KL{\zeta^{(t+1)}}{\best{\zeta}}  + (1-\eta\tau)\KL{\zeta^{(t+1)}}{\zeta^{(t)}}  = (1-\eta\tau) \KL{\best{\zeta}}{\zeta^{(t)}}.
\end{equation}
This immediately leads to $\KL{\best{\zeta}}{\zeta^{(t+1)}} \leq (1-\eta\tau) \KL{\best{\zeta}}{\zeta^{(t)}}$
by the nonnegativity of the KL divergence, as long as $1-\eta\tau\geq 0$. Therefore 
\begin{align*}
    \KL{\best{\zeta}}{\zeta^{(t)}} \le (1-\eta\tau)^t \KL{\best{\zeta}}{\zeta^{(0)}} \qquad \mbox{for all}~t\geq 0.
\end{align*}
 %We therefore obtain a linear convergence result for updates of the form~\eqref{eq:mirror_update}. 

\subsection{Proof of Theorem~\ref{thm:EG-gurantees-last-iterate} }
\label{Sec:pf-matrix-game}

\subsubsection{Proof of policy convergence in KL divergence~\eqref{eq:EG_conv_KL}}

First noticing that both PU and OMWU share the same update rule for $\mu^{(t+1)}$ and $\nu^{(t+1)}$, which takes the form 
\begin{align*}
    \begin{cases}
        {\mu}^{(t+1)}(a) \propto {\mu^{(t)}(a)}^{1-\eta\tau}\exp(\eta [A \bar{\nu}^{(t+1)}]_a),\\
        {\nu}^{(t+1)}(b) \propto {\nu^{(t)}(b)}^{1-\eta\tau}\exp(-\eta [A^\top \bar{\mu}^{(t+1)}]_b).
    \end{cases}    
\end{align*}
Regarding this sequence, Lemma~\ref{lemma:mirror_update} (cf.~\eqref{eq:ohhh}) gives  
\begin{equation}
        \innprod{\log \zeta^{(t+1)} - (1-\eta\tau)\log \zeta^{(t)} - \eta\tau \log \best{\zeta}, \bar{\zeta}^{(t+1)} - \best{\zeta}} = 0.
        \label{eq:EGohh2}
\end{equation}
In view of the similarity of \eqref{eq:implicitohh} and \eqref{eq:EGohh2}, we can expect similar convergence guarantees to that of the implicit updates established in Proposition~\ref{prop:iu} with the optimism that $\bar{\zeta}^{(t+1)}$ approximates $\zeta^{(t+1)}$ well.
    Following the same argument as \eqref{eq:little_alg_1}, we have
    \begin{equation}
        -\innprod{\log \zeta^{(t+1)} - (1-\eta\tau)\log \zeta^{(t)} - \eta\tau \log \best{\zeta},  \best{\zeta}}= -(1-\eta\tau) \KL{\best{\zeta}}{\zeta^{(t)}} + \KL{\best{\zeta}}{\zeta^{(t+1)}}.
        \label{eq:some_alg_1}
    \end{equation}
    On the other hand, it is easily seen that
    \begin{align}
        &\innprod{\log \zeta^{(t+1)} - (1-\eta\tau)\log \zeta^{(t)} - \eta\tau \log \best{\zeta}, \bar{\zeta}^{(t+1)}}\nonumber\\
        &=\innprod{\log \bar{\zeta}^{(t+1)} - (1-\eta\tau)\log \zeta^{(t)} - \eta\tau \log \best{\zeta}, \bar{\zeta}^{(t+1)}} - \innprod{\log \bar{\zeta}^{(t+1)} - \log \zeta^{(t+1)}, {\zeta}^{(t+1)}}\nonumber\\
        &\qquad- \innprod{\log \bar{\zeta}^{(t+1)} - \log \zeta^{(t+1)}, \bar{\zeta}^{(t+1)} - {\zeta}^{(t+1)}}\nonumber\\
        &= (1-\eta\tau)\KL{\bar{\zeta}^{(t+1)}}{\zeta^{(t)}} + \eta\tau\KL{\bar{\zeta}^{(t+1)}}{\best{\zeta}} + \KL{\zeta^{(t+1)}}{\bar{\zeta}^{(t+1)}} \nonumber\\
        &\qquad- \innprod{\log \bar{\zeta}^{(t+1)} - \log \zeta^{(t+1)}, \bar{\zeta}^{(t+1)} - \zeta^{(t+1)}}.
        \label{eq:some_alg_2}
    \end{align}
Combining equalities~\eqref{eq:some_alg_1}, \eqref{eq:some_alg_2} with \eqref{eq:EGohh2}, we are left with the following relation pertaining to bounding $\KLbig{\best{\zeta}}{\zeta^{(t)}}$:
\begin{align}
    (1-\eta\tau) \KL{\best{\zeta}}{\zeta^{(t)}} & =(1-\eta\tau)\KL{\bar{\zeta}^{(t+1)}}{\zeta^{(t)}} + \eta\tau\KL{\bar{\zeta}^{(t+1)}}{\best{\zeta}} + \KL{\zeta^{(t+1)}}{\bar{\zeta}^{(t+1)}} \nonumber \\
    & \qquad - \innprod{\log \bar{\zeta}^{(t+1)} - \log \zeta^{(t+1)}, \bar{\zeta}^{(t+1)} - \zeta^{(t+1)}} + \KL{\best{\zeta}}{\zeta^{(t+1)}}. \label{eq:EGohhh2}
\end{align}

In addition, to bound $ \KL{\best{\zeta}}{\bar{\zeta}^{(t+1)}} $, we will resort to the following three-point equality, which reads
\begin{align}
    \notag \KL{\best{\zeta}}{\bar{\zeta}^{(t+1)}} &=  \KL{\best{\zeta}}{\zeta^{(t+1)}}  - \innprod{\best{\zeta}, \log \bar{\zeta}^{(t+1)} - \log \zeta^{(t+1)}}\\
    &= \KL{\best{\zeta}}{\zeta^{(t+1)}}  - \KL{\bar{\zeta}^{(t+1)}}{\zeta^{(t+1)}} - \innprod{\best{\zeta} - \bar{\zeta}^{(t+1)}, \log \bar{\zeta}^{(t+1)} - \log \zeta^{(t+1)}},
    \label{eq:three_point}
\end{align}
which can be checked directly using the definition of the KL divergence. 

To proceed, we need to control $\innprod{\log \bar{\zeta}^{(t+1)} - \log \zeta^{(t+1)}, \bar{\zeta}^{(t+1)} - \zeta^{(t+1)}}$ on the right-hand side of inequality \eqref{eq:EGohhh2}, 
and
$\innprod{\best{\zeta} - \bar{\zeta}^{(t+1)}, \log \bar{\zeta}^{(t+1)} - \log \zeta^{(t+1)}}$ on the right-hand side of inequality \eqref{eq:three_point},
for which we continue the proofs for PU and OMWU separately as follows.

\paragraph{Bounding $\KLbig{\best{\zeta}}{\zeta^{(t)}}$ for PU.}
Following the update rule of $\bar{\zeta}^{(t+1)} =(\bar{\mu}^{(t+1)}, \bar{\nu}^{(t+1)}$) in PU, we have
\begin{equation} \label{eq:pu_bar_mu}
    \log \bar{\mu}^{(t+1)} - \log \mu^{(t+1)} = \eta A (\nu^{(t)} - \bar{\nu}^{(t+1)}) + c\cdot \mathbf{1}
\end{equation}
for some  normalization constant $c$. 
With this relation in place, one has 
\begin{align*}
    \innprod{\log \bar{\mu}^{(t+1)} - \log \mu^{(t+1)}, \bar{\mu}^{(t+1)} - \mu^{(t+1)}} &= \eta  (\bar{\mu}^{(t+1)} - \mu^{(t+1)} )^{\top}A   ( \nu^{(t)} - \bar{\nu}^{(t+1)}  ) \\
   & \le \eta \norm{A}_\infty \norm{\bar{\mu}^{(t+1)} - \mu^{(t+1)}}_1 \norm{\bar{\nu}^{(t+1)} - \nu^{(t)}}_1.
\end{align*}
Combined with Pinsker's inequality, it is therefore clear that 
\begin{align}
    \innprod{\log \bar{\mu}^{(t+1)} - \log \mu^{(t+1)}, \bar{\mu}^{(t+1)} - \mu^{(t+1)}} & \le \frac{1}{2}\eta \norm{A}_\infty\prn{ \norm{\bar{\mu}^{(t+1)} - \mu^{(t+1)}}_1^2 + \norm{\bar{\nu}^{(t+1)} - \nu^{(t)}}_1^2}\nonumber\\
  &  \le \eta \norm{A}_\infty\prn{ \KL{\mu^{(t+1)}}{\bar\mu^{(t+1)}} + \KL{\bar\nu^{(t+1)}}{\nu^{(t)}}}.
    \label{eq:EG_implicit_gap}
\end{align}
Analogously, one can achieve the same bound regarding the quantity $\innprod{\log \bar{\nu}^{(t+1)} - \log \nu^{(t+1)}, \bar{\nu}^{(t+1)} - \nu^{(t+1)}}$. 
Summing up these two inequalities, we end up with  
\begin{align*}
    \innprod{\log \bar{\zeta}^{(t+1)} - \log \zeta^{(t+1)}, \bar{\zeta}^{(t+1)} - \zeta^{(t+1)}}\le \eta \norm{A}_\infty\prn{ \KL{\zeta^{(t+1)}}{\bar{\zeta}^{(t+1)}} + \KL{\bar{\zeta}^{(t+1)}}{\zeta^{(t)}}}. 
\end{align*}
Plugging the above inequality into inequality~\eqref{eq:EGohhh2} leads to
\begin{align}
    \KL{\best{\zeta}}{\zeta^{(t+1)}} & \le (1-\eta\tau) \KL{\best{\zeta}}{\zeta^{(t)}} - \big(1-\eta\tau - \eta \norm{A}_\infty \big)\KL{\bar{\zeta}^{(t+1)}}{\zeta^{(t)}} - \eta\tau\KL{\bar{\zeta}^{(t+1)}}{\best{\zeta}} \nonumber\\
    & \qquad - (1-\eta \norm{A}_\infty)\KL{\zeta^{(t+1)}}{\bar{\zeta}^{(t+1)}} .
    \label{eq:EGohhh_final}
\end{align}
Therefore, as long as the learning rate $\eta$ satisfies $\eta \le \frac{1}{\tau + \norm{A}_\infty}$,
we are ensured that
\begin{align*}
    \KL{\best{\zeta}}{\zeta^{(t+1)}} & \le(1-\eta\tau) \KL{\best{\zeta}}{\zeta^{(t)}},
\end{align*}
which further implies inequality~\eqref{eq:EG_conv_KL} when applied recursively.

\paragraph{Bounding $ \KL{\best{\zeta}}{\bar{\zeta}^{(t+1)}} $ for PU.}
By similar tricks of arriving at \eqref{eq:EG_implicit_gap}, we have
\begin{align*}
    -{\innprod{\best{\mu} - \bar{\mu}^{(t+1)}, \log \bar{\mu}^{(t+1)} - \log \mu^{(t+1)}}} &= - \eta (\best{\mu} - \bar{\mu}^{(t+1)})^{\top} A ( \nu^{(t)} - \bar{\nu}^{(t+1)} ) \\
   & \le \frac{1}{2} \eta \norm{A}_\infty\prn{\norm{\best{\mu} - \bar{\mu}^{(t+1)}}_1^2 + \norm{\nu^{(t)} - \bar{\nu}^{(t+1)}}_1^2}\\
   & \le \eta \norm{A}_\infty\prn{\KL{\best{\mu}}{\bar{\mu}^{(t+1)}} + \KL{\bar{\nu}^{(t+1)}}{\nu^{(t)}}},
\end{align*}
following from \eqref{eq:pu_bar_mu} and Pinsker's inequality.
A similar inequality for $-{\innprod{\best{\nu} - \bar{\nu}^{(t+1)}, \log \bar{\nu}^{(t+1)} - \log \nu^{(t+1)}}}$ can be obtained by symmetry, and summing together the two leads to 
\[
    -{\innprod{\best{\zeta} - \bar{\zeta}^{(t+1)}, \log \bar{\zeta}^{(t+1)} - \log \zeta^{(t+1)}}} \le \eta \norm{A}_\infty\prn{\KL{\best{\zeta}}{\bar{\zeta}^{(t+1)}} + \KL{\bar{\zeta}^{(t+1)}}{\zeta^{(t)}}}.
\]
Plugging the above inequality into \eqref{eq:three_point} and rearranging terms, we reach at
\[
    (1-\eta\norm{A}_\infty)\KL{\best{\zeta}}{\bar{\zeta}^{(t+1)}} \le \KL{\best{\zeta}}{\zeta^{(t+1)}}  + \eta \norm{A}_\infty{\KL{\bar{\zeta}^{(t+1)}}{\zeta^{(t)}}}.    
\]
Along with \eqref{eq:EGohhh_final}, we have
\begin{align}
    (1-\eta\norm{A}_\infty)\KL{\best{\zeta}}{\bar{\zeta}^{(t+1)}} & \le(1-\eta\tau) \KL{\best{\zeta}}{\zeta^{(t)}} -(1-\eta\tau - 2\eta \norm{A}_\infty)\KL{\bar{\zeta}^{(t+1)}}{\zeta^{(t)}} \nonumber\\
    &\qquad - \eta\tau\KL{\bar{\zeta}^{(t+1)}}{\zeta^{\star}}  - (1-\eta \norm{A}_\infty)\KL{\zeta^{(t+1)}}{\bar{\zeta}^{(t+1)}}.
\end{align}
Therefore, with $\eta \le 1/(\tau + 2\norm{A}_\infty)$ we have
\[
    \KL{\best{\zeta}}{\bar{\zeta}^{(t+1)}} \le 2\KL{\best{\zeta}}{\zeta^{(t)}} \le 2 (1-\eta\tau)^{t} \KL{\best{\zeta}}{\zeta^{(0)}}.
\]

%Inequalities \eqref{eq:EG_conv_KL} and \eqref{eq:EG_conv_f_last} for PU hold for a slightly wider range of stepsize namely, $0 < \eta_{\mathsf{PU}} \le \frac{1}{\tau + \norm{A}_\infty}$ as can be seen in the proof. 

\paragraph{Bounding $\KLbig{\best{\zeta}}{\zeta^{(t)}}$  for OMWU.}
 Following the update rule of $\bar{\zeta}^{(t+1)} = (\bar{\mu}^{(t+1)}, \bar{\nu}^{(t+1)})$  for OMWU, 
we have
\begin{align}
    \log \bar{\mu}^{(t+1)} - \log \mu^{(t+1)} &= \eta A (\bar{\nu}^{(t)} - \bar{\nu}^{(t+1)}) + c\cdot \mathbf{1} \nonumber \\
    &= \eta A (\bar{\nu}^{(t)} - \nu^{(t)}) + \eta A (\nu^{(t)} - \bar{\nu}^{(t+1)}) + c\cdot \mathbf{1}, \label{eq:omwu_bar_mu}
\end{align}
where $c$ is some normalization constant. 
Similar to the proof of relation~\eqref{eq:EG_implicit_gap}, it can be easily demonstrated that 
\begin{align}
    &\innprod{\log \bar{\mu}^{(t+1)} - \log \mu^{(t+1)}, \bar{\mu}^{(t+1)} - \mu^{(t+1)}} \nonumber\\
    &= \eta ( \bar{\mu}^{(t+1)} - \mu^{(t+1)} )^{\top} A ( \bar{\nu}^{(t)} - \nu^{(t)}) + \eta  ( \bar{\mu}^{(t+1)} - \mu^{(t+1)})^{\top} A (\nu^{(t)} - \bar{\nu}^{(t+1)}) \nonumber\\
    &\le \eta \norm{A}_\infty \prn{\KL{\nu^{(t)}}{\bar{\nu}^{(t)}} + \KL{\bar{\nu}^{(t+1)}}{\nu^{(t)}} + 2\KL{\mu^{(t+1)}}{\bar{\mu}^{(t+1)}}}.
    \label{eq:OMWU_implicit_gap}
\end{align}
By symmetry, we can also establish a similar inequality for $\innprod{\log \bar{\nu}^{(t+1)} - \log \nu^{(t+1)}, \bar{\nu}^{(t+1)} - \nu^{(t+1)}} $, which in turns yields 
\begin{align*}
    &\innprod{\log \bar{\zeta}^{(t+1)} - \log \zeta^{(t+1)}, \bar{\zeta}^{(t+1)} - \zeta^{(t+1)}} \\
    &\le \eta \norm{A}_\infty \prn{\KL{\zeta^{(t)}}{\bar{\zeta}^{(t)}} + \KL{\bar{\zeta}^{(t+1)}}{\zeta^{(t)}} + 2\KL{\zeta^{(t+1)}}{\bar{\zeta}^{(t+1)}}}.
\end{align*}
Plugging the above inequality into equation~\eqref{eq:EGohhh2} and re-organizing terms, we arrive at
\begin{align}
\label{eqn:omwu-nice}
 \notag   \KL{\best{\zeta}}{\zeta^{(t+1)}} &\le (1-\eta\tau) \KL{\best{\zeta}}{\zeta^{(t)}} - (1-\eta\tau - \eta \norm{A}_\infty)\KL{\bar{\zeta}^{(t+1)}}{\zeta^{(t)}} - \eta\tau\KL{\bar{\zeta}^{(t+1)}}{\best{\zeta}} \\
    &- (1 - 2\eta \norm{A}_\infty) \KL{\zeta^{(t+1)}}{\bar{\zeta}^{(t+1)}} + \eta \norm{A}_\infty\KL{\zeta^{(t)}}{\bar{\zeta}^{(t)}}.
\end{align}
With the choice of the learning rate $\eta \le \min\{\frac{1}{2\norm{A}_\infty + 2\tau}, \frac{1}{4\norm{A}_\infty}\}$, it obeys
\[
    (1-\eta\tau)(1 - 2\eta \norm{A}_\infty) \ge \eta \norm{A}_\infty.
\]
Combining the above inequality with   \eqref{eqn:omwu-nice} gives 
\begin{align*}
    & \KL{\best{\zeta}}{\zeta^{(t+1)}} + (1 - 2\eta \norm{A}_\infty) \KL{\zeta^{(t+1)}}{\bar{\zeta}^{(t+1)}}\\
    &\le (1-\eta\tau) \KL{\best{\zeta}}{\zeta^{(t)}} + \eta \norm{A}_\infty\KL{\zeta^{(t)}}{\bar{\zeta}^{(t)}} - \eta\tau\KL{\bar{\zeta}^{(t+1)}}{\best{\zeta}}.\\
    &\le (1-\eta\tau) \brk{\KL{\best{\zeta}}{\zeta^{(t)}} + (1 - 2\eta \norm{A}_\infty)\KL{\zeta^{(t)}}{\bar{\zeta}^{(t)}}} - \eta\tau\KL{\bar{\zeta}^{(t+1)}}{\best{\zeta}}.
\end{align*}
For conciseness, let us introduce the shorthand notation 
\begin{equation}\label{eq:def_Lt}
    L^{(t)} := \KL{\best{\zeta}}{\zeta^{(t)}} + (1 - 2\eta \norm{A}_\infty) \KL{\zeta^{(t)}}{\bar{\zeta}^{(t)}}.
    \end{equation} 
As a result, the above inequality can be restated as 
\begin{equation}
    L^{(t+1)} \le (1-\eta\tau) L^{(t)} - \eta\tau\KL{\bar{\zeta}^{(t+1)}}{\best{\zeta}}.
    \label{eq:OMWUohhh_final}
\end{equation}
Since we initialize OMWU with $\bar{\zeta}^{(0)} = \zeta^{(0)}$, therefore $L^{(0)} = \KL{\best{\zeta}}{\zeta^{(0)}}$, which in turn gives 
\begin{align*}
    \KL{\best{\zeta}}{\zeta^{(t)}} &\le L^{(t)} \le (1-\eta\tau)^tL^{(0)} = (1-\eta\tau)^t \KL{\best{\zeta}}{\zeta^{(0)}}.
\end{align*}
We complete the proof of inequality~\eqref{eq:EG_conv_KL} for OMWU.

\paragraph{Bounding $ \KL{\best{\zeta}}{\bar{\zeta}^{(t+1)}} $ for OMWU.}
By similar tricks of arriving at \eqref{eq:OMWU_implicit_gap}, we have
\begin{align*}
    -{\innprod{\best{\mu} - \bar{\mu}^{(t+1)}, \log \bar{\mu}^{(t+1)} - \log \mu^{(t+1)}}} 
    &= \eta ( \bar{\mu}^{(t+1)} - \best{\mu} )^{\top} A ( \bar{\nu}^{(t)} - \nu^{(t)}) + \eta  ( \bar{\mu}^{(t+1)} - \best{\mu})^{\top} A (\nu^{(t)} - \bar{\nu}^{(t+1)}) \nonumber\\
    &\le \eta \norm{A}_\infty \prn{\KL{\nu^{(t)}}{\bar{\nu}^{(t)}} + \KL{\bar{\nu}^{(t+1)}}{\nu^{(t)}} + 2\KL{\best{\mu}}{\bar{\mu}^{(t+1)}}},
\end{align*}
where the first line follows from \eqref{eq:omwu_bar_mu}.
A similar inequality also holds for $-{\innprod{\best{\nu} - \bar{\nu}^{(t+1)}, \log \bar{\nu}^{(t+1)} - \log \nu^{(t+1)}}}$. Summing the two inequalities leads to 
\[
    -{\innprod{\best{\zeta} - \bar{\zeta}^{(t+1)}, \log \bar{\zeta}^{(t+1)} - \log \zeta^{(t+1)}}} \le \eta \norm{A}_\infty\prn{\KL{\zeta^{(t)}}{\bar{\zeta}^{(t)}} + \KL{\bar{\zeta}^{(t+1)}}{\zeta^{(t)}} + 2\KL{\best{\zeta}}{\bar{\zeta}^{(t+1)}}}.
\]
Plugging the above inequality into \eqref{eq:three_point} and rearranging terms, we reach at
%\[
%    \KL{\best{\zeta}}{\bar{\zeta}^{(t+1)}} \le \KL{\best{\zeta}}{\zeta^{(t+1)}} + \eta \norm{A}_\infty\prn{\KL{\zeta^{(t)}}{\bar{\zeta}^{(t)}} + \KL{\bar{\zeta}^{(t+1)}}{\zeta^{(t)}} + 2\KL{\best{\zeta}}{\bar{\zeta}^{(t+1)}}}.
%\]
\[
    (1-2\eta\norm{A}_\infty)\KL{\best{\zeta}}{\bar{\zeta}^{(t+1)}} \le \KL{\best{\zeta}}{\zeta^{(t+1)}} + \eta \norm{A}_\infty\prn{\KL{\zeta^{(t)}}{\bar{\zeta}^{(t)}} + \KL{\bar{\zeta}^{(t+1)}}{\zeta^{(t)}}}.
\]
Along with \eqref{eqn:omwu-nice},
we have
\begin{align*}
    (1-2\eta\norm{A}_\infty)\KL{\best{\zeta}}{\bar{\zeta}^{(t+1)}} &\le (1-\eta\tau) \KL{\best{\zeta}}{\zeta^{(t)}} - (1-\eta\tau - 2\eta \norm{A}_\infty)\KL{\bar{\zeta}^{(t+1)}}{\zeta^{(t)}} - \eta\tau\KL{\bar{\zeta}^{(t+1)}}{\best{\zeta}} \\
    &\qquad - (1 - 2\eta \norm{A}_\infty) \KL{\zeta^{(t+1)}}{\bar{\zeta}^{(t+1)}} + 2\eta \norm{A}_\infty\KL{\zeta^{(t)}}{\bar{\zeta}^{(t)}} \\
    &\le (1-\eta\tau) \KL{\best{\zeta}}{\zeta^{(t)}} + 2\eta \norm{A}_\infty\KL{\zeta^{(t)}}{\bar{\zeta}^{(t)}}\\
    &\le \KL{\best{\zeta}}{\zeta^{(t)}} + (1- 2\eta \norm{A}_\infty)\KL{\zeta^{(t)}}{\bar{\zeta}^{(t)}} =: L^{(t)},
\end{align*}
where we recall the shorthand notation $L^{(t)}$ in \eqref{eq:def_Lt}.
As the learning rate of OMWU satisfies $0< \eta < \min \left\{\frac{1}{2\norm{A}_\infty + 2\tau},\, \frac{1}{4\norm{A}_\infty} \right\}$, 
it is clear that 
\begin{align*}
    \KL{\best{\zeta}}{\bar{\zeta}^{(t+1)}} \le 2L^{(t)} 
    \stackrel{(\mathrm{i})}{\le} 2 (1-\eta\tau)^{t} L^{(0)} \le 2 (1-\eta\tau)^{t} \KL{\best{\zeta}}{\zeta^{(0)}},    
\end{align*}
where (i) follows from the recursive relation $L^{(t+1)} \leq (1-\eta\tau) L^{(t)}$ shown in inequality~\eqref{eq:OMWUohhh_final}.

\subsubsection{Proof of entrywise convergence of policy log-ratios~\eqref{eq:EG_conv_log}}

%In this subsection, we establish inequality~\eqref{eq:EG_conv_log} for both PU and OMWU together. 

To facilitate the proof, we introduce an auxiliary sequence $\{\xi^{(t)}\in \mathbb{R}^{|\asp|}\}$ constructed recursively by
\begin{subequations}
\begin{align}
    \xi^{(0)}(a) &= \norm{\exp(A \best{\nu}/\tau)}_1 \cdot \mu^{(0)}(a),\\
    \xi^{(t+1)}(a) &= {\xi^{(t)}}(a)^{1-\eta\tau}\exp(\eta [A \bar{\nu}^{(t+1)}]_a),\qquad \forall a \in \asp, t\ge 0. \label{eq:xi_def}
\end{align}
\end{subequations}
It is easily seen that $\mu^{(t)}(a) \propto \xi^{(t)}(a) = \exp(\log\xi^{(t)}(a))$ for $t\geq 0$. 
Noticing that $\mu_\tau^\star \propto \exp(A\best{\nu})$, one has 
\begin{align}
\label{eqn:brahms-reduction}
    \norm{\log \mu^{(t+1)} - \log \best{\mu}}_\infty &\le 2\norm{\log\xi^{(t+1)} - A\best{\nu}/\tau}_\infty,
\end{align}
where we make use of Lemma~\ref{lem:log_pi_gap}. 

Therefore it suffices for us to control the term $\norm{\log\xi^{(t+1)} - A\best{\nu}/\tau}_\infty$ on the right-hand side of inequality~\eqref{eqn:brahms-reduction}. Taking logarithm on both sides of \eqref{eq:xi_def} yields  
\begin{align*}
    \log\xi^{(t+1)} - A\best{\nu}/\tau &= (1-\eta \tau)\log\xi^{(t)} + \eta A \bar{\nu}^{(t+1)} - A\best{\nu}/\tau\\
    &= (1-\eta \tau)\prn{\log\xi^{(t)} - A\best{\nu}/\tau} + \eta A(\bar{\nu}^{(t+1)} - \best{\nu}), 
\end{align*}
which, when combined with Pinsker's inequality, implies 
\begin{align}
    \norm{\log\xi^{(t+1)} - A\best{\nu}/\tau}_\infty &\le (1-\eta \tau)\norm{\log\xi^{(t)} - A\best{\nu}/\tau}_\infty + \eta \norm{A}_\infty\norm{\bar{\nu}^{(t+1)} - \best{\nu}}_1\nonumber\\
    &\le (1-\eta \tau)\norm{\log\xi^{(t)} - A\best{\nu}/\tau}_\infty + \eta \norm{A}_\infty\brk{2\KL{\best{\nu}}{\bar{\nu}^{(t+1)}}}^{1/2}\nonumber\\
    &\le (1-\eta \tau)\norm{\log\xi^{(t)} - A\best{\nu}/\tau}_\infty + \eta \norm{A}_\infty\brk{2\KL{\best{\zeta}}{\bar{\zeta}^{(t+1)}}}^{1/2}.
    \label{eq:log_pi_recursive}
\end{align}
Plugging the bound of $\KL{\best{\zeta}}{\bar{\zeta}^{(t+1)}}$ from relation~\eqref{eq:EG_conv_KL} into \eqref{eq:log_pi_recursive} and invoking the inequality recursively leads to
\begin{align*}
    &\norm{\log\xi^{(t+1)} - A\best{\nu}/\tau}_\infty \\
    &\le (1-\eta\tau)^{t+1}\norm{\log\xi^{(0)} - A\best{\nu}/\tau}_\infty + 2\eta \norm{A}_\infty \sum_{s=1}^{t+1}(1-\eta\tau)^{t+1-s/2} {\KL{\best{\zeta}}{\zeta^{(0)}}}^{1/2}\\
    &\le (1-\eta\tau)^{t+1}\norm{\log\xi^{(0)} - A\best{\nu}/\tau}_\infty + 2\eta \norm{A}_\infty (1-\eta\tau)^{(t+1)/2} \frac{1}{1-(1-\eta\tau)^{1/2}} {\KL{\best{\zeta}}{\zeta^{(0)}}}^{1/2}\\
    &\le (1-\eta\tau)^{t+1}\norm{\log\xi^{(0)} - A\best{\nu}/\tau}_\infty + 4\tau^{-1}\norm{A}_\infty (1-\eta\tau)^{(t+1)/2} {\KL{\best{\zeta}}{\zeta^{(0)}}}^{1/2},
\end{align*}
where the last line results from the fact that $(1-\eta\tau)^{1/2} \le 1 - \eta\tau/2$.   
Combining pieces together, we end up with 
\begin{align*}
    \norm{\log \mu^{(t+1)} - \log \best{\mu}}_\infty &\le 2\norm{\log\xi^{(t+1)} - A\best{\nu}/\tau}_\infty\\
    &\le 2(1-\eta\tau)^{t+1}\norm{\log \xi^{(0)} - A\best{\nu}/\tau}_\infty + 8\tau^{-1}\norm{A}_\infty (1-\eta\tau)^{(t+1)/2} {\KL{\best{\zeta}}{\zeta^{(0)}}}^{1/2}\\
    &\le 2(1-\eta\tau)^{t+1}\norm{\log \mu^{(0)} - \log \best{\mu}}_\infty + 8\tau^{-1}\norm{A}_\infty (1-\eta\tau)^{(t+1)/2} {\KL{\best{\zeta}}{\zeta^{(0)}}}^{1/2}.
\end{align*}
Similarly, one can establish the corresponding inequality for $\norm{\log \nu^{(t+1)} - \log \best{\nu}}_\infty$, therefore completing the proof of inequality~\eqref{eq:EG_conv_log}.

% \begin{proof}
    
% \end{proof}

\subsubsection{Proof of convergence of optimality gap~\eqref{eq:EG_conv_f_last}}
 % and \eqref{eq:EG_conv_f}}
\label{sec:pf-eg_conv_f}
To streamline our discussions, we only provide the proof of inequality~\eqref{eq:EG_conv_f_last} 
% and \eqref{eq:EG_conv_f} 
concerning upper bounding $ f_\tau(\bar{\mu}^{(t)}, \bar{\nu}^{(t)}) - f_\tau(\best{\mu}, \best{\nu})$ without taking the absolute value; the other direction of the inequality can be established in the similar manner and hence is omitted.
%without absolute values on the left hand sides and the other side of both inequalities .

We first make note of an important relation that holds both for PU and OMWU. 
 Consider the update rule of $(\mu^{(t+1)}, \nu^{(t+1)})$, which is the same in PU and OMWU. Lemma~\ref{lemma:mirror_update} inequality~\eqref{eq:ohhh_mu} gives 
\begin{equation}
    \innprod{\log \mu^{(t+1)} - (1-\eta\tau)\log \mu^{(t)} - \eta\tau \log \best{\mu}, \bar{\mu}^{(t+1)} - \best{\mu}} = \eta (\best{\mu} - \bar{\mu}^{(t+1)})^{\top} A ( \best{\nu} - \bar{\nu}^{(t+1)} ).
\end{equation} 
Similar to what we have done in the proof of \eqref{eq:EG_conv_KL} (cf.~\eqref{eq:EGohhh2}), based on the above relation, we can therefore rearrange terms and conclude that
\begin{align}
 \notag   & \eta\prn{\tau\KL{\bar{\mu}^{(t+1)}}{\best{\mu}} -     (\best{\mu} - \bar{\mu}^{(t+1)})^{\top} A ( \best{\nu} - \bar{\nu}^{(t+1)} ) }\\
    &=(1-\eta\tau) \KL{\best{\mu}}{\mu^{(t)}} - (1-\eta\tau)\KL{\bar{\mu}^{(t+1)}}{\mu^{(t)}} - \KL{\mu^{(t+1)}}{\bar{\mu}^{(t+1)}}   \notag\\
    & \qquad + \innprod{\log \bar{\mu}^{(t+1)} - \log \mu^{(t+1)}, \bar{\mu}^{(t+1)} - \mu^{(t+1)}} - \KL{\best{\mu}}{\mu^{(t+1)}}.
    \label{eq:fix_this}
\end{align}

In conjunction with Lemma~\ref{lemma:gap_KL_diff_decomp} (cf.~\eqref{eq:f_diff_decomp}),
we can further derive 
% Along with \eqref{eq:EG_implicit_gap} and \eqref{eqn:piazzolla-decomposition} 
\begin{align}
 &   \eta \big(f_{\tau}(\best{\mu}, \best{\nu})  - f_{\tau}(\bar{\mu}^{(t+1)}, \bar{\nu}^{(t+1)}) \big)
    \le\eta\prn{\tau\KL{\bar{\mu}^{(t+1)}}{\best{\mu}} - ( \best{\mu} - \bar{\mu}^{(t+1)})^{\top} A ( \best{\nu} - \bar{\nu}^{(t+1)}) }  \nonumber \\
    &= (1-\eta\tau) \KL{\best{\mu}}{\mu^{(t)}} - (1-\eta\tau)\KL{\bar{\mu}^{(t+1)}}{\mu^{(t)}} - \KL{\best{\mu}}{\mu^{(t+1)}} \nonumber \\
    &\qquad - \KL{\mu^{(t+1)}}{\bar{\mu}^{(t+1)}}+ \innprod{\log \bar{\mu}^{(t+1)} - \log \mu^{(t+1)}, \bar{\mu}^{(t+1)} - \mu^{(t+1)}}, \label{eqn:piazzolla-decomposition}
%
    % & \quad - (1 - \eta \norm{A}_\infty) \KL{\mu^{(t+1)}}{\bar{\mu}^{(t+1)}} + \eta \norm{A}_\infty \KL{\bar\nu^{(t+1)}}{\nu^{(t)}}\\
%
\end{align}
where the second line follows from \eqref{eq:fix_this}. From this point, we shall continue the proofs for PU and OMWU separately but follow similar strategies.

\paragraph{Remaining steps for PU.} 
Plugging relation~\eqref{eq:EG_implicit_gap} into \eqref{eqn:piazzolla-decomposition},
we arrive at
\begin{align}
    &\eta \big(f_{\tau}(\best{\mu}, \best{\nu}) - f_{\tau}(\bar{\mu}^{(t+1)}, \bar{\nu}^{(t+1)}) \big)\nonumber\\
     &\leq (1-\eta\tau) \KL{\best{\mu}}{\mu^{(t)}} - (1-\eta\tau)\KL{\bar{\mu}^{(t+1)}}{\mu^{(t)}} - \KL{\best{\mu}}{\mu^{(t+1)}} \nonumber \\
&\qquad   - (1 - \eta \norm{A}_\infty) \KL{\mu^{(t+1)}}{\bar{\mu}^{(t+1)}} + \eta \norm{A}_\infty \KL{\bar\nu^{(t+1)}}{\nu^{(t)}} \nonumber \\
    &\le (1-\eta\tau) \KL{\best{\mu}}{\mu^{(t)}} - \KL{\best{\mu}}{\mu^{(t+1)}} - (1-\eta\tau)\KL{\bar{\mu}^{(t+1)}}{\mu^{(t)}} + \eta \norm{A}_\infty \KL{\bar\nu^{(t+1)}}{\nu^{(t)}} ,
    \label{eq:EG_conv_f_step_mu}
\end{align}
where the last line holds since $\eta (\tau+\norm{A}_\infty) \leq 1 $.
Similarly, from Lemma~\ref{lemma:mirror_update} inequality~\eqref{eq:ohhh_nu}, one can establish the following inequality in parallel
\begin{align}
    &\eta \big(f_{\tau}(\bar{\mu}^{(t+1)}, \bar{\nu}^{(t+1)}) - f_{\tau}(\best{\mu}, \best{\nu}) \big)\nonumber\\
    &\le (1-\eta\tau) \KL{\best{\nu}}{\nu^{(t)}} - \KL{\best{\nu}}{\nu^{(t+1)}} - (1-\eta\tau)\KL{\bar{\nu}^{(t+1)}}{\nu^{(t)}} + \eta \norm{A}_\infty \KL{\bar\mu^{(t+1)}}{\mu^{(t)}}.
    \label{eq:EG_conv_f_step_nu}
\end{align}

We are ready to establish inequality~\eqref{eq:EG_conv_f_last} for PU. 
% \paragraph{Proof of inequality~\eqref{eq:EG_conv_f_last}.}
Computing \eqref{eq:EG_conv_f_step_mu} $+ \frac{2}{3} \cdot$ \eqref{eq:EG_conv_f_step_nu} gives
\begin{align}
    & \frac{\eta}{3} \big(f_{\tau}(\best{\mu}, \best{\nu}) - f_{\tau}(\bar{\mu}^{(t+1)}, \bar{\nu}^{(t+1)}) \big)\nonumber\\
    &\le (1-\eta\tau)\brk{\KL{\best{\mu}}{\mu^{(t)}}+\frac{2}{3}\KL{\best{\nu}}{\nu^{(t)}}} - \brk{\KL{\best{\mu}}{\mu^{(t+1)}}+\frac{2}{3}\KL{\best{\nu}}{\nu^{(t+1)}}}\nonumber\\
    &\qquad  - \brk{(1-\eta\tau)  - \frac{2}{3}\eta\norm{A}_\infty}\KL{\bar{\mu}^{(t+1)}}{\mu^{(t)}} + \brk{\eta\norm{A}_\infty - \frac{2}{3}(1-\eta\tau)}\KL{\bar{\nu}^{(t+1)}}{\nu^{(t)}}\nonumber\\
    &\le (1-\eta\tau)\brk{\KL{\best{\mu}}{\mu^{(t)}}+\frac{2}{3}\KL{\best{\nu}}{\nu^{(t)}}} - \brk{\KL{\best{\mu}}{\mu^{(t+1)}}+\frac{2}{3}\KL{\best{\nu}}{\nu^{(t+1)}}}
\label{eq:EG_conv_f_upper}
\end{align}
Here, the last step is due to the fact that $(1-\eta\tau)  - \frac{2}{3}\eta\norm{A}_\infty \ge 0$ and $\eta\norm{A}_\infty - \frac{2}{3}(1-\eta\tau) \le 0$ when $0< \eta \le \frac{1}{\tau + 2\norm{A}_\infty}$. 
As a direct consequence, the difference $f_{\tau}(\best{\mu}, \best{\nu}) - f_{\tau}(\bar{\mu}^{(t)}, \bar{\nu}^{(t)})$ satisfies 
\begin{align*}
  &  \frac{\eta}{3}\big(f_{\tau}(\best{\mu}, \best{\nu}) - f_{\tau}(\bar{\mu}^{(t)}, \bar{\nu}^{(t)}) \big) \\
    &\le (1-\eta\tau)\brk{(1-\eta\tau)\KL{\best{\mu}}{\mu^{(t-1)}}+\eta\norm{A}_\infty\KL{\best{\nu}}{\nu^{(t-1)}}}\\
    &\le (1-\eta\tau) \KL{\best{\zeta}}{\zeta^{(t-1)}} \le (1-\eta\tau)^{t} \KL{\best{\zeta}}{\zeta^{(0)}}.
\end{align*}
We conclude by noting that the other side of \eqref{eq:EG_conv_f_last} can be shown by considering $\frac{2}{3}\cdot$ \eqref{eq:EG_conv_f_step_mu} $+$ \eqref{eq:EG_conv_f_step_nu} combined with similar arguments, and are therefore omitted.

\paragraph{Remaining steps for OMWU.} 
 Similar to the case of PU, plugging \eqref{eq:OMWU_implicit_gap}
into \eqref{eqn:piazzolla-decomposition} gives
\begin{align}
    &\eta \big(f_{\tau}(\best{\mu}, \best{\nu}) - f_{\tau}(\bar{\mu}^{(t+1)}, \bar{\nu}^{(t+1)}) \big)
    \nonumber\\
    &\leq (1-\eta\tau) \KL{\best{\mu}}{\mu^{(t)}} - (1-\eta\tau)\KL{\bar{\mu}^{(t+1)}}{\mu^{(t)}} - \KL{\best{\mu}}{\mu^{(t+1)}} - (1-2\eta\norm{A}_\infty)\KL{\mu^{(t+1)}}{\bar{\mu}^{(t+1)}}\nonumber\\
    & \qquad\qquad  + \eta \norm{A}_\infty \left[ \KL{\nu^{(t)}}{\bar{\nu}^{(t)}} + \KL{\bar\nu^{(t+1)}}{\nu^{(t)}}\right].
    \label{eq:OMWU_conv_f_step_mu}
\end{align}
Similarly, one can establish a symmetric inequality as follows  
\begin{align}
    &\eta \big( f_{\tau}(\bar{\mu}^{(t+1)}, \bar{\nu}^{(t+1)}) - f_{\tau}(\best{\mu}, \best{\nu}) \big)
    \nonumber\\
    &\leq (1-\eta\tau) \KL{\best{\nu}}{\nu^{(t)}} - (1-\eta\tau)\KL{\bar{\nu}^{(t+1)}}{\nu^{(t)}} - \KL{\best{\nu}}{\nu^{(t+1)}} - (1-2\eta\norm{A}_\infty)\KL{\nu^{(t+1)}}{\bar{\nu}^{(t+1)}}\nonumber\\
    & \quad  + \eta \norm{A}_\infty \left[ \KL{\mu^{(t)}}{\bar{\mu}^{(t)}} + \KL{\bar\mu^{(t+1)}}{\mu^{(t)}}\right].
    \label{eq:OMWU_conv_f_step_nu}
\end{align}
% \paragraph{Proof of inequality~\eqref{eq:EG_conv_f_last}.}
 
Directly computing \eqref{eq:OMWU_conv_f_step_mu} $+ \frac{2}{3} \cdot$ \eqref{eq:OMWU_conv_f_step_nu} gives
\begin{align}
    \label{enq:OMWU-finale}
    &\frac{\eta}{3}\cdot(f_{\tau}(\best{\mu}, \best{\nu}) - f_{\tau}(\bar{\mu}^{(t+1)}, \bar{\nu}^{(t+1)}))\nonumber\\
    &\le (1-\eta\tau)\brk{\KL{\best{\mu}}{\mu^{(t)}}+\frac{2}{3}\KL{\best{\nu}}{\nu^{(t)}}} - \brk{\KL{\best{\mu}}{\mu^{(t+1)}}+\frac{2}{3}\KL{\best{\nu}}{\nu^{(t+1)}}}\nonumber\\
    &\quad - \brk{(1-\eta\tau)  - \frac{2}{3}\eta\norm{A}_\infty}\KL{\bar{\mu}^{(t+1)}}{\mu^{(t)}} + \brk{\eta\norm{A}_\infty - \frac{2}{3}(1-\eta\tau)}\KL{\bar{\nu}^{(t+1)}}{\nu^{(t)}}\nonumber\\
    &\quad+ \eta\norm{A}_\infty\left[\frac{2}{3}\KL{\mu^{(t)}}{\bar{\mu}^{(t)}}+ \KL{\nu^{(t)}}{\bar{\nu}^{(t)}}\right]\nonumber\\
    &\quad- (1-2\eta\norm{A}_\infty)\left[\KL{\mu^{(t+1)}}{\bar{\mu}^{(t+1)}}+\frac{2}{3}\KL{\nu^{(t+1)}}{\bar{\nu}^{(t+1)}}\right] .
\end{align}
With our choice of the learning rate $\eta \le \min\{\frac{1}{2\norm{A}_\infty + 2\tau}, \frac{1}{4\norm{A}_\infty}\}$, it is guarantees that
\[
\eta\norm{A}_\infty - \frac{2}{3}(1-\eta\tau) \leq 0,\quad (1-\eta\tau)  - \frac{2}{3}\eta\norm{A}_\infty \geq 0 \quad \mbox{and} \quad    (1-\eta\tau)(1 - 2\eta \norm{A}_\infty) \ge \frac{3}{2}\eta \norm{A}_\infty.
\]
To proceed, let us introduce the shorthand notation 
\begin{align*}
    G^{(t)} &:= \KL{\best{\mu}}{\mu^{(t)}}+\frac{2}{3}\KL{\best{\nu}}{\nu^{(t)}} \\
    &\qquad+ \frac{2}{3}(1-2\eta\norm{A}_\infty)\left[\KL{\mu^{(t)}}{\bar{\mu}^{(t)}}+\KL{\nu^{(t)}}{\bar{\nu}^{(t)}}\right].
\end{align*}
With this piece of notation, we can write inequality~\eqref{enq:OMWU-finale} as
\begin{equation}
    \frac{\eta}{3}(f_{\tau}(\best{\mu}, \best{\nu}) - f_{\tau}(\bar{\mu}^{(t+1)}, \bar{\nu}^{(t+1)})) \le (1-\eta\tau) G^{(t)} - G^{(t+1)},
    \label{eq:OMWU_conv_f_upper}
\end{equation}
which in turn implies 
\begin{align*}
    &\frac{\eta}{3}(f_{\tau}(\best{\mu}, \best{\nu}) - f_{\tau}(\bar{\mu}^{(t)}, \bar{\nu}^{(t)})) \\
    &\le (1-\eta\tau) G^{(t-1)} \le (1-\eta\tau)L^{(t-1)} \le (1-\eta\tau)^tL^{(0)} = (1-\eta\tau)^t \KL{\best{\zeta}}{\zeta^{(0)}},
\end{align*}
with $L^{(t)}$ defined in \eqref{eq:def_Lt}. This finishes the proof of \eqref{eq:EG_conv_f_last} for OMWU.
 
\subsubsection{Proof of convergence of duality gap~\eqref{eq:EG_conv_gap_last} }
% {\color{dacong}

The proof of inequality~\eqref{eq:EG_conv_gap_last} is built upon the following lemma whose proof is deferred to Appendix~\ref{sec:pf_dual_gap}.
\begin{lemma}
    \label{lemma:dual_gap}
    The duality gap at $\zeta = (\mu, \nu)$ can be bounded as
    \begin{align*}
     \max_{\mu'\in \Delta(\mathcal{A})} f_\tau(\mu', \nu) - \min_{\nu'\in\Delta(\mathcal{B})}f_\tau(\mu, \nu')   \le \tau \KL{\zeta}{\best{\zeta}} + \tau^{-1}{\norm{A}_\infty^2} \KL{\best{\zeta}}{\zeta}. 
    \end{align*}    
    % or
    % \begin{align*}
    %     \max_{\mu'\in \Delta(\mathcal{A})} f_\tau(\mu', \nu) - \min_{\nu'\in\Delta(\mathcal{B})}f_\tau(\mu, \nu') \le \tau\KL{\zeta}{\best{\zeta}} + 2\norm{A}_\infty \KL{\best{\zeta}}{\zeta} + 4\norm{A}_\infty \KL{\best{\zeta}}{\zeta}^{1/2}.
    % \end{align*}    TBD
\end{lemma}
% }

Applying Lemma \ref{lemma:dual_gap} to $\bar{\zeta}^{(t)} = (\bar{\mu}^{(t)}, \bar{\nu}^{(t)})$ yields
\begin{align}
    \mathsf{DualGap}_{\tau} (\bar{\zeta}^{(t)}) 
    &\le \tau\KL{\bar{\zeta}^{(t)}}{\best{\zeta}} + \tau^{-1}{\norm{A}_\infty^2} \KL{\best{\zeta}}{\bar{\zeta}^{(t)}}\nonumber\\
    &\le\tau\KL{\bar{\zeta}^{(t)}}{\best{\zeta}} + 2\tau^{-1}{\norm{A}_\infty^2} (1-\eta\tau)^{t-1} \KL{\best{\zeta}}{{\zeta}^{(0)}},
    \label{eq:EG_conv_gap_step1}
\end{align}
where the second step results from \eqref{eq:EG_conv_KL}.
It remains to bound $\tau\KL{\bar{\zeta}^{(t)}}{\best{\zeta}}$, which we proceed separately for PU and OMWU. 

\paragraph{Remaining steps for PU.} From inequality \eqref{eq:EGohhh_final}, we are ensured that 
\[
    \eta\tau\KL{\bar{\zeta}^{(t)}}{\best{\zeta}} \le (1-\eta\tau) \KL{\best{\zeta}}{\zeta^{(t-1)}} - \KL{\best{\zeta}}{\zeta^{(t)}}.
    % \label{eq:EG_conv_gap_step1}
\]
It thus follows that
\begin{align*}
  \tau\KL{\bar{\zeta}^{(t)}}{\best{\zeta}}\le \eta^{-1}(1-\eta\tau) \KL{\best{\zeta}}{\zeta^{(t-1)}}\le \eta^{-1}(1-\eta\tau)^{t-1}\KL{\best{\zeta}}{\zeta^{(0)}}, 
\end{align*}
where the last inequality is due to inequality~\eqref{eq:EG_conv_KL}.
Plugging the above inequality into \eqref{eq:EG_conv_gap_step1} completes the proof of inequality~\eqref{eq:EG_conv_gap_last} for PU. 
% which completes the proof of inequality~\eqref{eq:EG_conv_gap_last}. 

\paragraph{Remaining steps for OMWU.} 
From inequality \eqref{eq:OMWUohhh_final}, we are ensured that
\begin{align*}
    \tau\KL{\bar{\zeta}^{(t)}}{\best{\zeta}}\le \eta^{-1}(1-\eta\tau) L^{(t-1)}
    &\le \eta^{-1}(1-\eta\tau)^tL^{(0)} = \eta^{-1}(1-\eta\tau)^t\KL{\best{\zeta}}{\zeta^{(0)}},
\end{align*}
where the last equality follows from $L^{(0)} = \KL{\best{\zeta}}{\zeta^{(0)}}$. Plugging the above inequality into \eqref{eq:EG_conv_gap_step1} finishes the proof of inequality~\eqref{eq:EG_conv_gap_last} for OMWU.

% Next, we move on to consider the average convergence. 
% First, notice that $\eta\tau \sum_{i=0}^{t-1} (1-\eta\tau)^i = 1- (1-\eta\tau)^{t}$, we have
% \begin{align}
%     \big(1- (1-\eta\tau)^{t} \big)\ex{k}{\KL{\bar{\zeta}^{(k)}}{\best{\zeta}}}=\eta\tau \sum_{i=0}^{t-1} (1-\eta\tau)^i \KL{\bar{\zeta}^{(t-i)}}{\best{\zeta}},  
%     % \le (1-\eta\tau)^{t} \KL{\best{\zeta}}{\zeta^{(0)}},
%     \label{eq:EG_conv_gap_step2}
% \end{align}
% where we expand the expectation w.r.t. the distribution $p(k=i) \propto (1-\eta\tau)^{t-i}$, $0\leq i\leq t$. 
% To further control the right-hand side, as a result of inequality~\eqref{eq:EG_conv_gap_step1} plus some direct calculations, one obtains 
% \begin{align*}
%    \tau \sum_{i=0}^{t-1} \eta(1-\eta\tau)^i \KL{\bar{\zeta}^{(t-i)}}{\best{\zeta}} 
%   %
%   &\le (1-\eta\tau)^{t} \KL{\best{\zeta}}{\zeta^{(0)}} - \KL{\best{\zeta}}{\zeta^{(t)}}.
% \end{align*}
% Putting things together gives 
% \begin{align*}
%     \tau^{-1} \big(1- (1-\eta\tau)^{t} \big)\ex{k}{f_{\tau}(\best{\mu}, \bar{\nu}^{(k)}) - f_{\tau}(\bar{\mu}^{(k)}, \best{\nu})}
%     \leq 
%     (1-\eta\tau)^{t} \KL{\best{\zeta}}{\zeta^{(0)}}.
% \end{align*}
% Reorganizing terms finishes the proof of inequality~\eqref{eq:EG_conv_gap}. 

% where the last step uses inequality~\eqref{eq:EG_conv_gap_step1}. 

	 \subsection{Proof of Theorem~\ref{thm:no_regret}}
\label{sec:pf-regret}

To begin, we note that in the no-regret setting, upon receiving $A^\top \bar{\mu}^{(t)}$ (which  is possibly adversarial), the update rule of player 2 is given by
            \begin{subequations} \label{eq:single_omwu}
            \begin{align}
                    {\nu}^{(t)}(b) & \propto {\nu^{(t-1)}(b)}^{1-\eta_{t-1}\tau}\exp(-\eta_{t-1} [A^\top \bar{\mu}^{(t)}]_b),  \label{eq:single_omwu_nu} \\
                    \bar{\nu}^{(t+1)}(b) & \propto {\nu^{(t)}(b)}^{1-\eta_{t}\tau}\exp(-\eta_{t} [A^\top \bar{\mu}^{(t)}]_b). \label{eq:single_omwu_bar_nu}
            \end{align}
            \end{subequations}

Recalling $f_\tau^{(t)}(\nu) = \bar{\mu}^{(t)\top} A \nu  + \tau \mathcal{H}(\bar{\mu}^{(t)}) - \tau \mathcal{H}(\nu)$, we introduce an important quantity which is the gradient of $f_\tau^{(t)}(\nu) $ at $\bar\nu^{(t)}$:
\begin{equation}\label{eq:nabla_t}
\overline{\nabla}^{(t)} := \nabla_\nu f_\tau^{(t)}(\nu)\Big|_{\nu=\bar\nu^{(t)}} = A^\top \bar{\mu}^{(t)} + \tau(\log \bar\nu^{(t)} + \mathbf{1}).
\end{equation}
The following lemma presents an $\ell_\infty$ bound on the size of $\overline\nabla^{(t)}$, whose proof can be found in Appendix~\ref{proof:lemma_OMWU_grad_bound}.
\begin{lemma}
    \label{lemma:OMWU_grad_bound}
    It holds for all $t \ge 0$ that
    \[
        \left\| \overline\nabla^{(t)} \right\|_\infty \le \tau\log |\mathcal{B}| + 3\norm{A}_\infty.
    \]
\end{lemma}

\paragraph{Regret decomposition.}
By the definition of $f_\tau^{(t)}(\nu) $, we have
\begin{align}
    f_\tau^{(t)}(\bar{\nu}^{(t)}) - f_\tau^{(t)}(\nu) &= \langle{A^\top \bar{\mu}^{(t)}}, \bar{\nu}^{(t)} - \nu\rangle + \tau \bar{\nu}^{(t)\top} \log {\bar\nu^{(t)}} - \tau\nu^\top \log \nu\notag\\
    &= (A^\top \bar{\mu}^{(t)} + \tau\log \bar\nu^{(t)})^\top(\bar\nu^{(t)} - \nu) - \tau\KL{\nu}{\bar\nu^{(t)}}\notag\\
    &= \langle \overline\nabla^{(t)}, \bar\nu^{(t)} - \nu\rangle - \tau\KL{\nu}{\bar\nu^{(t)}},\label{eq:single_regret_decomp}
\end{align}
where the last line follows from the definition of $\overline\nabla^{(t)}$ (cf. \eqref{eq:nabla_t}).
To continue, by the update rule in \eqref{eq:single_omwu}, we have 
\begin{align}
    \log \bar{\nu}^{(t+1)} &= (1-\eta_t \tau) \log \nu^{(t)} - \eta_t A^\top \bar\mu^{(t)} + c_1\cdot \mathbf{1}\notag\\
    &= (1-\eta_t\tau) \left[(1-\eta_{t-1}\tau) \log \nu^{(t-1)} - \eta_{t-1} A^\top \bar\mu^{(t)} + c_2 \cdot \mathbf{1}\right] - \eta_t A^\top \bar\mu^{(t)} + c_1\cdot \mathbf{1}\notag\\
    &= (1-\eta_t\tau) \left[ \log \bar{\nu}^{(t)} + \eta_{t-1} A^\top \bar\mu^{(t-1)}- \eta_{t-1} A^\top \bar\mu^{(t)} \right] - \eta_t A^\top \bar\mu^{(t)} + c^{(t)}\cdot\mathbf{1}\notag\\
    &= \log \bar{\nu}^{(t)} - \eta_t \overline{\nabla}^{(t)} + (1-\eta_t\tau)\eta_{t-1}A^\top (\bar\mu^{(t-1)} - \bar\mu^{(t)}) + c^{(t)}\cdot\mathbf{1},\label{eq:OMWU_log_update}
\end{align}
where the first three steps result from the update rule of $\bar{\nu}^{(t+1)}$, $\nu^{(t)}$ and $\bar{\nu}^{(t)}$, respectively, and the last line follows from \eqref{eq:nabla_t}. Here, $c_1$, $c_2$, and $c^{(t)}$ are some normalization constants.\footnote{We shall set $\eta_{-1} = 0$ and $\bar{\mu}^{(-1)} = \bar\mu^{(0)}$ to accommodate the case when $t=0$ in \eqref{eq:OMWU_log_update}.} Rearranging terms allows us to rewrite $\overline{\nabla}^{(t)} $ as
\begin{equation} \label{eq:grad_OMWU}
    \overline{\nabla}^{(t)} = \frac{1}{\eta_t} \left(\log \bar{\nu}^{(t)} - \log \bar{\nu}^{(t+1)} + c^{(t)}\cdot \mathbf{1}\right)+\frac{\eta_{t-1}}{\eta_t}(1-\eta_t\tau)  A^\top(\bar\mu^{(t-1)}-\bar\mu^{(t)}).
\end{equation}
Plugging \eqref{eq:grad_OMWU} into \eqref{eq:single_regret_decomp}, we have
\begin{align*}
   f_\tau^{(t)}(\bar{\nu}^{(t)})  &- f_\tau^{(t)}(\nu) \\
 & = \frac{1}{\eta_t}\langle \log\bar{\nu}^{(t)} - \log\bar{\nu}^{(t+1)}, \bar{\nu}^{(t)} - \nu\rangle - \tau\KL{\nu}{\bar{\nu}^{(t)}}  +\frac{\eta_{t-1}}{\eta_t}(1-\eta_t\tau) \langle A^\top(\bar\mu^{(t-1)}-\bar\mu^{(t)}), \bar{\nu}^{(t)} - \nu\rangle\\
    &= \frac{1}{\eta_t}\left[ \KL{\nu}{\bar{\nu}^{(t)}} - \KL{\nu}{\bar{\nu}^{(t+1)}} + \KL{\bar{\nu}^{(t)}}{\bar{\nu}^{(t+1)}}\right] - \tau\KL{\nu}{\bar{\nu}^{(t)}}\\
    &\qquad +\frac{\eta_{t-1}}{\eta_t}(1-\eta_t\tau) \langle A^\top(\bar\mu^{(t-1)}-\bar\mu^{(t)}), \bar{\nu}^{(t)} - \nu\rangle.
\end{align*}
Summing the equality over $t=0,\ldots, T$ gives
% \begin{align*}
%     \mathsf{Regret}(T) &= \sum_{t=0}^T f_\tau^{(t)}(\nu^{(t)}) - \min_{\nu \in \Delta(\mathcal{B})}\sum_{t=0}^T f_\tau^{(t)}(\nu)\\
%     &= \max_{\nu \in \Delta(\mathcal{B})}\Big\{\sum_{t=0}^T f_\tau^{(t)}(\nu^{(t)}) - \sum_{t=0}^T f_\tau^{(t)}(\nu)\Big\}\\
%     &= \max_{\nu \in \Delta(\mathcal{B})}\left\{ \left(\frac{1}{\eta^{(0)}} - \tau\right) \KL{\nu}{\nu^{(0)}} +  \sum_{t=1}^T \left(\frac{1}{\eta_t} - \frac{1}{\eta_{t-1}}-\tau\right)\KL{\nu}{\nu^{(t)}}\right\} \\
%     &\quad + \sum_{t=0}^T \frac{1}{\eta_t}\KL{\nu^{(t)}}{\nu^{(t+1)}}.
% \end{align*}
\begin{align}
  \sum_{t=0}^T \left[f_\tau^{(t)}(\bar{\nu}^{(t)}) - f_\tau^{(t)}(\nu)\right]    & = \left(\frac{1}{\eta_{0}} - \tau\right) \KL{\nu}{\bar\nu^{(0)}} +  \sum_{t=1}^T \left(\frac{1}{\eta_t} - \frac{1}{\eta_{t-1}}-\tau\right)\KL{\nu}{\bar\nu^{(t)}} \nonumber\\
    &\quad + \sum_{t=0}^T \frac{1}{\eta_t}\KL{\bar\nu^{(t)}}{\bar\nu^{(t+1)}} + \sum_{t=1}^{T} \frac{\eta_{t-1}}{\eta_t}(1-\eta_t\tau) \langle A^\top(\bar\mu^{(t-1)}-\bar\mu^{(t)}), \bar{\nu}^{(t)} - \nu\rangle.\label{eq:regret_decomp_pre}
\end{align}
With the choice of the learning rate 
$$\eta_t = \frac{1}{(t+1)\tau},$$ 
one has
\begin{equation}
 \frac{1}{\eta_{0}} - \tau = 0 ,\qquad  \frac{1}{\eta_t} - \frac{1}{\eta_{t-1}} - \tau = 0, \qquad \frac{\eta_{t-1}}{\eta_t}(1-\eta_t\tau) = 1,\; \forall\;t\ge 1.
    \label{eq:magic_etas}
\end{equation}
Plugging the above relations into \eqref{eq:regret_decomp_pre} leads to
\begin{align}
  \sum_{t=0}^T \left[f_\tau^{(t)}(\bar{\nu}^{(t)}) - f_\tau^{(t)}(\nu)\right] &= \sum_{t=0}^T \frac{1}{\eta_t}\KL{\bar\nu^{(t)}}{\bar\nu^{(t+1)}} + \sum_{t=1}^T \langle A^\top(\bar\mu^{(t-1)}-\bar\mu^{(t)}), \bar{\nu}^{(t)} - \nu\rangle.\label{eq:regret_decomp}
\end{align}
Next, we seek to bound the two terms in \eqref{eq:regret_decomp} separately. 

\paragraph{Bounding the first term of the regret.}
According to Lemma \ref{lem:log_pi_gap} and \eqref{eq:OMWU_log_update}, we have
\begin{align*}
    \KL{\bar\nu^{(t)}}{\bar\nu^{(t+1)}} &\le \frac{1}{2}\norm{\log \bar\nu^{(t+1)} - \log \bar\nu^{(t)} - c^{(t)}\cdot \mathbf{1}}_\infty^2\\
    &= \frac{1}{2}\norm{- \eta_t \overline{\nabla}^{(t)} + (1-\eta_t\tau)\eta_{t-1}A^\top (\bar\mu^{(t-1)} - \bar\mu^{(t)})}_\infty^2\\
    &\le \norm{-\eta_t\overline\nabla^{(t)}}_\infty^2 + \norm{(1-\eta_t\tau)\eta_{t-1}A^\top (\bar\mu^{(t-1)} - \bar\mu^{(t)})}_\infty^2\\
    &\le {\eta_t}^2 (\tau\log |\mathcal{B}| + 3\norm{A}_\infty)^2 + 4{\eta_t}^2\norm{A}_\infty^2,
\end{align*}
where the final step results from Lemma \ref{lemma:OMWU_grad_bound} and the last equality in \eqref{eq:magic_etas}. Summing the above inequality over $t = 0,\cdots, T$ yields
\begin{align}
    \sum_{t=0}^T \frac{1}{\eta_t}\KL{\bar\nu^{(t)}}{\bar\nu^{(t+1)}} &\le \sum_{t=0}^T \eta_t (\tau\log |\mathcal{B}| + 3\norm{A}_\infty)^2 + 4\sum_{t=0}^T {\eta_t}\norm{A}_\infty^2\notag\\
    &\le \tau^{-1}(\log T + 1)\left[(\tau\log |\mathcal{B}| + 3\norm{A}_\infty)^2 + 4\norm{A}_\infty^2\right],\label{eq:regret_part1}
\end{align}
where the last line follows from
$ \sum_{t=0}^T \eta_t \leq\tau^{-1}(\log T + 1) $
due to the choice of the learning rate.

\paragraph{Bounding the second term of the regret.} Observe that by the telescoping relation, we have
\begin{align}
    &\sum_{t=1}^T \langle A^\top(\bar\mu^{(t-1)}-\bar\mu^{(t)}), \bar{\nu}^{(t)} - \nu\rangle \notag\\
    &= \sum_{t=1}^T \langle A^\top(\bar\mu^{(t-1)}-\bar\mu^{(t)}), \bar{\nu}^{(t)}\rangle - \langle A^\top(\mu^{(0)} - \bar\mu^{(T)}), \nu\rangle \notag\\
%    &= \sum_{t=1}^T \left[\langle A^\top\bar\mu^{(t-1)}, \bar{\nu}^{(t)}\rangle - \langle A^\top \bar\mu^{(t)}, \bar{\nu}^{(t)}\rangle \right] - \langle A^\top(\mu^{(0)} - \bar\mu^{(T)}), \nu\rangle \notag\\
    &= \sum_{t=1}^T \langle A^\top \bar\mu^{(t-1)}, \bar{\nu}^{(t)} - \bar{\nu}^{(t-1)} \rangle + \sum_{t=1}^T \left[\langle A^\top\bar\mu^{(t-1)}, \bar{\nu}^{(t-1)}\rangle - \langle A^\top \bar\mu^{(t)}, \bar{\nu}^{(t)}\rangle \right] - \langle A^\top(\mu^{(0)} - \bar\mu^{(T)}), \nu\rangle \notag\\
    &= \sum_{t=1}^T \langle A^\top \bar\mu^{(t-1)}, \bar{\nu}^{(t)} - \bar{\nu}^{(t-1)} \rangle + \langle A^\top\mu^{(0)}, \bar{\nu}^{(0)}\rangle - \langle A^\top \bar\mu^{(T)}, \bar{\nu}^{(T)}\rangle  - \langle A^\top(\mu^{(0)} - \bar\mu^{(T)}), \nu\rangle \notag\\
    &\le \norm{A}_\infty \sum_{t=0}^{T-1} \norm{\bar{\nu}^{(t+1)} - \bar{\nu}^{(t)}}_1 + 4\norm{A}_\infty. \label{eq:regret_part2_tmp}
\end{align}
Due to Pinsker's inequality and Lemma \ref{lem:log_pi_gap}, we have 
\[
    \frac{1}{2}\norm{\bar{\nu}^{(t+1)} - \bar{\nu}^{(t)}}_1^2 \le \KL{\bar{\nu}^{(t+1)}}{\bar{\nu}^{(t)}} \le \frac{1}{2} \norm{\log\bar{\nu}^{(t+1)} - \log\bar{\nu}^{(t)} - c^{(t)}\cdot \mathbf{1}}_{\infty}^2,
\]
which further ensures that
\begin{align*}
    \norm{\bar{\nu}^{(t+1)} - \bar{\nu}^{(t)}}_1 &\le \norm{\log\bar{\nu}^{(t+1)} - \log\bar{\nu}^{(t)} - c^{(t)}\cdot \mathbf{1}}_{\infty}\\
    & = \norm{- \eta_t \overline{\nabla}^{(t)} + (1-\eta_t\tau)\eta_{t-1}A^\top (\bar\mu^{(t-1)} - \bar\mu^{(t)})}_\infty\\
    &\le \eta_t \norm{\overline{\nabla}^{(t)}}_{\infty} + \eta_t \norm{A^\top (\bar\mu^{(t-1)} - \bar\mu^{(t)})}_\infty\\
    &\le \eta_t \left(\tau\log |\mathcal{B}| + 5\norm{A}_\infty\right),
\end{align*}
where the second line follows from \eqref{eq:OMWU_log_update}, the third line follows from the triangle inequality, and the last line follows from Lemma \ref{lemma:OMWU_grad_bound}.
Plugging the above inequality into \eqref{eq:regret_part2_tmp} leads to
\begin{align}
\sum_{t=1}^T \langle A^\top(\bar\mu^{(t-1)}-\bar\mu^{(t)}), \bar{\nu}^{(t)} - \nu\rangle   &\le \tau^{-1}(\log T + 1)\norm{A}_\infty\left(\tau\log |\mathcal{B}| + 5\norm{A}_\infty\right) + 4\norm{A}_\infty,\label{eq:regret_part2}
\end{align}
where we use again $ \sum_{t=0}^T \eta_t \leq\tau^{-1}(\log T + 1) $.

\paragraph{Putting things together.} Combining \eqref{eq:regret_part1} and \eqref{eq:regret_part2} into \eqref{eq:regret_decomp}, we have
\begin{align*}
    \mathsf{Regret}_\lambda(T) &= \max_{\nu \in \Delta(\mathcal{B})} \sum_{t=0}^T \left[ f_\tau^{(t)}(\nu^{(t)}) -  f_\tau^{(t)}(\nu) \right] \le \tau^{-1}(\log T + 1)(\tau\log |\mathcal{B}| + 5\norm{A}_\infty)^2 + 4\norm{A}_\infty.
\end{align*}

% First, we show that
% \begin{align}
%     f_\tau^{(t)}(\nu^{(t)}) - f_\tau^{(t)}(\nu) &= {\mu^{(t)}}^\top A(\bar{\nu}^{(t)} - \nu) + \tau{\nu^{(t)}}^\top \log {\nu^{(t)}} - \tau\nu^\top \log \nu\notag\\
%     &= (A^\top \mu^{(t)} + \tau\log \nu^{(t)})^\top(\nu^{(t)} - \nu) - \tau\KL{\nu}{\nu^{(t)}}\notag\\
%     &= \langle \nabla^{(t)}, \nu^{(t)} - \nu\rangle - \tau\KL{\nu}{\nu^{(t)}},\label{eq:single_regret_decomp}
% \end{align}
% where $\nabla^{(t)} := \nabla_\nu f_\tau^{(t)}(\nu)\Big|_{\nu=\nu^{(t)}} = A^\top \mu^{(t)} + \tau\log \nu^{(t)}$. According to the update rule of MWU, we have
% \[
%     \eta_t \nabla^{(t)} = \log\nu^{(t)} - \log\nu^{(t+1)} + c\cdot\mathbf{1}
% \]
% for some constant $c$. Therefore, we can write 
% \begin{align*}
%     f_\tau^{(t)}(\nu^{(t)}) - f_\tau^{(t)}(\nu) &= \langle \nabla^{(t)}, \nu^{(t)} - \nu\rangle - \tau\KL{\nu}{\nu^{(t)}}\\
%     &= \frac{1}{\eta_t}\langle \log\nu^{(t)} - \log\nu^{(t+1)}, \nu^{(t)} - \nu\rangle - \tau\KL{\nu}{\nu^{(t)}}\\
%     &= \frac{1}{\eta_t}\left(\KL{\nu}{\nu^{(t)}} - \KL{\nu}{\nu^{(t+1)}} + \KL{\nu^{(t)}}{\nu^{(t+1)}}\right) - \tau\KL{\nu}{\nu^{(t)}}.\\
% \end{align*}

	 \section{Analysis for entropy-regularized Markov games}
\label{sec:markov-games-proof-main}

\subsection{Proof of Proposition~\ref{prop:MG-cont}}

For each $t$, let 
\begin{align*}
   V^{(t)}(s) :=   \max_{\mu(s) \in \Delta(\asp)}\min_{\nu(s) \in \Delta(\cB)} 
   \; f_{\tau}(Q^{(t)}(s); \mu(s),\nu(s)),
\end{align*}
which is, in other words, the minimax value of the associated matrix game using a payoff matrix $Q^{(t)}(s)$. 
We start by making a simple observation that for $\mu(s) \in \Delta(\asp), \nu(s) \in \Delta(\cB)$, 
    \begin{align*}
        \left|f_{Q^{(t)}(s)}(\mu(s), \nu(s)) - f_{Q_\tau^{\star}(s)}(\mu(s), \nu(s)) \right| &= \left |\mu(s)^\top (Q^{(t)}(s) - Q_\tau^\star(s)) \nu(s) \right| \\
        & \le \norm{Q^{(t)}(s) - Q_\tau^\star(s)}_\infty \le \norm{Q^{(t)} - Q_\tau^\star}_\infty.
\end{align*}
As a direct consequence, we can control $V^{(t)}(s) - V_\tau^{\star}(s)$ by 
    \begin{align*}
        \left|V^{(t)}(s) - V_\tau^{\star}(s) \right| &= \left|\max_{\mu(s) \in \Delta(\asp)}\min_{\nu(s) \in \Delta(\cB)} f_{Q^{(t)}(s)}(\mu(s), \nu(s)) - \max_{\mu(s) \in \Delta(\asp)}\min_{\nu(s) \in \Delta(\cB)}f_{Q_\tau^{\star}(s)}(\mu(s), \nu(s)) \right|\\
        &\le \norm{Q^{(t)} - Q_\tau^\star}_\infty.
    \end{align*}
Recalling the definition of the soft Bellman operator $\mathcal{T}_{\tau}$ in \eqref{eq:soft_bellman_op}, it then follows that
    \begin{align*}
        \norm{Q^{(t+1)} - Q_\tau^\star}_\infty 
        = \norm{ \mathcal{T}_{\tau}(Q^{(t)}) - \mathcal{T}_{\tau}(Q_{\tau}^{\star})} 
        &= \gamma \left|\mathbb{E}_{s'\sim P(\cdot|s,a,b)} \Big[V^{(t)}(s') - V_\tau^\star(s') \Big]\right|\\
        &\le \gamma \norm{V^{(t)} - V_\tau^\star}_\infty  \le \gamma  \norm{Q^{(t)} - Q_\tau^\star}_\infty.
    \end{align*}
Recursively invoking the above inequality proves inequality~\eqref{eqn:MG-cont}.

\subsection{Proof of Theorem~\ref{Thm:MG-result}}
\label{Sec:pf-mg-result}

% \ytw{use $\epsilon_{\rm opt}$ or simply $\epsilon$?}
The inner loop of Algorithm \ref{alg:approx_value_iter} aims to solve an entropy-regularized matrix game indexed by $Q^{(t)}(s)$, which is done by running the proposed PU or OMWU methods. 
To analyze the efficacy of the inner loop, let us denote the exact minimax game value on state $s$ at $t$-th iteration by 
\begin{align}
\label{eqn:inner-matrix-game}
  \breve{V}^{(t+1)}(s) := \max_{\mu(s) \in \Delta(\asp)}\min_{\nu(s) \in \Delta(\asp)} f_\tau(Q^{(t)}(s); \mu(s), \nu(s)),  
\end{align}
which is adopted in the exact value iteration analyzed in Proposition~\ref{prop:MG-cont}, and achieved by the equilibrium $\zeta^{\star(t)} = (\mu^{\star(t)}, \nu^{\star(t)})$ of \eqref{eqn:inner-matrix-game}.

\begin{itemize} 
\item Denote the output of the inner loop as $\bar{\zeta}^{(t, T_{\text{sub}})} = (\bar\mu^{(t, T_{\text{sub}})}(s), \bar\nu^{(t, T_{\text{sub}})}(s) )$, which the entropy-regularized matrix game~\eqref{eqn:inner-matrix-game} is approximately solved by executing PU\,/\,OMWU for $T_{\text{sub}}$ iterations. 
Theorem~\ref{thm:EG-gurantees-last-iterate} (cf.~\eqref{eq:EG_conv_f_last}) guarantees that for every $s\in\ssp$, one has 
\begin{align*}
  \left| V^{(t+1)}(s) - \breve{V}^{(t+1)}(s) \right|&= \left| f^\tau_{Q^{(t)}(s)}(\bar\mu^{(t, T_{\text{sub}})}(s), \bar\nu^{(t, T_{\text{sub}})}(s)) - f^\tau_{Q^{(t)}(s)}(\mu^{\star(t)}(s), \nu^{\star(t)}(s)) \right|\\
    &\le \eta^{-1}\frac{1}{1 - (\tau + \| Q^{(t)}(s)\|_\infty)\eta} \cdot \frac{(1-\eta\tau)^{T_{\text{sub}}}}{1- (1-\eta\tau)^{T_{\text{sub}}}} \cdot \KL{\zeta^{\star(t)}}{\zeta^{(0)}}\\
    &\le 2\eta^{-1} \cdot \frac{(1-\eta\tau)^{T_{\text{sub}}}}{1- (1-\eta\tau)^{T_{\text{sub}}}} \cdot 2 \log|\mathcal{A}|,
\end{align*}
where the last step makes use of the choice of the learning rate
\[
    \eta = \frac{1-\gamma}{2(1+\tau(\log|\mathcal{A}| + 1))} \le \frac{1}{2(\tau + \| Q^{(t)}(s)\|_\infty)}
\]
and $\KL{\zeta^{\star(t)}}{\zeta^{(0)}}\leq \log|\mathcal{A}| + \log|\mathcal{B}| \leq 2\log|\mathcal{A}|$. As a consequence, setting 
\begin{equation}\label{eq:T_sub_bound}
T_{\text{sub}} = O\prn{\frac{1}{\eta\tau}\prn{\log \frac{1}{\epsilon} + \log\frac{1}{1-\gamma} + \log\log|\asp| + \log \frac{1}{\eta}}}
\end{equation} 
yields 
\begin{equation}\label{eq:V_inner_bd}
    {\left|V^{(t+1)}(s) - \breve{V}^{(t+1)}(s) \right|} \le (1-\gamma)\epsilon, \quad\mbox{for all}\quad s\in\cS.
\end{equation}

\item We now move to monitor the progress of the outer loop.
Combining \eqref{eq:V_inner_bd} with some basic calculations, we arrive at 
\begin{align*}    
    {\norm{Q^{(t+1)} - Q_\tau^\star}_\infty} \le \gamma{\norm{V^{(t+1)} - V_\tau^\star}_\infty} &\le \gamma{\norm{V^{(t+1)} - \breve{V}^{(t+1)}}_\infty} + \gamma\norm{\breve{V}^{(t+1)} - V_\tau^\star}_\infty\\
    &\le (1-\gamma)\epsilon + \gamma\norm{Q^{(t)} - Q_\tau^\star}_\infty.
\end{align*}
Now invoking the above relation recursively, it is ensured that
\begin{align*}
    {\norm{Q^{(t+1)} - Q_\tau^\star}_\infty} \le \epsilon + \gamma^{t+1}\norm{Q^{(0)} - Q_\tau^\star}_\infty.
\end{align*}
In view of the above relation, if one takes 
\begin{equation}\label{eq:T_main_bd}
T_{\text{\rm main}} = O\prn{\frac{1}{1-\gamma}\prn{\log \frac{1}{\epsilon} + \log\frac{1 + \tau \log|\asp|}{1-\gamma}}}
\end{equation}
iterations of the outer loop in Algorithm \ref{alg:approx_value_iter}, we have $\norm{Q^{(T_{\text{\rm main}})} - Q_\tau^\star}_\infty \le 2\epsilon$ as desired. 
\end{itemize}

Putting things together, the total iteration complexity sufficient to achieve $\epsilon$-accuracy equals to  
\begin{align*}
    &T_{\text{\rm main}}T_{\rm sub} \\
	&=O\prn{\frac{1}{\eta\tau(1-\gamma)}\prn{\log \frac{1}{\epsilon} + \log\log|\asp| + \log\frac{1}{1-\gamma} + \log\tau}\prn{\log \frac{1}{\epsilon} + \log\log|\asp| + \log \frac{1}{1-\gamma} + \log \frac{1}{\eta}}}.	
\end{align*}
Therefore the advertised iteration complexity in Theorem~\ref{Thm:MG-result} holds true by simply noticing that $\eta = \frac{1-\gamma}{2(1+\tau(\log(|\mathcal{A}|) + 1))}$ and $\tau < 1$, and hence
\begin{align*}
   \log\left(\frac{1}{\eta}\right) \leq \log\left(\frac{2\log|\mathcal{A}|+4}{1-\gamma}\right),
   \qquad \text{and}
   \qquad \log(\tau) \leq \log \left(\frac{1}{1-\gamma} \right).
\end{align*}

\subsection{Proof of Corollary~\ref{cor:unreg_MG}}

We begin by recording two supporting lemmas whose proofs are deferred to Appendix~\ref{proof:QRE_payoff_lip} and Appendix~\ref{proof:entropy_PDL}.
 
 	\begin{lemma}\label{lem:QRE_payoff_lip}
		Let $\best{\zeta} = (\best{\mu}, \best{\nu})$ be the QRE of payoff matrix $A$ and $\widetilde{\zeta}^\star_\tau = (\widetilde{\mu}^\star_\tau, \widetilde{\nu}^\star_\tau)$ be that of $\widetilde{A}$. We have
%		\begin{equation}
%			\norm{\best{\zeta} - \zeta'^\star_\tau}_1= \norm{\best{\mu} - \mu'^\star_\tau}_1 + \norm{\best{\nu} - \nu'^\star_\tau}_1 \le \frac{2}{\tau} \norm{A - A'}_\infty,
%			\label{eq:optima_shift1}
%		\end{equation}
%		\begin{equation}
%			\KL{\best{\zeta}}{\zeta'^\star_\tau} + \KL{\zeta'^\star_\tau}{\best{\zeta}}\le \frac{2}{\tau^2}\norm{A - A'}_\infty^2.
%			\label{eq:optima_shift2}        
%		\end{equation}
%		and
		\[
		\norm{\log\best{\zeta} - \log\widetilde{\zeta}_\tau^\star}_\infty \le \frac{2}{\tau}\cdot\prn{\frac{2}{\tau}\norm{A}_\infty + 1}\norm{A-\widetilde{A}}_\infty.
		\]
	\end{lemma}
	
\begin{lemma}\label{lem:entropy_PDL}
	For any single-agent MDP $(\ssp, \asp, P, r, \gamma)$ with bounded reward $0 \le r \le R$, the entropy-regularized value function satisfies
	\[
	\abs{\soft{V}^{\pi'}(\rho) - \soft{V}^{\pi}(\rho)}\le \frac{1}{1-\gamma}\left[ \frac{R}{1-\gamma} + \tau\prn{\frac{\log|\asp|}{1-\gamma} + 1 + \log(|\asp|+1)}\right] \norm{\log \pi' - \log \pi}_\infty
	\]
	for any two policies $\pi$ and $\pi'$.
\end{lemma}

We are now ready to prove Corollary~\ref{cor:unreg_MG}, which we break into a few steps.
	
\paragraph{Step 1: iteration complexity to obtain an approximate $\best{Q}$.}
It is immediate from Theorem~\ref{Thm:MG-result} that $T_{\mathsf{Q}} = \widetilde{O}\big(\frac{1}{\tau (1-\gamma)^2}\log^2(\frac{1}{\epsilon_{\mathsf{Q}}}) \big)$ iterations are sufficient to get a $\widetilde{Q}_\tau^\star \in \real^{|\mathcal{S}|\times|\mathcal{A}|\times|\mathcal{B}|}$ that achieves
	\[
	\norm{\widetilde{Q}_\tau^\star - \best{Q}}_\infty \le \epsilon_{\mathsf{Q}}.
	\]

	\paragraph{Step 2: iteration complexity to obtain an approximate $\best{\zeta}$.}
	Denote the QRE of the matrix game induced by $\widetilde{Q}_\tau^\star$ by $\widetilde{\zeta}_\tau^\star = (\widetilde{\mu}_\tau^\star, \widetilde{\nu}_\tau^\star)$. Specifically, for every $s \in \ssp$, $(\widetilde{\mu}_\tau^\star(s), \widetilde{\nu}_\tau^\star(s))$ solves the entropy-regularized matrix game induced by $\widetilde{Q}_\tau^\star(s)$. 
	Invoking PU or OMWU with $\eta = \frac{1-\gamma}{2(1+\tau(\log|\mathcal{A}| + 1 - \gamma))}$ ensures that within $T_{\mathsf{policy}} = \widetilde{O}\big(\frac{1+\tau\log|\mathcal{A}|}{(1-\gamma)\tau}\log \frac{1}{\epsilon_{\mathsf{policy}}}\big)$ iterations, we can find a policy pair $\zeta = (\mu, \nu)$ such that
	
	\[
	\norm{\log \zeta - \log  \widetilde{\zeta}_\tau^\star}_\infty \le \epsilon_{\mathsf{policy}}.
	\]
	This taken together with Lemma \ref{lem:QRE_payoff_lip} gives
	\begin{align*}
		\norm{\log \zeta - \log\best{\zeta}}_\infty &\le \norm{\log \zeta - \log  \widetilde{\zeta}_\tau^\star}_\infty + \norm{\log\best{\zeta} - \log \widetilde{\zeta}_\tau^\star}_\infty \\
		&\le \epsilon_{\mathsf{policy}} + \frac{2}{\tau}\cdot\prn{\frac{2}{\tau}\norm{\best{Q}}_\infty + 1}\norm{\best{Q} - \widetilde{Q}_\tau^\star}_\infty \\
		&\le \epsilon_{\mathsf{policy}} + {\frac{4 + 4\tau \log|\mathcal{A}| + 2\tau}{\tau^2(1-\gamma)} }\cdot\epsilon_{\mathsf{Q}}.
	\end{align*}
	Therefore, we can get $\norm{\log \zeta - \log\best{\zeta}}_\infty \le \epsilon \cdot C^{-1}$ within 
	\begin{align*}
		T_{\mathsf{Q}} + T_{\mathsf{policy}} &= \widetilde{O}\prn{\frac{1}{\tau (1-\gamma)^2}\log^2(\frac{1}{\epsilon_{\mathsf{Q}}})} + \widetilde{O}\prn{\frac{1+\tau\log|\mathcal{A}|}{\tau(1-\gamma)}\log \frac{1}{\epsilon_{\mathsf{policy}}}}\\
		&= \widetilde{O} \prn{\frac{1}{\tau (1-\gamma)^2}\log^2(\frac{1}{\epsilon})}
	\end{align*}
	 iterations as long as $C = \textsf{poly}\prn{(1-\gamma)^{-1}, \tau^{-1}, \log|\mathcal{A}|}$.
	
	\paragraph{Step 3: bounding the dual gap.}
	
	For any $\mu', \nu'$ we have
	\begin{align}
		V_\tau^{\mu', \nu}(\rho) - V_\tau^{\mu, \nu'}(\rho) &= \prn{V_\tau^{\mu', \nu}(\rho) - V_\tau^{\best\mu, \best\nu}(\rho)} + \prn{V_\tau^{\best\mu, \best\nu}(\rho) - V_\tau^{\mu, \nu'}(\rho)} \notag\\
		&\le \prn{V_\tau^{\mu', \nu}(\rho) - V_\tau^{\mu', \best\nu}(\rho)} + \prn{V_\tau^{\best\mu, \nu'}(\rho) - V_\tau^{\mu, \nu'}(\rho)},\label{eq:MDP_dual_decomp}
	\end{align}
	where in each term, only one of the policies is varied. Consequently, it is possible to invoke the well-known performance difference lemma for single-agent MDP. 
	
	We note that the same policy $\mu' $ appears in both $V_\tau^{\mu', \nu}(\rho)$ and $V_\tau^{\mu', \best\nu}(\rho)$, and it is therefore possible to invoke performance difference lemma for single-agent MDP to characterize $V_\tau^{\mu', \nu}(\rho) - V_\tau^{\mu', \best\nu}(\rho)$, we construct a MDP $(\ssp, \cB, \overline{P}, \overline{r}, \gamma)$ with
	\begin{align*}
			\overline{P}(s'|s,b) &=\sum_{a\in \cA} \mu'(a|s) P(s'|s,a,b),\\
			\overline{r}(s,b) &= \sum_{a\in \cA} \mu'(a|s) (r(s,a,b) - \tau\log\mu'(a|s)),
	\end{align*}
    and denote the associated entropy-regularized value function by $\overline{V}_\tau$.
	This allows us to write $V_\tau^{\mu', \nu}(\rho) = \overline{V}_\tau^{\nu}(\rho)$ and $V_\tau^{\mu', \best{\nu}}(\rho) = \overline{V}_\tau^{\best{\nu}}(\rho)$ (cf.~\eqref{defn:V-tau}).
	Applying Lemma \ref{lem:entropy_PDL} with $R = 1+\tau\log|\cA|$ gives 
	\begin{align*}
		\abs{V_\tau^{\mu', \nu}(\rho) - V_\tau^{\mu', \best\nu}(\rho)} &= \abs{\overline{V}_\tau^{\nu}(\rho) - \overline{V}_\tau^{\best{\nu}}(\rho)}\\
		&\le \frac{1}{1-\gamma}\left[ \frac{1+\tau\log|\cA|}{1-\gamma} + \tau\prn{\frac{\log|\cB|}{1-\gamma} + 1 + \log(|\cB|+1)}\right] \norm{\log \nu - \log \best{\nu}}_\infty.
	\end{align*}
	Similarly, one can derive 
	\begin{align*}
		\abs{V_\tau^{\mu, \nu'}(\rho) - V_\tau^{\best{\mu}, \nu'}(\rho)} &\le \frac{1}{1-\gamma}\left[ \frac{1+\tau\log|\cB|}{1-\gamma} + \tau\prn{\frac{\log|\cA|}{1-\gamma} + 1 + \log(|\cA|+1)}\right] \norm{\log \mu - \log \best{\mu}}_\infty
	\end{align*}
	for the second term in \eqref{eq:MDP_dual_decomp}.
	Plugging the above two inequalities into \eqref{eq:MDP_dual_decomp} gives 
	\begin{align*}
		V_\tau^{\mu', \nu}(\rho) - V_\tau^{\mu, \nu'}(\rho) &\le C \norm{\log \zeta - \log \best{\zeta}}_\infty 
		\leq \epsilon,
	\end{align*}
	where the last inequality holds as long as $\norm{\log \zeta - \log\best{\zeta}}_\infty \le \epsilon \cdot C^{-1}$ with 
	$$C = \frac{2}{1-\gamma}\left[ \frac{1}{1-\gamma} + \tau\prn{\frac{\log|\cA|+\log|\cB|}{1-\gamma} + 1 + \log(|\cA|+1) + \log(|\cB|+1)}\right ] .$$
	Therefore it takes $\widetilde{O} \big(\frac{1}{\tau (1-\gamma)^2}\log^2(\frac{1}{\epsilon}) \big)$ iterations to achieve
	\[
		\max_{\mu'}V_\tau^{\mu', \nu}(\rho) - \min_{\nu'}V_\tau^{\mu, \nu'}(\rho) \le \epsilon.
	\]
	% Similar to eq (10), setting $\tau = \frac{\epsilon (1-\gamma)}{4 \log|\mathcal{A}|}$ leads to $\max_{\mu'}V^{\mu', \nu}(\rho) - \min_{\nu'}V^{\mu, \nu'}(\rho) \le \epsilon$. The final iteration complexity is $\widetilde{O} (\frac{1}{ (1-\gamma)^3 \epsilon})$.

	 \section{Proof of auxiliary lemmas}

\subsection{Proof of Lemma \ref{lemma:mirror_update}} 
\label{Sec:Pf-key-lemma}
Lemma~\ref{lemma:mirror_update} follows directly from the update sequence~\eqref{eq:mirror_update} and the form of the optimal solution pair $(\best{\mu}, \best{\nu})$, provided in \eqref{eq:QRE-matrix}.
Given the update sequence~\eqref{eq:mirror_update}, taking logarithm of both sides of the first equation gives 
\[
    \log \mu^{(t+1)} = (1-\eta\tau)\log \mu^{(t)} +\eta A z_2 + c\cdot\mathbf{1},
\]
where $c$ is the corresponding normalization constant. By rearranging terms and taking the inner product with $z_1 - \best{\mu}$, we have
\begin{equation}
    \innprod{\log \mu^{(t+1)} - (1-\eta\tau)\log \mu^{(t)}, z_1 - \best{\mu}} = \eta z_1^\top A z_2 - {\eta \best{\mu}}^\top A z_2,
    \label{eq:step1_1} 
\end{equation}
Similarly, one can derive 
\begin{equation}
    \innprod{\log \nu^{(t+1)} - (1-\eta\tau)\log \nu^{(t)}, z_2 - \best{\nu}} = - \eta z_1^\top A z_2 + \eta z_1^\top A\best{\nu}.
    \label{eq:step1_2} 
\end{equation}
By summing up equations~\eqref{eq:step1_1} and \eqref{eq:step1_2}, it is guarantee that 
\begin{equation}
    \innprod{\log \zeta^{(t+1)} - (1-\eta\tau)\log \zeta^{(t)}, \zeta_z - \best{\zeta}} = - \eta {\best{\mu}}^\top Az_2 + \eta z_1^\top A\best{\nu},
    \label{eq:step1} 
\end{equation}
where $\zeta(z)=(z_1, z_2)$.

On the other hand, recall the optimal policy pair $(\best{\mu}, \best{\nu})$ satisfies the following fixed point equation 
\[
    \begin{cases}
        \best{\mu}(a) \propto \exp( [A \best{\nu}]_a/\tau), \quad \forall a \in \mathcal{A}, \\
        \best{\nu}(b) \propto \exp( -[A^\top \best{\mu}]_b/\tau), \quad \forall b \in \mathcal{B}.\\      
    \end{cases} 
\]
Taking logarithm of both sides of the first relation gives 
\begin{equation}
    \eta\tau\log \best{\mu} =  \eta A \best{\nu} + c\cdot\mathbf{1},
    \label{eq:log_st_1}
\end{equation}
for some normalization constant $c$. Again, by taking the inner product with $z_1 - \best{\mu}$, we have
\begin{equation}
    \innprod{\eta\tau\log \best{\mu}, z_1 - \best{\mu}} 
    = \eta (z_1-\best{\mu})^{\top} A \best{\nu},
    % \innprod{z_1}^\top A\best{\nu} - {\eta \best{\mu}}^\top A\best{\nu}.
    \label{eq:step2_1} 
\end{equation}
and similarly
\begin{equation}
    \innprod{\eta\tau\log \best{\nu}, z_2 - \best{\nu}} 
    = \eta {\best{\mu}}^{\top} A (z_2 - \best{\nu} ). 
    % - \eta {\best{\mu}}^\top A z_2 + \eta {\best{\mu}}^\top A\best{\nu}.
    \label{eq:step2_2} 
\end{equation}

Combining inequalities \eqref{eq:step1_1} and \eqref{eq:step2_1}, we arrive at inequality~\eqref{eq:ohhh_mu}; combining inequalities \eqref{eq:step1_2} and \eqref{eq:step2_2} gives inequality~\eqref{eq:ohhh_nu}.
Moreover, putting together inequalities \eqref{eq:step1}, \eqref{eq:step2_1} and \eqref{eq:step2_2} leads to
\[
    \innprod{\log \zeta^{(t+1)} - (1-\eta\tau)\log \zeta^{(t)} - \eta\tau \log \best{\zeta}, \zeta(z) - \best{\zeta}} = 0.
\]

    % Therefore, we have
    % \[
    %     \innprod{\best{\mu} - \mu, \best{\nu}-\nu}_A - \tau \KL{\mu}{\best{\mu}} \le f_{\tau}(\mu, \nu) - f_{\tau}(\best{\mu}, \best{\nu}) \le \innprod{\best{\mu} - \mu, \best{\nu}-\nu}_A + \tau \KL{\nu}{\best{\nu}}.
    % \]
% \end{proof}

\subsection{Proof of Lemma \ref{lemma:gap_KL_diff_decomp}} \label{proof_lemma:gap_KL_diff_decomp}
 
We begin with establishing \eqref{eq:gap_KL}.
By the definition of $f_{\tau}(\mu, \nu)$, direct calculations yield 
    \begin{align}
        f_{\tau}(\best{\mu}, \best{\nu}) - f_{\tau}(\mu, \best{\nu}) &= (\best{\mu} - \mu)^\top A \best{\nu} + \tau\mu^\top \log \mu - \tau {\best{\mu}}^\top \log \best{\mu} \nonumber \\
        &= \tau\prn{\innprod{\best{\mu} - \mu, \log \best{\mu}} + \mu^\top \log \mu -  {\best{\mu}}^\top \log \best{\mu}}= \tau \KL{\mu}{\best{\mu}}. \label{eq:baozi}
    \end{align}
    Here, the second equality is obtained by plugging in \eqref{eq:log_st_1}. Similarly, we have
    \begin{align}\label{eq:jiaozi}
        f_{\tau}(\best{\mu}, \nu) - f_{\tau}(\best{\mu}, \best{\nu}) = \tau \KL{\nu}{\best{\nu}}. 
    \end{align}
    Summing these two equalities completes the proof of \eqref{eq:gap_KL}.

Turning to \eqref{eq:f_diff_decomp},  we first write 
    \begin{align*}
        f_{\tau}(\mu, \nu) + f_{\tau}(\best{\mu}, \best{\nu}) &= \mu^\top A \nu + {\best{\mu}}^\top A \best{\nu} + \tau \mathcal{H}(\mu) - \tau \mathcal{H}(\nu) + \tau \mathcal{H}(\best{\mu}) - \tau \mathcal{H}(\best{\nu}),\\
        f_{\tau}(\best{\mu}, \nu) + f_{\tau}(\mu, \best{\nu}) &= {\best{\mu}}^\top A \nu + {\mu}^\top A \best{\nu} + \tau \mathcal{H}(\best{\mu}) - \tau \mathcal{H}(\nu) + \tau \mathcal{H}(\mu) - \tau \mathcal{H}(\best{\nu}).
    \end{align*}
    As a consequence, taking the difference of the above two equations leads to
    \[
        f_{\tau}(\mu, \nu) + f_{\tau}(\best{\mu}, \best{\nu})- f_{\tau}(\best{\mu}, \nu) - f_{\tau}(\mu, \best{\nu}) = ( \best{\mu} - \mu )^{\top} A (\best{\nu}-\nu ).
    \]
    This in turn allows us to write $f_{\tau}(\mu, \nu) - f_{\tau}(\best{\mu}, \best{\nu})$ as follows
    \begin{equation}
        f_{\tau}(\mu, \nu) - f_{\tau}(\best{\mu}, \best{\nu})
        =  ( \best{\mu} - \mu)^{\top}A ( \best{\nu}-\nu ) + f_{\tau}(\best{\mu}, \nu) + f_{\tau}(\mu, \best{\nu})- 2f_{\tau}(\best{\mu}, \best{\nu}).
        \label{eq:mantou}
    \end{equation}
    % Direct calculations yield 
    % \begin{align}
    %     f_{\tau}(\best{\mu}, \best{\nu}) - f_{\tau}(\mu, \best{\nu}) &= (\best{\mu} - \mu)^\top A \best{\nu} + \tau\mu^\top \log \mu - \tau {\best{\mu}}^\top \log \best{\mu} \nonumber \\
    %     &= \tau\brk{\innprod{\best{\mu} - \mu, \log \best{\mu}} + \mu^\top \log \mu -  {\best{\mu}}^\top \log \best{\mu}} \nonumber\\
    %     &= \tau \KL{\mu}{\best{\mu}}. \label{eq:baozi}
    % \end{align}
    % Here, the second equality is obtained by plugging in \eqref{eq:log_st_1}. Similarly, we have
    % \begin{align}\label{eq:jiaozi}
    %     f_{\tau}(\best{\mu}, \nu) - f_{\tau}(\best{\mu}, \best{\nu}) = \tau \KL{\nu}{\best{\nu}}. 
    % \end{align}
   Finally, plugging \eqref{eq:baozi} and \eqref{eq:jiaozi} into \eqref{eq:mantou} reveals the desired relation \eqref{eq:f_diff_decomp}.
%    \begin{align*}
%        f_{\tau}(\mu, \nu) - f_{\tau}(\best{\mu}, \best{\nu})
%        % &=  ( \best{\mu} - \mu)^{\top}A ( \best{\nu}-\nu ) + f_{\tau}(\best{\mu}, \nu) + f_{\tau}(\mu, \best{\nu})- 2f_{\tau}(\best{\mu}, \best{\nu})\\
%        &= ( \best{\mu} - \mu)^{\top} A( \best{\nu}-\nu )  + \tau \KL{\nu}{\best{\nu}} - \tau \KL{\mu}{\best{\mu}}.
%    \end{align*}

\subsection{Proof of Lemma~\ref{lem:log_pi_gap}} \label{proof_lem:log_pi_gap}
The second inequality follows directly from \cite[Lemma 27]{mei2020global}. The first inequality has appeared, e.g., in \cite{cen2020fast}. We reproduce a short proof for self-completeness. By straightforward calculations, the gradient of the function $\log (\norm{\exp(x)}_1)$ is given by
    \[
        \nabla_x \log (\norm{\exp(x)}_1) = \exp(x)/\norm{\exp(x)}_1,
    \]
    which implies $\norm{\nabla_x \log (\norm{\exp(x)}_1)}_1 = 1$, $\forall x \in \mathbb{R}^{|\mathcal{A}|}$. Therefore, we have
    \begin{align*}
        \norm{\log \mu_1 - \log \mu_2}_\infty &= \norm{x_1 - x_2 - \log (\norm{\exp(x_1)}_1)\cdot \mathbf{1} + \log (\norm{\exp(x_2)}_1)\cdot \mathbf{1}}_\infty\\
        &\le \norm{x_1 - x_2}_\infty + \Big|- \log (\norm{\exp(x_1)}_1) + \log (\norm{\exp(x_2)}_1)\Big|\\
        &= \norm{x_1 - x_2}_\infty + \Big|\left\langle x_1 - x_2, \nabla_x \log (\norm{\exp(x)}_1)|_{x = x_c} \right\rangle\Big|\\
        &\le \norm{x_1 - x_2}_\infty + \Big| \norm{x_1 - x_2}_\infty \norm{\nabla_x \log (\norm{\exp(x)}_1)|_{x = x_c}}_1 \Big|\\
        &=2\norm{x_1 - x_2}_\infty,
    \end{align*}
    where $x_c$ is a certain convex combination of $x_1$ and $x_2$.

\subsection{Proof of Lemma \ref{lemma:dual_gap}}
\label{sec:pf_dual_gap}
 Since
 $$  \max_{\mu'\in \Delta(\mathcal{A})} f_\tau(\mu', \nu) - \min_{\nu'\in\Delta(\mathcal{B})}f_\tau(\mu, \nu') = \max_{\mu'\in \Delta(\mathcal{A}),\nu'\in\Delta(\mathcal{B})} f_\tau(\mu', \nu) -  f_\tau(\mu, \nu') , $$
 it boils down to control $f_\tau(\mu', \nu) -  f_\tau(\mu, \nu')$ for any $(\mu', \nu') \in \Delta(\cA)\times \Delta(\cB)$. Towards this, we have
\begin{align}
    f_\tau(\mu', \nu) - f_\tau(\mu, \nu') &= \prn{f_\tau(\mu', \nu) - f_{\tau}(\mu', \best{\nu}) - f_\tau(\mu, \nu') + f_{\tau}(\best{\mu}, \nu')}- \prn{f_{\tau}(\best{\mu}, \nu') - f_{\tau}(\mu', \best{\nu})}\nonumber\\
    &= \prn{f_\tau(\mu', \nu) - f_{\tau}(\mu', \best{\nu}) - f_\tau(\mu, \nu') + f_{\tau}(\best{\mu}, \nu')}- \tau \KL{\zeta'}{\best{\zeta}},
    \label{eq:gap_KL_decomp_step1}
\end{align}
where the last step is due to $f_{\tau}(\mu, \best{\nu}) - f_{\tau}(\best{\mu}, \nu) = \tau \KL{\zeta}{\best{\zeta}}$, as revealed in Lemma~\ref{lemma:gap_KL_diff_decomp} (cf.~\eqref{eq:gap_KL}).

To continue, observe that
\begin{align*}
    f_\tau(\mu', \nu) - f_{\tau}(\mu', \best{\nu}) &= {\mu'}{}^\top A(\nu - \best{\nu}) + \nu^\top \log \nu - {\best{\nu}}^\top \log \best{\nu}\\
    &= (\mu' - \best{\mu})^\top A(\nu - \best{\nu}) + f_\tau(\best{\mu}, \nu) - f_{\tau}(\best{\mu}, \best{\nu}).
\end{align*}
Similarly, we have
\[
    - f_\tau(\mu, \nu') + f_{\tau}(\best{\mu}, \nu') = -\prn{\mu - \best{\mu}}^\top A(\nu' - \best{\nu}) + f_{\tau}(\best{\mu}, \best{\nu}) - f_\tau({\mu}, \best{\nu}).
\]
Plugging the above two equalities into \eqref{eq:gap_KL_decomp_step1} gives 
\begin{align*}
    f_\tau(\mu', \nu) - f_{\tau}(\mu, \nu') &= \prn{\mu' - \best{\mu}}^\top A(\nu - \best{\nu})-\prn{\mu - \best{\mu}}^\top A(\nu' - \best{\nu}) + f_\tau(\best{\mu}, \nu) - f_\tau({\mu}, \best{\nu}) - \tau \KL{\zeta'}{\best{\zeta}}\\
    & = \prn{\mu' - \best{\mu}}^\top A(\nu - \best{\nu})-\prn{\mu - \best{\mu}}^\top A(\nu' - \best{\nu}) + \tau \KL{\zeta}{\best{\zeta}} - \tau \KL{\zeta'}{\best{\zeta}}\\ 
    & \leq \norm{A}_{\infty}\prn { \norm{\mu' - \best{\mu}}_1 \norm{\nu - \best{\nu}}_1 +  \norm{\nu' - \best{\nu}}_1 \norm{\mu - \best{\mu}}_1 } + \tau \KL{\zeta}{\best{\zeta}} - \tau \KL{\zeta'}{\best{\zeta}} \\
    &\overset{\mathrm{(i)}}{\le} \frac{1}{2}\norm{A}_\infty \brk{\frac{\tau}{\norm{A}_\infty}\prn{\norm{\mu' - \best{\mu}}_1^2+ \norm{\nu' - \best{\nu}}_1^2} + \frac{\norm{A}_\infty}{\tau}\prn{\norm{\mu - \best{\mu}}_1^2+ \norm{\nu - \best{\nu}}_1^2}}\\
    &\qquad + \tau \KL{\zeta}{\best{\zeta}} - \tau \KL{\zeta'}{\best{\zeta}}\\
    &\overset{\mathrm{(ii)}}{\le} \tau \KL{\zeta'}{\best{\zeta}} + \frac{\norm{A}_\infty^2}{\tau} \KL{\best{\zeta}}{\zeta}+ \tau \KL{\zeta}{\best{\zeta}} - \tau \KL{\zeta'}{\best{\zeta}}\\
    &=\frac{\norm{A}_\infty^2}{\tau} \KL{\best{\zeta}}{\zeta}+ \tau \KL{\zeta}{\best{\zeta}},
\end{align*}
where the second step invokes Lemma~\ref{lemma:gap_KL_diff_decomp} (cf.~\eqref{eq:gap_KL}), (i) follows from Young's inequality, namely $ab\leq \frac{a^2}{2\varepsilon}+\frac{\varepsilon b^2}{2}$ with $\varepsilon =  \frac{\norm{A}_\infty}{\tau}$, and (ii)
 results from Pinsker's inequality.\
Taking maximum over $\mu', \nu'$ finishes the proof.

\subsection{Proof of Lemma~\ref{lemma:OMWU_grad_bound}}
\label{proof:lemma_OMWU_grad_bound}
 
    First, we show that the update of $\bar\nu^{(t)}$ (cf. \eqref{eq:single_omwu_bar_nu}) satisfies
    \begin{equation}\label{eq:exp_nu_induction}
    \bar\nu^{(t)}(b) \propto \exp([A\widetilde{\mu}^{(t)}]_b/\tau) \quad \mbox{for some } \widetilde{\mu}^{(t)} \in \Delta(\mathcal{A})
    \end{equation}
    by induction.
    \begin{itemize}
    \item For $t = 0$, it is easily seen that $\bar\nu^{(0)}(b) = \frac{1}{|\mathcal{B}|} \propto \exp([A^\top\widetilde{\mu}^{(0)}]_b/\tau)$ with $\widetilde{\mu}^{(0)} = 0$. 
    \item Now assume \eqref{eq:exp_nu_induction} holds for all steps up to $t$.
    The update rule (cf. \eqref{eq:single_omwu_bar_nu}) implies
    \begin{align*}
        \bar\nu^{(t+1)}(b) &\propto \bar\nu^{(t)}(b)^{1-\eta_t\tau}\exp(\eta_t[A^\top \bar\mu^{(t)}])\\
        &\propto \exp((1-\eta_t\tau)[A^\top\widetilde{\mu}^{(t)}]_b/\tau)\exp(\eta_t[A^\top \bar\mu^{(t)}]) = \exp([A^\top \widetilde{\mu}^{(t+1)}]_b/\tau),
    \end{align*}
    with $\widetilde{\mu}^{(t+1)} = (1-\eta_t \tau)\widetilde{\mu}^{(t)} + \eta_t \tau   \bar\mu^{(t)} \in \Delta(\mathcal{B})$. 
    \end{itemize}
    Therefore, the claim \eqref{eq:exp_nu_induction} holds for all $t > 0$.

    It then follows from \eqref{eq:exp_nu_induction} straightforwardly that
    \[
        \frac{\bar\nu^{(t)}(b_1)}{\bar\nu^{(t)}(b_2) }  = \frac{\exp([A\widetilde{\mu}^{(t)}]_{b_1}/\tau)}{\exp([A\widetilde{\mu}^{(t)}]_{b_2}/\tau)}   \leq \exp(2\norm{A}_\infty/\tau)
    \]
    for any $b_1, b_2 \in \mathcal{B}$. Therefore, we have
    \[
        1 = \sum_{b \in \mathcal{B}} \bar\nu^{(t)}(b) \le |\mathcal{B}| \exp(2\norm{A}_\infty/\tau) \cdot\min_{b\in\mathcal{B}} \bar\nu^{(t)}(b) ,
    \]
    which gives $\min_{b\in\mathcal{B}}\bar\nu^{(t)}(b) \ge |\mathcal{B}|^{-1}\exp(-2\norm{A}_\infty/\tau)$, or equivalently
    \[
        \norm{\log\bar\nu^{(t)}}_\infty \le 2\norm{A}_\infty/\tau + \log |\mathcal{B}|.
    \]
We conclude the proof in view of the expression of $\overline\nabla^{(t)}$ in \eqref{eq:nabla_t}:
    \[
        \norm{\overline\nabla^{(t)}}_\infty = \norm{A^\top \bar\mu^{(t)} + \tau\log \bar\nu^{(t)}}_\infty \le \norm{A^\top \bar\mu^{(t)}}_\infty + \tau \norm{\log\nu^{(t)}}_\infty \le \tau\log |\mathcal{B}| + 3\norm{A}_\infty.
    \]

\subsection{Proof of Lemma~\ref{lem:QRE_payoff_lip}}
	\label{proof:QRE_payoff_lip}

Instantiating \eqref{eq:step2_1} at $z_1 := \mu'^\star_\tau $, we have
	\begin{equation*}
			\innprod{\tau\log \best{\mu},  \mu'^\star_\tau - \best{\mu}} = { \mu'^\star_\tau}^\top A\best{\nu} - { \best{\mu}}^\top A\best{\nu}.
		\end{equation*}
Similarly, instantiating \eqref{eq:step2_2} with $z_2 :=\nu'^\star_\tau$, we obtain 
		\begin{equation*}
			\innprod{\tau\log \best{\nu}, \nu'^\star_\tau - \best{\nu}} = - {\best{\mu}}^\top A \nu'^\star_\tau + {\best{\mu}}^\top A\best{\nu}.
		\end{equation*}
		Summing the above two equalities then leads to
		\[
		\tau \innprod{\log \best{\zeta}, \zeta'^\star_\tau-\best{\zeta}} = - {\best{\mu}}^\top A \nu'^\star_\tau +{\mu'^\star_\tau}^\top A\best{\nu}.
		\]
In view of symmetry,  the following equality holds as well
		\[
		\tau \innprod{\log \zeta'^\star_\tau, \best{\zeta}-\zeta'^\star_\tau} = - {\mu'^\star_\tau}^\top A' \best{\nu} +{\best{\mu}}^\top A'\nu'^\star_\tau,
		\]

Taken the above relations collectively allows us to arrive at 
		\begin{align}
			\KL{\best{\zeta}}{\zeta'^\star_\tau} + \KL{\zeta'^\star_\tau}{\best{\zeta}} & = \innprod{\log \zeta'^\star_\tau- \log \best{\zeta}, \zeta'^\star_\tau-\best{\zeta}} \nonumber\\
			& =  \frac{1}{\tau}\prn{{\best{\mu}}^\top (A -A') \nu'^\star_\tau - {\mu'^{\star}}^\top (A -A') \best{\nu} }\nonumber\\
			& = \frac{1}{\tau}\prn{{\best{\mu}}^\top (A -A')(\nu'^\star_\tau - \best{\nu}) - \prn{\mu'^\star_\tau - \best{\mu}}^\top (A -A') \best{\nu}}\nonumber\\
			& \le \frac{1}{\tau} \norm{A-A'}_\infty(\norm{\best{\mu} - \mu'^\star_\tau}_1 + \norm{\best{\nu} - \nu'^\star_\tau}_1).
			\label{eq:optima_shift_step1}
		\end{align}
		On the other hand, we have
		\begin{equation}
			\frac{1}{2}(\norm{\best{\mu} - \mu'^\star_\tau}_1 + \norm{\best{\nu} - \nu'^\star_\tau}_1)^2\le {\norm{\best{\mu} - \mu'^\star_\tau}_1^2 + \norm{\best{\nu} - \nu'^\star_\tau}_1^2} \le \KL{\best{\zeta}}{\zeta'^\star_\tau} + \KL{\zeta'^\star_\tau}{\best{\zeta}},
			\label{eq:optima_shift_step2}
		\end{equation}
		where the second step results from Pinsker's inequality. Combining \eqref{eq:optima_shift_step1} and \eqref{eq:optima_shift_step2} leads to 
	\begin{equation}
			\norm{\best{\zeta} - \zeta'^\star_\tau}_1= \norm{\best{\mu} - \mu'^\star_\tau}_1 + \norm{\best{\nu} - \nu'^\star_\tau}_1 \le \frac{2}{\tau} \norm{A - A'}_\infty.
			\label{eq:optima_shift1}
		\end{equation}
%Plugging the result into \eqref{eq:optima_shift_step1} leads to \eqref{eq:optima_shift2}.
This in turn allows us to bound
		\begin{align*}
			\norm{A\best{\nu} - A'\nu'^\star}_\infty & \le \norm{A(\best{\nu} - \nu'^\star)}_\infty + \norm{(A - A')\nu'^\star}_\infty\\
			&\le \norm{A}_\infty\norm{\best{\nu} - \nu'^\star}_1 + \norm{A-A'}_\infty\\
			&\le \norm{A}_\infty\norm{\best{\zeta} - \zeta'^\star_\tau}_1 + \norm{A-A'}_\infty\\
			&\le \prn{\frac{2}{\tau}\norm{A}_\infty + 1}\norm{A-A'}_\infty,
		\end{align*}
		where the final step invokes \eqref{eq:optima_shift1}. Since $\best{\mu} \propto \exp (-A\best{\nu}/\tau)$, $\mu'^\star \propto \exp (-A\nu'^\star/\tau)$, we invoke Lemma \ref{lem:log_pi_gap} to arrive at 
		\[
		\norm{\log\best{\mu} - \log\mu'^\star}_\infty \le \frac{2}{\tau}\norm{A\best{\nu} - A'\nu'^\star}_\infty \le \frac{2}{\tau}\cdot\prn{\frac{2}{\tau}\norm{A}_\infty + 1}\norm{A-A'}_\infty,
		\]
        which establishes the desired bound.

 \subsection{Proof of Lemma~\ref{lem:entropy_PDL}}
 \label{proof:entropy_PDL}
 
Using the regularized version of the performance difference lemma (see \cite[Lemma 7]{zhan2021policy} or \cite[Lemma 2]{lan2021policy}), we have:
	\begin{align}
\soft{V}^{\pi'}(\rho) - \soft{V}^{\pi}(\rho) 
		&=\frac{1}{1-\gamma} \ex{s \sim d_\rho^{\pi'}}{\innprod{Q^{\pi'}_\tau(s), \pi'(\cdot|s) - \pi(\cdot|s)} -\tau\brk{\mathcal{H}(\pi'(\cdot|s)) - \mathcal{H}(\pi(\cdot|s))}}.\label{eq:entropy_PDL}
	\end{align}
We then bound the two terms separately.
\begin{itemize}
\item To control the first term, we notice that $\soft{Q}(s, a)$ is bounded by $0 \le \soft{Q}(s,a) \le \frac{R}{1-\disct} + \frac{\tau\log|\asp|}{1-\disct}$ for all $(s,a)$. Hence, we have
	\begin{align}
\abs{\innprod{Q^{\pi'}_\tau(s), \pi'(\cdot|s) - \pi(\cdot|s )}} 		&\le \norm{Q^{\pi'}_\tau(s)}_\infty \norm{\pi'(\cdot|s) - \pi(\cdot|s)}_1 \notag\\
		&\le \norm{Q^{\pi'}_\tau(s)}_\infty\norm{\log \pi'(\cdot|s) - \log \pi(\cdot|s)}_\infty\notag\\
		&\le \frac{1}{1-\gamma}\prn{R+\tau\log|\asp|}\norm{\log \pi'(\cdot|s) - \log \pi(\cdot|s)}_\infty,
		\label{eq:entropy_PDL_first}
	\end{align}
	where the second step is due to \cite[Lemma 24]{mei2020global}. 
\item Turning to the entropy difference term in \eqref{eq:entropy_PDL}, let  
	\[
		F(\log\pi(\cdot|s)) = -\innprod{\exp(\log \pi(\cdot|s)), \log \pi(\cdot|s)} = \mathcal{H}(\pi(\cdot|s)).
	\]
	We can then invoke mean value theorem to show
	\begin{align}
\abs{-\mathcal{H}(\pi'(\cdot|s)) + \mathcal{H}(\pi(\cdot|s))} 
		&=\abs{-F(\log\pi'(\cdot|s)) + F(\log\pi(\cdot|s))}\notag\\
		&=\abs{\innprod{\log \pi'(\cdot|s) - \log \pi(\cdot|s), \nabla F(\log \xi)}} \notag\\
		&\le \norm{ \log\pi'(\cdot|s)- \log\pi(\cdot|s)}_\infty \norm{\nabla F(\log \xi)}_1 \notag\\
		&\leq  \norm{ \log\pi'(\cdot|s)- \log\pi(\cdot|s)}_\infty ( \norm{\xi}_1  + \norm{ \xi\odot\log \xi }_1) ,
		% &\le \norm{ \log\pi'(\cdot|s)- \log\pi(\cdot|s)}_\infty \norm{\xi + \xi\log \xi}_1\\	
		% &\le \norm{ \log\pi'(\cdot|s)- \log\pi(\cdot|s)}_\infty \prn{ 1 + \log(|\asp|+1)},
		\label{eq:entropy_log_mean_value}
	\end{align}
	where $\log \xi = c \log \pi(\cdot|s) + (1-c)\log \pi'(\cdot|s)$ for some constant $0<c<1$, and the last line follows from  $\nabla F(\log \xi) = -\xi - \xi\odot\log \xi$, where $\odot$ denotes point-wise multiplication. H\"{o}lder's inequality guarantees that
	\[
	\norm{\xi}_1 = \abs{\innprod{\pi(\cdot|s)^c, \pi'(\cdot|s)^{1-c}}} \le \norm{\pi(\cdot|s)^c}_{1/c} \norm{\pi'(\cdot|s)^{1-c}}_{1/(1-c)} = 1.
	\]
	Introduce $\overline{\xi}$ which appends a scalar $1-\norm{\xi}_1$ to $\xi$ as $\overline{\xi} = [
		\xi^{\top} ,
		1-\norm{\xi}_1 ]^{\top}$, so that $\overline{\xi}$ is a probability vector. It is straightforward to get 
	\begin{align*}
		\norm{\xi\odot\log \xi}_1 &= -\sum_{a\in \cA} \xi(a)\log\xi(a) \\
		&\le -\sum_{a\in \cA} \xi(a)\log\xi(a) - (1-\norm{\xi}_1)\log(1-\norm{\xi}_1) 
		= \mathcal{H}(\overline{\xi}) \le \log(|\cA| + 1).
	\end{align*}
	Substitution of the above two inequalities into \eqref{eq:entropy_log_mean_value} gives
	\begin{equation} \label{eq:entropy_smooth}
		\abs{-\mathcal{H}(\pi'(\cdot|s)) + \mathcal{H}(\pi(\cdot|s))} \le \norm{ \log\pi'(\cdot|s)- \log\pi(\cdot|s)}_\infty (1 + \log(|\cA| + 1)).
	\end{equation}
\end{itemize}
	Plugging \eqref{eq:entropy_smooth} and \eqref{eq:entropy_PDL_first} into \eqref{eq:entropy_PDL} completes the proof.

\end{document}